\documentclass[11pt]{article}
\usepackage{tipa}
\usepackage{latexsym,amssymb,amsfonts}
\usepackage{mathrsfs}
\usepackage{graphicx,color}
\usepackage[pdftex,colorlinks=true,bookmarks=false,citecolor=blue]{hyperref}
\usepackage{psfrag}
\usepackage{amsmath, amsbsy}
\usepackage{amsopn, amstext}
\usepackage{amsmath,multicol,amsopn}
\usepackage{latexsym,amssymb}
\usepackage{amsbsy, amstext,makeidx}

\def\sqr#1#2{{\vcenter{\vbox{\hrule height.#2pt
              \hbox{\vrule width.#2pt height#1pt \kern#1pt \vrule width.#2pt}
              \hrule height.#2pt}}}}
\def\signed #1{{\unskip\nobreak\hfil\penalty50
              \hskip2em\hbox{}\nobreak\hfil#1
              \parfillskip=0pt \finalhyphendemerits=0 \par}}
\def\endpf{\signed {$\sqr69$}}

\def\dbR{{\mathop{\rm l\negthinspace R}}}

\def\3n{\negthinspace \negthinspace \negthinspace }
\def\1n{\negthinspace }

\def\dbN{{\mathop{\rm l\negthinspace N}}}

\def\dbR{\mathbb{R}}
\def\deq {\buildrel \triangle \over =}

\def\ds{\displaystyle}

\def\ns{\noalign{\ss}}

\def\rank{{\rm rank}}
\def\a{\alpha}

\def\d{\delta}
\def\e{\varepsilon}

\def\si{\sigma}

\def\f{\varphi}

\def\o{\omega}

\def\bd{{\mathbf d}}
\def\G{\Gamma}
\def\D{\Delta}

\def\O{\Omega}

\def\cC{{\cal C}}

\def\cG{{\cal G}}

\def\cM{{\cal M}}
\def\cN{{\cal N}}

\def\cP{{\cal P}}

\def\cT{{\cal T}}

\def\no{\noindent}

\def\ss{\smallskip}
\def\ms{\medskip}
\def\bs{\bigskip}
\def\q{\quad}
\def\qq{\qquad}

\def\liminf{\mathop{\underline{\rm lim}}}

\def\pa{\partial}

\def\wt{\widetilde}
\def\cd{\cdot}
\def\cds{\cdots}

\def\span{\hbox{\rm span$\,$}}

\def\co{\mathop{{\rm co}}}
\def\coh{\mathop{\overline{\rm co}}}
\def\deq{\mathop{\buildrel\D\over=}}

\def\Im{{\mathop{\rm Im}\,}}

\def\|{\Big |}
\def\({\Big (}
\def\){\Big )}
\def\[{\Big[}
\def\]{\Big]}
\def\bde{\begin{definition}}
\def\ede{\end{definition}}
\def\be{\begin{equation}}
\def\bel{\begin{equation}\label}
\def\ee{\end{equation}}
\def\bt{\begin{theorem}}
\def\et{\end{theorem}}
\def\bc{\begin{corollary}}
\def\ec{\end{corollary}}
\def\bl{\begin{lemma}}
\def\el{\end{lemma}}
\def\bp{\begin{proposition}}
\def\ep{\end{proposition}}
\def\bas{\begin{assumption}}
\def\eas{\end{assumption}}
\def\br{\begin{remark}}
\def\er{\end{remark}}
\def\ba{\begin{array}}
\def\ea{\end{array}}
\def\ed{\end{document}}

\def\square#1{\vbox{\hrule\hbox{\vrule height#1%
     \kern#1\vrule}\hrule}}
\def\rectangle#1#2{\vbox{\hrule\hbox{\vrule height#1%
     \kern#2\vrule}\hrule}}

\font\tenbb=msbm10 \font\sevenbb=msbm7 \font\fivebb=msbm5

\newfam\bbfam
\scriptscriptfont\bbfam=\fivebb \textfont\bbfam=\tenbb
\scriptfont\bbfam=\sevenbb

\newtheorem{lemma}{Lemma}[section]
\newtheorem{remark}{Remark}[section]
\newtheorem{example}{Example}[section]
\newtheorem{theorem}{Theorem}[section]
\newtheorem{corollary}{Corollary}[section]

\newtheorem{definition}{Definition}[section]
\newtheorem{proposition}{Proposition}[section]
\newtheorem{assumption}{Assumption}[section]

\makeatletter
   
   \@addtoreset{equation}{section}
\makeatother

\begin{document}

\title{\bf  Verifiable constraint qualifications for infinite-dimensional optimization problems and their applications to optimal control problems}

\author{Xu Liu\thanks{School of Mathematics and Statistics, Northeast Normal
University, Changchun 130024, China.  The
research of this author is partially supported
by National Key R\&D Program of China under grant 2023YFA1009002 and
NSF of China under grant 12371444. {\small\it E-mail:} {\small\tt
liux216@nenu.edu.cn}.},~~~~Qi
L\"{u}\thanks{School of Mathematics, Sichuan
University, Chengdu 610064, China. The research
of this author is partially supported by NSF of
China under grants 12025105  and 11931011, and by the Chang
Jiang Scholars Program from the Ministry of
Education of China.
{\small\it E-mail:} {\small\tt
lu@scu.edu.cn}.},~~~~Haisen Zhang\thanks{School
of Mathematical Sciences, Sichuan Normal
University, Chengdu 610068, China.  The research
of this author is partially supported by NSF of
China under grant  12071324 and the Key Project of Sichuan Science and Technology Education Joint Fund under grant  2024NSFSC1963. {\small\it E-mail:} {\small\tt
haisenzhang@yeah.net}.}~~~~and~~~~Xu Zhang\thanks{School of Mathematics and New Cornerstone Science Laboratory, Sichuan
University, Chengdu 610064, China. The research
of this author is partially supported by  the National Key R\&D Program of China
under grant  2024YFA1013400  and the
New Cornerstone Investigator Program.
{\small\it E-mail:} {\small\tt
zhang$\_$xu@scu.edu.cn}.}}

\date{}

\maketitle

\begin{abstract}
This paper employs a finite codimensionality condition  to establish an enhanced Fritz John
 condition for general constrained nonlinear  infinite-dimensional  optimization problems.
 In applications, this approach provides a new, unified framework for deriving first-order necessary conditions across a broad class of optimal control problems involving both deterministic and stochastic systems. Compared to  existing constraint qualifications,
 our finite codimensionality condition, which is equivalent to the validity of certain \textit{a priori} estimates, yields a direct and analytically tractable verification method for each application. Furthermore, the core ideas of this method can be further extended to investigate the KKT conditions for infinite-dimensional optimization problems.
\end{abstract}

\bs

\no{\bf Key Words}.   Infinite-dimensional optimization
problem, Fritz John condition,  finite codimensionality, {\it a priori} estimate.
\ms

\no{\bf AMS subject classifications}.
49K27, 93C25,  35Q93, 93E20.

\section{Introduction}

Infinite-dimensional optimization is not merely an abstract extension of finite-dimensional optimization; it is the natural language for a wide range of fundamental problems, including optimal control, calculus of variations, partial differential equation (PDE)-constrained optimization, stochastic control, and inverse problems, where the optimization variable is typically a function, measure, or operator rather than a finite-dimensional vector. This perspective is especially indispensable for machine learning, because the canonical learning problem is to minimize the expected risk over a hypothesis space of functions, which is typically an infinite-dimensional function space, such as a  reproducing kernel Hilbert space (RKHS), a Sobolev or Besov space, or a space of probability measures. Although practical models are parameterized by finitely many weights, the parameterization is merely a computational representation: the true object being optimized is a function, and many different parameter vectors can induce the same function, as is strikingly evident in overparameterized neural networks. The infinite-dimensional viewpoint therefore reveals the intrinsic geometry of learning, provides a unified framework for kernel methods, mean-field neural networks, and gradient flows in Hilbert or Wasserstein spaces, and yields necessary conditions, regularity, and generalization insights that remain obscured or are easily misattributed when training is viewed only as finite-dimensional parameter optimization. In this sense, recognizing machine learning as function-space optimization is not a philosophical preference but a mathematical clarification of what is actually being learned and optimized.

Constrained optimization is equally indispensable, because most meaningful learning and control problems are not unconstrained. Constraints encode prior knowledge, physical laws, safety requirements, fairness criteria, resource budgets, and structural properties such as smoothness, sparsity, or stability. In infinite-dimensional settings, they arise naturally as norm-ball constraints in RKHSs, PDE state equations, probability-measure constraints, Lipschitz or bounded-variation conditions, and mean-field or transport constraints. Ignoring them may yield mathematically elegant but practically inadmissible solutions. Thus a rigorous theory of constrained infinite-dimensional optimization---and of the conditions under which nontrivial multipliers exist and first-order necessary conditions such as Fritz John and KKT hold---is essential for turning function-space learning into a principled, reliable, and deployable methodology.

More precisely, let $V$ be a complete metric space, $X$ be a
Banach space with its dual space $X^{'}$, and $E$ be
a nonempty closed convex subset of $X$.  Given continuous mappings
 $f_0: V\rightarrow
\dbR$ and  $f:
V\rightarrow X$, and a closed subset $\mathcal{K}\subseteq V$, we consider the following optimization
problem:
$$
\text{\bf (P)}\q\q\mbox{Minimize }  f_0(u),
\q\mbox{ subject to } f(u)\in E\mbox{  and  } u\in   \mathcal{K}.
$$

When  $V=\dbR^n$,  $X=\dbR^m$ (for $n, m \in
\dbN$),  $f=(g_1, g_2, \cdots, g_{m_1}, h_{m_1+1}, \cdots, h_{m})^\top$ with $m_1 \in
\dbN$ and $m_1\le m$, $g_{i}: \dbR^{n}\to \dbR$ ($i=1,2, \cdots, m_1$) and $h_{j}: \dbR^{n}\to \dbR$  ($j=m_1+1, m_1+2, \cdots, m$),  $\mathcal{K}=V$,  and
$$
E=\Big\{ (x_1, x_2, \cdots, x_{m_1}, 0, \cdots,
0)^\top\in\dbR^m\  \Big|\  x_{i}\le 0, i=1,  2,
\cdots, m_1  \Big\},
$$
the problem {\bf (P)}  can be rewritten in the classical
finite-dimensional form:
\begin{eqnarray*}\label{MP} \text{\bf (FP)}\q\q
\left\{
\begin{array}{ll}
\mbox{Minimize }  f_0(u),\quad   u\in\dbR^n,\\
\ns\ds \mbox{subject to } g_{i}(u)\le 0, i=1,2, \cdots, m_1;\\
\ns\ds\quad\quad\mbox{ and}\  \    h_{j}(u)=0, j=m_1+1, m_1+2, \cdots, m.
\end{array}
\right.
\end{eqnarray*}
Assuming  that  $f_0$,  $g_{i}$ $(i=1, 2, \cdots,
m_1)$,  and $h_{j}$ $(j=m_{1}+1,  \cdots,
m)$ are sufficiently smooth, a well-known first-order necessary optimality condition for the
problem {\bf (FP)}  (e.g., \cite[Page 199]{BSS})
states  that  if $\bar u\in\dbR^n$ solves {\bf
(FP)}, then there exist $z_{i}\ge 0$ ($i=0, 1,
\cdots, m_1$) and $z_{j}\in \dbR$ ($j=m_1
+1, \cdots, m$), not  all zero, satisfying  the complementary slackness
condition  $z_{i}g_{i}(\bar u)=0$
 $(i=1, 2, \cdots, m_1)$,   such that
\begin{equation}\label{1st condition for mp}
z_0  f_0^{\prime}(\bar
u)+\sum_{i=1}^{m_1}z_{i}  g_{i}^{\prime}(\bar
u)+\sum_{j=m_1 +1}^{m}z_{j} h_{j}^{\prime}(\bar
u)={\bf 0}.
\end{equation}
Here, $\varphi^{\prime}(\bar
u)$ denotes the gradient of  $\varphi:\dbR^n\to \dbR$ at $\bar u$ for $\varphi\in\{f_0, g_{i}, h_{j}\}$, where $i=1, 2, \cdots,
m_1$ and $j=m_{1}+1, \cdots,
m$.

This first-order necessary condition (\ref{1st condition for mp}), known as the Fritz John condition,
is widely used  in optimization theory.  This condition was originally
  introduced by F. John in \cite{John1948}  for the case where $m_1 = m$. In the  condition (\ref{1st condition for mp}), the non-zero vector $${\bf z}=(z_0, z_1,\cdots, z_{m_1}, z_{m_1
+1}, \cdots, z_m )^{\top}$$ is referred to as a  {\it nontrivial
Lagrange multiplier}.  The requirement that  ${\bf z}\not=0$
is essential   to ensure the validity of this first-order necessary condition; otherwise,  the condition (\ref{1st condition for mp}) trivially reduces to ``{\bf 0}={\bf0}", which  provides no
 useful information about the solution $\bar u$.

If  $z_0=0$,   the vector  ${\bf z}$ is referred to as a {\it
singular Lagrange multiplier}. If $z_0\not=0$,
 we may  assume without loss of generality  that $z_0=1$.
 In this case,  ${\bf z}$  is  called a {\it normal
Lagrange multiplier} and (\ref{1st condition for mp}) becomes
 the {\it KKT condition}, which was independently derived by W. Karush in \cite{K1939}, and by H. W. Kuhn and A. W. Tucker in \cite{KK1951}.

In the general setting  where $V$, $X$, $E$, $f_0$, $f$ and $\mathcal{K}$ are defined as for the problem
{\bf (P)}, if analogous first-order necessary conditions hold, they are likewise referred to as the Fritz John condition and the KKT condition, respectively.

The KKT condition is stronger than the Fritz John condition
 because it not only establishes  that $z_0=1$, but also
 explicitly incorporates
  information about the objective
   functional $f_0$. However,  the KKT condition does not generally
    hold for the problem {\bf
(FP)} (see \cite[Page 185]{BSS} for a counterexample).  It can be shown that the KKT condition holds
 for  {\bf
(FP)},  if and only if
\begin{eqnarray}\label{1.20-eq1}
\begin{array}{rl}
&\displaystyle  -f_0^{\prime}(\bar u)\in  f^{\prime}(\bar
u)^*E_1,\\[2mm]
&\displaystyle\quad  \mbox{  with   } E_1=\Big\{
(z_1, z_2, \cdots,
 z_m)^{\top}\in\dbR^m \ \Big|\  z_{i}\geq 0 \mbox{ and }  z_i g_i(\bar{u})=0,  i=1,  2,
\cdots, m_1 \Big\}.
\end{array}
\end{eqnarray}
Here and henceforth,   for  a bounded linear operator $F$,  we
denote its conjugate operator by $F^*$.  Therefore,  in
\eqref{1.20-eq1},  $  f^{\prime}(\bar u)^*=\big(
g_{1}^{\prime}(\bar u),\cds,  g_{m_1}^{\prime}(\bar u),   h_{m_1+1}^{\prime}(\bar
u),\cds,   h_{m}^{\prime}(\bar u)\big)$.   This implies
 that certain  constraint qualifications for the equality and inequality constraints, such as the Mangasarian-Fromovitz constraint qualification in \cite{MF1967},  are used  to
 guarantee  the validity of the  KKT condition (or, equivalently,  the existence of a normal Lagrange multiplier). On the other hand,
when $X$ is finite-dimensional,
 no additional  conditions are required for the Fritz John condition to hold,  even if $V$ is an infinite-dimensional Banach space (see \cite{Clarke1976}).

Generally speaking, there are two primary  methods for deriving the Fritz John condition in optimization theory:
 \begin{itemize}
   \item  The  first approach, known as the separation method,  utilizes
   the separation theorem in the range space $X$ of $f$ or the
   product space $\mathbb{R}\times X$ to find  Lagrange multipliers.
   This method has been widely explored in the literature \cite{BS, Kurcyusz1976, Robinson1976}.

 \item The second approach, the penalty functional method, involves introducing suitable penalty functionals to construct a sequence of unconstrained approximating optimization problems. By passing to the limit in the first-order necessary conditions satisfied by the solutions of these approximating problems, one obtains the desired Lagrange multipliers. This method has been used in \cite{Borgens2020, Clarke1976, FF, YXQ2007}, among others.
\end{itemize}

The penalty functional method offers distinct  advantages over the
separation method. While
the separation method only confirms the existence of Lagrange multipliers,  the penalty functional method constructively identifies a specific   multiplier as the
 limit of a  sequence, thereby yielding richer information.
  Further enhanced results regarding  the Fritz John condition can be found in \cite[Proposition 5.2.1, Page 284]{Bertsekas2003}, \cite{Bertsekas2006} and \cite{guoYezhang2013}.

Both methods work effectively when $X$ is finite-dimensional. By the separation theorem, any nonempty convex set in a finite-dimensional space has a supporting hyperplane at any of its boundary points. Likewise, any sequence on the unit sphere in a finite-dimensional space has a subsequence that converges to a nonzero vector. These two properties ensure the existence of nontrivial Lagrange multipliers via the separation and penalty functional methods, respectively. However, the situation changes drastically when the range space
$X$ of the constraint map $f$ in the problem {\bf (P)} is infinite-dimensional.
%
This
stems  from the inherent complexities of infinite-dimensional spaces:
\begin{itemize}
  \item The existence  of a supporting  hyperplane
   for a nonempty convex set is not guaranteed unless the set satisfies  additional
   conditions,  such as possessing
    a nonempty relative interior   (see \cite{klee1969}).

  \item A sequence residing on the unit sphere of an infinite-dimensional space
  may only  exhibit a weak (or weak*)
    convergent  subsequence, which could  potentially  converge to zero.
\end{itemize}

Hence,
for the general infinite-dimensional optimization problem {\bf(P)}, the existence of a nontrivial Lagrange multiplier cannot be guaranteed without additional assumptions  (see \cite{Brokate1980}). To derive   first-order necessary conditions for  constrained
infinite-dimensional optimization problems  {\bf (P)}, various
 constraint qualifications have been proposed in the  literature:
\begin{itemize}
  \item

To  apply the separation method when  $X$
is a Banach space,  specific  requirements are imposed on the set
\begin{equation}\label{c(bar u)}
\mathcal{Z}(\bar u) \deq f^\prime(\bar
u)\big(\mathcal{R}_{\mathcal{K}}(\bar u)\big)-
\mathcal{R}_{E}(f(\bar u)).
\end{equation}
For instance,
 $\mathcal{Z}(\bar u)$ is often required to have a nonempty
 relative interior,  or its   closure must not encompass
 the entire space   $ X$.
Here and throughout,  $f'(\bar
u)$ denotes the Fr\'echet derivative of  $f$ at $\bar u$, and
$\mathcal{R}_{D}(v)$ represents  the radial cone (see Definition \ref{llzd1}
in the next section) of a set $D\subseteq X$ at $v\in D$.
 For  further  details, see \cite{BS,Kurcyusz1976,Robinson1976,Kurcyusz1979}
and the references  therein.

\item For the penalty functional  method, the existence of nontrivial
Lagrange multipliers for solutions to the problem {\bf (P)} was
established  in \cite{FF}. This requires a   constraint qualification
formulated using  properties  of variation sets on a sequence of approximate optimal solutions (see \cite[Corollary 6.4.5, Page 267]{FF} or the condition  ${\bf (H_4)}$
in Remark \ref{propsition 2.1} of this paper).
\end{itemize}

While the aforementioned
 constraint qualifications can guarantee
  the existence of nontrivial Lagrange multipliers for
  the infinite-dimensional optimization problem {\bf (P)},
    it is often challenging to verify them directly.
For example, consider a special case where $V$ is a Banach space,
$\mathcal{K}= V$   and $E=\{0\}$.
In this scenario, $\mathcal{Z}(\bar u) =\Im (f^\prime(\bar u))$ in (\ref{c(bar u)}),
where $\Im (f^\prime(\bar u))$  denotes
 the range of the bounded linear operator $f^\prime(\bar u):V\to X$.
 Since $\Im (f^\prime(\bar u))$  is a subspace of $X$,
 its relative interior is nonempty if and only if it is closed.
Although  alternative characterizations for the closedness of  $\Im (f^\prime(\bar u))$
exist  (e.g., \cite[Theorem 5.2, Chapter IV]{Kato} or Proposition \ref{prop21}
in  this paper), verifying  them directly in infinite-dimensional spaces is generally quite difficult.
  We refer
to Section \ref{sec-5} for  a more detailed  discussion  on this aspect.
 Furthermore,   verifying the constraint qualification in \cite[Corollary 6.4.5, Page 267]{FF}
 poses an even  greater challenge, as it involves the  properties of a sequence of
  approximate optimal solutions. Therefore,
 introducing  new,  verifiable  conditions for the existence of nontrivial Lagrange multipliers
  in the  problem {\bf (P)} is highly desirable.

To this end, in this work we introduce a finite codimensionality condition,
rooted in functional analysis,  to establish an enhanced Fritz John (first-order necessary)
condition for the problem {\bf (P)}. Notably, our finite codimensionality condition
 is equivalent to suitable {\it a priori} estimates, making it more practical to  verify
  than  conventional constraint qualifications (see Section \ref{sec-4} for  detailed explanations).
This approach provides a unified framework for establishing first-order
necessary conditions for various constrained  optimal control problems,  which
can be viewed  as specific
 instances of  infinite-dimensional optimization problems.
 To demonstrate the verifiability of this method, we present some concrete
  examples in Section  \ref{sec-4}.
Additionally,  unlike typical nonlinear optimization problems
formulated  on topological vector spaces, the space $V$ in the problem {\bf (P)}   is
only  assumed to be a complete metric space.  Consequently, our
  results remain  applicable even if $V$ lacks a linear structure.

It is important to note that various finite codimensionality conditions have been employed in \cite{FF1987, FF, ly, LY91, 2, LLZ} (and the references therein) to establish Pontryagin maximum principles for infinite-dimensional optimal control problems with state constraints. However, in each of these works, the finite codimensionality condition was either assumed directly without verification, or verified only for specific control systems (e.g., in \cite{LLZ}). In a related vein, \cite[Chapter 6]{Troltzsch2010} discussed the KKT conditions for abstract constrained optimization problems under the Robinson constraint qualification, together with their applications to optimal control problems for PDEs with state constraints. Another related work, \cite{i}, employed Clarke tangent cones and subdifferentials to derive optimality conditions for optimization problems with inequality   constraints  via a finite codimensionality condition.

Despite these contributions, the constraint qualifications proposed in these earlier works generally remain difficult---and often prohibitively so---to verify in practice. Indeed, to the best of our knowledge, no systematic investigation has yet been undertaken into verifiable sufficient conditions for constraint qualifications in infinite-dimensional optimization problems. This gap is precisely what the present paper addresses.

Moreover, the core ideas of the finite codimensionality method developed here can be further extended to study the KKT conditions for infinite-dimensional optimization problems. In particular, in Section \ref{sec-4}, we introduce a constraint qualification for the KKT condition of the problem {\bf (P)} that is equivalent to a suitable \textit{a priori} estimate, thereby providing a verifiable sufficient condition.

The remainder of this paper is organized as follows. In Section \ref{sec-ex}, we introduce the finite codimensionality condition to establish an enhanced Fritz John condition for the problem {\bf (P)}, together with an application to optimal control problems. Section \ref{sec-4} provides several equivalent characterizations of the finite codimensionality condition. These general theoretical results are then applied to specific optimal control examples, including a deterministic evolution equation with an end-point constraint, an elliptic equation with a pointwise state constraint, and a stochastic differential equation with an end-point constraint. In Section \ref{sec-5}, we present equivalent \textit{a priori} estimates for closed-range operators and discuss further perspectives. Finally, the proofs of several preliminary results are deferred to Appendices A--D.

\section{First-order necessary  conditions  for the optimization problem}\label{sec-ex}

In this section, we derive a Fritz John
(first-order necessary) condition  for the problem
{\bf (P)} under a finite
codimensionality condition.
As an immediate application, this result
 is  utilized to study constrained optimal control problems.

\subsection{Preliminaries}

First, we introduce some standard notation. Let $\bd_V(\cdot, \cdot)$ denote  the metric on $V$.
We use $|\cdot|_X$ to represent the norm of the
Banach space $X$ and $\langle\cdot,
\cdot\rangle_{X',  X}$  to represent  the duality pairing
 between $X'$ and $X$.  In particular,  when $X$ is a Hilbert space, we simply
write $(\cdot,\cdot)_X$ for the inner
product on $X$.  For another Banach space $Z$,
$\mathcal{L}(X; Z)$  denotes
the space  of all bounded linear operators from $X$ to $Z$.
For simplicity, we write  $\mathcal{L}(X)$ in place of  $\mathcal{L}(X; X)$.

Let $D$ be a subset of $X$. We  denote  by  $\overline{D}$ the
closure of $D$,  by $\overline{\mbox{co}}D$ the closed convex hull
of $D$, and by $\mbox{Int} D$ the interior of $D$. For any $\alpha\in \mathbb{R}$,
we define $\alpha D\deq\big\{\alpha x\in X\;\big|\;x\in D\big\}$. The
subspace spanned by $D$ is defined as
$$ \mbox{span}D\deq\Big\{\sum_{i=1}^{n}\alpha_i  x_i \in X\ \Big|\ \alpha_i\in \mathbb{R}, x_i\in D \mbox{ for } i=1,  \cdots, n \mbox{ and } n\in \dbN\Big\}.$$
We write $\overline{\mbox{span}}D$ for the closure of span$D$.

For any two subsets $D_1$ and $D_2$ of $X$,  we set
$$
D_1\pm D_2\deq\big\{ x_1\pm
x_2\in X \;\big|\; x_1\in D_1,
x_2\in D_2 \big\}\
\mbox{ and }\
D_1\setminus D_2\deq\big\{ x\in X\; \big|\; x\in
D_1 \mbox{ and }x\notin D_2 \big\}.
$$

Next, we recall the definition of a topological complement for a closed subspace of $X$.
\begin{definition}
Let $X_1$ be a  closed subspace of $X$. A subspace
$X_2$ of $X$ is said to be a topological complement  of $X_1$,  if $X_2$ is closed,
$X_1\cap X_2=\{0\}$ and $X_1 + X_2=X$.
\end{definition}

Note that if $X_1$ and $X_2$ are complementary subspaces of $X$,
then every $x\in X$ can be uniquely written as $x = x_1 + x_2$, with $x_1 \in X_1$ and $x_2\in X_2$.
We define the projection operator $\Pi_{X_1}$ from $X$ onto the closed subspace $X_1$ by $\Pi_{X_1}(x)=x_1$  for any $x=x_1+x_2\in X$ with $x_1 \in X_1$ and $x_2\in X_2$.
The operator $\Pi_{X_2}$ can be defined in a similar manner. In particular, when $X_1$ is a finite-dimensional  subspace of $X$, $X_1$ always has a topological complement,  and  $\Pi_{X_1}$ is well-defined.   When $X$ is  a Hilbert  space, we denote
the orthogonal complement of a closed subspace $X_1$
by $X_1^\bot$.

Let $ X\times Z$ denote the product space of two Banach spaces $X$ and $Z$.  For a set $S\subseteq X\times Z$, we define the projection onto $X$ as
$$\pi_{X}(S)\deq \big\{ x\in X \ \big|\  (x,z)\in S
\mbox{ for some }z\in Z \big\}.$$
The projection $\pi_{Z}(S)$  is  defined analogously.

Furthermore, for any $x_0\in X$ and
$\rho>0$, let
$
B_X(x_0,\rho)\deq\big\{x\in X \  \big|\
|x-x_0|_X\leq \rho\big\}
$ denote the closed ball.
Similarly, we  define $B_{X'}(x'_0,\rho)$ for $x'_0\in X'$,
 and  $B_V(v_0,\rho)\deq\big\{v\in V \  \big|\
 \bd_V(v,v_0)\leq \rho\big\}$ for $v_0\in V$.

Recall that a Banach space $X$ is strictly convex,  if   any $x,y\in X$ with $|x|_{X}=|y|_{X}=1$ and
 $x\neq y$ satisfy
\begin{equation}\label{2}
\Big|\frac{x+y}{2}\Big|_{X}<1.
\end{equation}
Clearly, any Hilbert space is  strictly convex.

For two Banach spaces $X$ and $Z$,  a map  $\varphi:  X\to Z$
 is said to be directionally  differentiable at $x$ along a
 direction $v\in X$, if  the following limit exists:
$$
\displaystyle\lim\limits_{h\rightarrow 0^+}\frac{\varphi(x+hv)-\varphi(x)}{h}.
$$
This limit is denoted by $\varphi'(x;  v)$ and
is called the directional derivative of  $\varphi$ at $x$ along the direction $v$.
 In addition,  if  $\varphi$ is directionally
  differentiable at $x$ along every direction   $v\in X$,  then $\varphi$ is said to be directionally differentiable at $x$.
Furthermore, if  there exists  a bounded linear operator $\varphi'(x)\in \mathcal{L}(X; Z)$
such that
$$\lim_{y\rightarrow x}\frac{|\varphi(y)-\varphi(x)-\varphi'(x)(y-x)|_Z}{|y-x|_{X}}=0,$$
then $\varphi$ is said to be Fr\'{e}chet differentiable at $x$,
 and $\varphi'(x)$ is called the Fr\'{e}chet derivative of $\varphi$  at $x$.
Clearly, when $\varphi$ is  Fr\'{e}chet differentiable at $x$,
it is also directionally differentiable at $x$,  and
$\varphi'(x;  v)=\varphi'(x)v$ for any $v\in X$. If
$\varphi$ is Fr\'{e}chet differentiable on a neighborhood of  $x$ and its Fr\'{e}chet derivative  is continuous (with respect to the operator  norm topology) at $x$,
then $\varphi$  is said to be
continuously differentiable at $x$.
When $\varphi$ is directionally differentiable (respectively,
Fr\'{e}chet differentiable, or continuously differentiable)  at every $x\in X$,
we simply  say that  $\varphi$ is directionally differentiable (respectively,
Fr\'{e}chet differentiable, or continuously differentiable).

\medskip

Furthermore,  we recall some fundamental  concepts  from
 variational analysis;   more
detailed discussions   can be found in \cite{00,BS, FF}.
\begin{definition}\label{llzd1}
For a nonempty closed convex subset $E$ of  $X$ and  $e\in E$, the set
$$\mathcal{R}_{E}(e)=\bigcup_{\a\in [0,\infty)} \a(E-e)$$
is called the radial cone of $E$ at $e$. The
closure of $\mathcal{R}_{E}(e)$, denoted by
$\cT_{E}(e)$, is called the tangent cone of $E$ at
$e$.
\end{definition}
\begin{definition}\label{llzd3}
The  normal cone  of  $E$  at
$e\in  E$,   denoted by   $\cN_E(e)\subseteq X'$,  is  defined as
$$
\cN_E(e)=\big\{ w\in X' \ \big|\  \langle w, \tilde
e -e\rangle_{X',  X}\leq 0,\ \   \forall\; \tilde
e\in E \big\}.
$$
\end{definition}
It is easy to
verify  that $\cN_E(e)$ can equivalently  be defined   as:
$$
\cN_E(e)=\big\{ w\in X' \ \big|\   \langle w,
v\rangle_{X',  X}\leq 0,\ \  \forall\;   v\in
\cT_E(e)  \big\}.
$$

Let $K$ be a subset of $X$.
The
distance function $\mbox{dist}(\cd, K):X\to \dbR$  is defined as
$$
\mbox{dist}(x, K)\deq \inf\limits_{y\in
K}\big|y-x\big|_X,\q \forall\;x\in X.
$$
For any $x\in X$,   the metric projection of $x$ onto $K$ is defined as
\begin{equation}\label{metric proj}
\cP_{K}(x)\deq\big\{ \bar e\in K \  \big|\
\mbox{dist}(x, K)=|x-\bar e|_X   \big\}.
\end{equation}
When $X$ is a reflexive Banach space and $K$ is a nonempty closed convex
subset of $X$, $\cP_{K}(x)$ is nonempty for any $x\in X$. In addition, when $X$ is strictly convex,
 it follows  that $\cP_{K}(x)$ is a singleton.

For    a nonempty closed convex subset $E$ of a Banach space
$X$, the distance
function $\mbox{\rm dist}(\cd,E)$ is Lipschitz continuous and convex on $X$.
   Its  subdifferential   at $x\in X$ is defined as
$$
\partial \mbox{\rm dist}( x,E)\deq\big\{
w\in X^{\prime} \;\big|\; \mbox{\rm dist}(\tilde
x,E) -\mbox{\rm dist}( x,E)\ge  \langle w,
\tilde x -x\rangle_{X^{\prime},  X},\ \ \forall\;
\tilde x\in X  \big\}.
$$

We recall the following properties of $\partial \mbox{\rm dist}(\cd,E)$.
\begin{lemma}\label{4.3-prop11}{\rm\cite[Proposition 3.11, Page 146]{2}}
For any $x\in X$,  the set $\partial \mbox{\rm
	dist}( x,E)$ is a nonempty,  convex,
weak$^*$-compact subset of $X^{\prime}$,  and
$\partial \mbox{\rm dist}( x,E)\subseteq
B_{X^{\prime}}(0,1)$.  Moreover,  $\mbox{\rm dist}( \cdot,E)$  is directionally differentiable,
 and  for any $x,v\in X$,
\begin{equation}\label{4.3-eq67}
		\mbox{\rm dist}^{\prime}(x, E; v)=
		\max\big\{\langle\zeta, v\rangle_{X', X}\in \mathbb R\ \big|\  \zeta\in\pa
		{\rm dist}(x,
		E)\big\},
	\end{equation}
where $\mbox{\rm dist}^{\prime}(x, E; v)$ denotes the directional derivative of  $\mbox{\rm dist}( \cdot,E)$ at $x$ along the direction $v$.
\end{lemma}

When  $X$ is
a reflexive Banach space,
\begin{equation}\label{partial dist}
\partial \mbox{\rm dist}( x,E)=
\cN_E(\bar e)\bigcap \big\{ w\in
B_{X^{\prime}}(0,1)\ \big|\  \langle w,  x -\bar
e\rangle_{X^{\prime},  X}=|x-\bar e|_{X}\big\}, \q \forall\;\bar e\in \cP_{E}(x).
\end{equation}
Here, $\cP_{E}(x)$ is the metric projection of $x$ onto $E$ defined by (\ref{metric proj}). Note that  even if
the set $\cP_{E}(x)$ is not   a singleton,  the
set on the right-hand side of  $(\ref{partial
dist})$ remains invariant  for any choice of  $\bar e\in
\cP_{E}(x)$.  Furthermore, when $X^{\prime}$ is strictly convex and $x\notin E$,  $\partial \mbox{\rm dist}( x,E)$ is a singleton. In this case,
 the unique element $w\in \partial \mbox{\rm dist}( x,E)$ satisfies  $|w|_{X'}=1$.
In particular, when $X$ is a Hilbert
space and $x\notin E$,
 the unique $w\in \partial\mbox{\rm dist}(x,
E)$ is the element (via the Riesz-Fr\'echet
isomorphism)  corresponding     to
\begin{equation}\label{gradent dist}
\frac{x-\cP_{E}(x)}{|x-\cP_{E}(x)|_{X}}\in X.
\end{equation}
We refer to $\cite[\mbox{Example } 2.130,
\mbox{Page } 89]{BS}$ for further  details.

\medskip

The following concept of variation originates
 from  \cite[Page 94]{FF} and \cite[Page 43]{F-F1991}.
\begin{definition}\label{llzd4}
We call $\xi\in X$ a variation of  $f: V\rightarrow X$ with respect to
a subset $\mathcal{K}\subseteq V$ at $e\in \mathcal{K}$,  if
there exists a sequence $\{h_k\}_{k=1}^\infty\subseteq
(0,+\infty)$  with $\lim\limits_{k\rightarrow\infty}h_k=0$  and a
sequence  $\{e_k\}_{k=1}^\infty\subseteq \mathcal{K}$  such that
$$
\bd_V(e_k,e)\leq h_k\q \mbox{ and }\q
\lim\limits_{k\rightarrow\infty}\frac{f(e_k)-f(e)}{h_k}=\xi.
$$
The set of all
variations of $f$  with respect to the subset $\mathcal{K}\subseteq V$ at $e\in \mathcal{K}$ is  denoted by ${\rm Var}_{\mathcal{K}} f(e)$.
\end{definition}

In the standard definition of  variation (e.g., \cite[Page 94]{FF} and \cite[Page 43]{F-F1991}),
the domain  $\mathcal{D}(f)$ of $f$ may not be the entire
 metric space $V$. The  variation of $f$ at
$e\in \mathcal{D}(f)$ defined in those works
  coincides with  Definition \ref{llzd4} by taking $\mathcal{K}=\mathcal{D}(f)$.
The following example demonstrates
 that  the variation of $f$ generalizes  the notion of the  directional
derivative.
\begin{example}\label{remark2.2}
If  $V$ is  a  reflexive  Banach space,
$\mathcal{K}$ is a closed convex subset of V, and  $f: V\rightarrow X$ is locally Lipschitz
continuous and directionally  differentiable,  then for any
$e\in \mathcal{K}$,
\begin{equation}\label{appendix}
\big\{ f'(e; v)\in X\   \big|\  v\in  \cT_{\mathcal{K}}(e) \cap
B_V(0, 1) \big\}\subseteq \mbox{\rm Var}_{\mathcal{K}}f(e).
\end{equation}
In addition,  assume that either $f$ is Fr\'echet differentiable,   or $V$ is finite-dimensional.   Then
\begin{equation}\label{appendix_add}
\mbox{\rm Var}_{\mathcal{K}}f(e)=\big\{ f'(e;v)\in X\    \big|\  v\in  \cT_{\mathcal{K}}(e) \cap
B_V(0, 1) \big\} .
\end{equation}
As a special case, when $f$ is Fr\'echet differentiable and $e\in
{\rm Int}\mathcal{K}$, we have
$$
\mbox{\rm Var}_{\mathcal{K}}f(e)=\big\{ f'(e)v\in X    \;\big|\;  v\in  B_V(0, 1)
\big\}.
$$
For the convenience of the readers,   detailed proofs of
\eqref{appendix} and \eqref{appendix_add} are provided  in Appendix \ref{A}.
\end{example}

For a set-valued map $F:  \mathcal{D}(F)\subseteq V\rightsquigarrow X$,
the domain $\mathcal{D}(F)$ of $F$ is the set of all elements $e\in V$ such that $F(e)$ is nonempty;
that is,
 $\mathcal{D}(F)\deq\big\{e\in V~\big|~ F(e)\neq \emptyset \big\}$.
The graph $\mbox{Gph}(F)$ of $F$,   viewed as a subset of the product space $V\times X$,
 is defined as
  $\mbox{Gph}(F)\deq\big\{(e,\xi)\in V\times X~\big|~\xi\in F(e)\big\}$.
\begin{definition}
A  set-valued map $F: V\rightsquigarrow X$ is said to be
 locally bounded at $ e_0\in \mathcal{D}(F)$,  if there  exist
  two positive constants $\delta $ and $C $
   such that  for any $e\in \mathcal{D}(F)$ with  $\bd_{V}(e,e_0)\le\delta$
    and any $\xi\in F(e)$, it holds that $|\xi|_X \le C$.
\end{definition}
\begin{remark}\label{RRR}
For a subset  $\mathcal{K}\subseteq V$,
if  $f:V\rightarrow X$ is locally Lipschitz
continuous at $e_0\in\mathcal{K}$, meaning
there exist two positive constants $L$ and  $\delta$ such that
$$|f(e_1)-f(e_2)|_{X}\le L \bd_{V}(e_1,e_2),
 \q \forall\; e_1,e_2\in V \mbox{ with } \bd_{V}(e_1,e_0)\le\delta\mbox{ and }
 \bd_{V}(e_2,e_0)\le\delta, $$
  then the set-valued map
 $\mbox{\rm Var}_{\mathcal{K}}f(\cdot) $ is locally bounded at   $e_0\in\mathcal{K}$.

 Indeed, for any $e\in \mathcal{K}$ with  $\bd_{V}(e,e_0)\le\frac{\delta}{2}$ and any $\xi\in \mbox{\rm Var}_{\mathcal{K}} f(e)$, it follows from Definition $\ref{llzd4}$ that
  there exists  a sequence $\{h_k\}_{k=1}^\infty\subseteq
(0,+\infty)$  with $\lim\limits_{k\rightarrow\infty}h_k=0$ and a
sequence  $\{e_k\}_{k=1}^\infty\subseteq \mathcal{K}$ with $\bd_V(e_k,e)\leq h_k$, such
that
$$\xi=
\lim\limits_{k\rightarrow\infty}\frac{f(e_k)-f(e)}{h_k}.$$
For sufficiently large $k$, we have $h_k\leq \frac{\delta}{2}$, which implies that
$\bd_V(e_k, e_0)\leq \bd_V(e_k, e)+\bd_V(e, e_0)\leq\delta$. By the  local Lipschitz continuity of $f$ at $e_0$,
it follows that $$
|\xi|_{X}=\Big|
\lim\limits_{k\rightarrow\infty}\frac{f(e_k)-f(e)}{h_k}\Big|_X
\le \varlimsup\limits_{k\rightarrow\infty} \frac{L \bd_{V}(e_k,e) }{h_k}\le L.
$$This establishes  the local boundedness of $\mbox{\rm Var}_{\mathcal{K}}f(\cdot) $
 at   $e_0\in\mathcal{K}$.
\end{remark}

To conclude  this subsection, we recall several  concepts concerning
finite codimensionality (see,  e.g., \cite{Brezis, 2}).
\begin{definition}\label{llzd5}
A  linear subspace $X_1$ of  $X$  is said to be finite codimensional
\footnote{In \cite[Page 351]{Brezis},  a linear subspace $X_1$ is considered
 finite codimensional, if  there exist  $m\in\dbN$ and linearly independent $x_1, x_2, \cdots,
x_m\in X\setminus X_1$ such that $X_1+\mbox{span}\big\{x_1,
x_2,  \cdots, x_m \big\}=X$. Clearly, this formulation
 in \cite{Brezis}
  is a purely  algebraic concept  involving only  linear operations.
 Throughout this paper,  however, we  adopt
  Definition \ref{llzd5} for the  finite codimensionality of  a
   linear subspace $X_1$, which explicitly
    incorporates  the topological closure of $X_1$.  When $X_1$ is closed,
    Definition \ref{llzd5} precisely coincides with the definition provided  in \cite{Brezis}.},  if
there exist  $m\in\dbN$ and linearly independent $x_1, x_2, \cdots,
x_m\in X\setminus X_1$  such that $$\overline{X_1}+\mbox{span}\big\{x_1,
x_2,  \cdots, x_m \big\}=X. $$
\end{definition}

\begin{definition}\label{llzd6}
A  subset   $D$ of  $X$ is said to be finite codimensional  in $X$,  if there exists
an $x_0\in\coh D$  such that

\smallskip

 ${\bf (i)}$
$\overline{\mbox{span}}\big\{D-x_0\big\}$  is a finite
codimensional subspace of $X$; and

\smallskip

 $\bf  (ii)$ $\coh \big(D-x_0\big)$
has at least  one interior point  relative to $\overline{\mbox{span}}\big\{D-x_0\big\}$.
\end{definition}

The following result demonstrates that the  finite
codimensionality of a subset $D$ is independent of the choice of $x_0 \in\coh D$ in Definition \ref{llzd6}.
\begin{proposition}\label{add**}
If the conditions
${\bf (i)}$ and ${\bf (ii)}$ stated in  Definition $\ref{llzd6}$  hold for some $x_0\in \coh D$,
then  they hold for all $x_0\in \coh D$.
\end{proposition}

Proposition \ref{add**} can be deduced from \cite[Proposition 3.1, Page 137]{2}. For the readers' convenience,
a proof of Proposition \ref{add**} is provided  in Appendix \ref{B}.

\medskip

Finally,  the following  lemma plays a crucial role in
establishing   the main result of this paper.

\begin{lemma}\label{codim converg}
{\rm ($\cite[\mbox{Lemma } 3.6,  \mbox{Page } 142]{2}$)}
Let $D$ be a finite codimensional subset of $X$,  and let
  $\{\varepsilon_{k}\}_{k=1}^{\infty}$ be a sequence of  nonnegative real numbers
   such that $\lim\limits_{k\to\infty}\varepsilon_{k}= 0$. Assume that $\{\Lambda_{k}\}_{k=1}^{\infty}\subseteq
X^\prime$ and that  $\Lambda_{k}$
converges $\mbox{weakly}^{*}$ to $\Lambda\in  X^\prime$ as $k\to \infty$.  If
there exists  a $\delta>0$ such that
\begin{enumerate}
  \item [${\bf (i)}$] $\left|\Lambda_{k}\right|_{X'} \geq \delta$ for all $k\in \dbN$; and
  \item [${\bf (ii)}$] $\left\langle \Lambda_{k}, x\right\rangle_{X', X}
\geq-\varepsilon_{k}$ for all $(x,k) \in D\times \dbN$,
\end{enumerate}
then $\Lambda\neq 0$.
\end{lemma}
For the readers' convenience, a detailed proof of Lemma \ref{codim converg}
is deferred to   Appendix \ref{C}.

\subsection{First-order necessary  conditions for  the problem {\bf (P)}}\label{subsection 2.2}

In this subsection, we establish   first-order necessary conditions for the  problem {\bf (P)}. We say
 $\bar{u}\in \mathcal{K}$ is a solution to the problem {\bf (P)},  if
$$
f_0(\bar{u})=\min\big\{  f_0(u)\in \mathbb R\  \big
|\  f(u)\in E\mbox{ and }u\in \mathcal{K} \big\}.
$$

In what follows, we assume the existence of   a  set-valued
map $\mathcal V: \mathcal{K}\rightsquigarrow \dbR\times X$ satisfying $\mathcal V(u)
\subseteq \mbox{\rm Var}_{\mathcal{K}}(f_0,f)(u)$ for any $
u\in\mathcal{K}$,  alongside    a  modulus of continuity
 $\ell: [0, +\infty) \to [0, +\infty)$
  with  $\lim\limits_{s\rightarrow 0^+}\ell(s)=0$, such that
\begin{enumerate}
\item [${\bf  (H_1)}$] For any
 $u^{\e}  \in\mathcal{K}$ with $\bd_V(\bar u,  u^{\e})\to 0$ as $\e\to 0^+$,  there exists an $\varepsilon_0>0$ such that  for any $(\xi_1,
\xi_2)^\top\in \mathcal V(\bar{u})$   there exists
  $(\xi^{\e}_1, \xi^{\e}_2)^\top\in \mathcal V(u^{\e})$ satisfying
$\xi^{\e}_1\rightarrow \xi_1$ in $\dbR$ as $\e\to 0^+$, and for all  $0<\varepsilon<\varepsilon_0$,
\begin{equation}\label{modulus-contin}
\big|\xi_2-\xi^{\e}_2\big|_X\leq \ell\big(\bd_V(\bar  u,  u^{\e})\big);
\end{equation}
\item [${\bf  (H_2)}$] $\pi_\dbR\big(\mathcal V(\cdot)\big)$   is   locally bounded at $\bar{u}$;
\item [${\bf  (H_3)}$] $\pi_{X}
\big(\mathcal{V}(\bar u)\big)-E$
is finite codimensional  in   $X$.
\end{enumerate}

The main result of this subsection, stated below,
provides a   necessary condition for solutions to the problem {\bf
(P)} under the assumptions  ${\bf  (H_1)}$--${\bf  (H_3)}$. Further
interpretations of these conditions  will be provided later.
\begin{theorem} \label{tt1}
Let $X$ be a reflexive Banach space such that its dual  $X^{\prime}$ is strictly convex,
and let $\bar u\in\mathcal{K}$ be a solution to the problem {\bf
(P)} with   $\bar{e}=f(\bar{u})$.  If the conditions  ${\bf  (H_1)}$--${\bf  (H_3)}$ are satisfied,  then
there exists  a non-zero pair  $(z_0, z)^\top\in
[0,+\infty)\times X'$ with $z\in
\cN_E(\bar{e})$,  such that the following Fritz John condition holds:
\begin{equation}\label{FJ condition}
z_0\xi_1+\langle z,  \xi_2\rangle_{X', X}\geq 0,
\q  \forall\; (\xi_1, \xi_2)^\top\in \mathcal
V(\bar{u}).
\end{equation}
Moreover, the following  enhanced results
hold true:

$(i)$ If $z_0\neq 0$, then  the following KKT condition holds:
$$
 \xi_1+\langle \tilde z,  \xi_2\rangle_{X', X}\geq 0, \q  \forall\; (\xi_1, \xi_2)^\top\in
\mathcal V(\bar{u})
$$
with $\tilde z=z/ z_0$;

$(ii)$  If $z\neq0$,  then there exists a
sequence $\{ u^k\}_{k=1}^\infty \subseteq \mathcal{K}$,
converging to $\bar u$  as $k\rightarrow
\infty$, such that
\begin{equation}\label{enhanced FJ eq1}
 \begin{cases}\ds
\ f(u^k)\notin E \text{ and } f_0(u^k)<f_0(\bar u) \hbox{ for sufficiently large $k$}, \\
\ns\ds\
\lim\limits_{k\rightarrow\infty}\text{\rm
dist}(f(u^k), E)=0,  \\
\ns\ds\
\lim\limits_{k\to\infty}
f_0(u^k)=f_0(\bar{u}).
\end{cases}
\end{equation}
Furthermore, when $X$ is a Hilbert space, the
sequence $\big\{\cP_{E}(f(u^{k}))\big\}_{k=1}^{\infty}$
satisfies  $\cP_{E}(f(u^{k}))\to f(\bar u)$
as $k\to \infty$, and
\begin{equation}\label{enhanced FJ eq3}
\big(\hat z,  f(u^k)-\cP_{E}(f(u^{k}))\big)_{ X}>0 \q \hbox{ for sufficiently large } k,
\end{equation}
where $\hat z\in X$ is the element corresponding
to $z \in X'$ via the Riesz-Fr\'echet
isomorphism.
\end{theorem}

\noindent  {\bf  Proof. }  The whole proof is divided
into  four steps.

\medskip

\noindent {\bf Step 1.  } For any $\e\in(0,  1)$ and
$u\in \mathcal{K}$, set
$$
\begin{array}{ll}
\Phi_\e(u)=
\Big\{\big[\mbox{dist}\big(f(u ),E\big)\big]^2+\big[
\big(f_0(u)-f_0(\bar{u})+\e\big)^+
\big]^2\Big\}^{1/2}&\end{array}
$$
and $ \mathcal F(u)=\big(f_0(u),
f(u)\big)^\top$,  where $s^+=s$ if $s\ge 0$,  and $s^+=0$ if $s<0$.
Then  $\Phi_\e(\cdot)$  is continuous  on $\mathcal{K}$, $\Phi_\e(u)>0$ for any $u\in\mathcal{K}$,
and
$$
\Phi_\e(\bar{u})=\e\leq
\inf\limits_{u\in \mathcal{K}} \Phi_\e(u)+\e.
$$
Since $\mathcal{K}$ is a complete metric
space equipped  with the metric $\mathbf{d}_V(\cdot, \cdot)$,   the Ekeland
variational principle guarantees   the existence of a $u^\e\in \mathcal{K}$  such that
$$
\Phi_\e(u^\e)\leq
\Phi_\e(\bar{u}),\quad\quad
\bd_V(\bar{u},u^\e) \leq \sqrt\e,
$$
and
\begin{equation}\label{66}
-\sqrt{\e}\bd_V(u^\e, u)  \leq \Phi_\e(u)
-\Phi_\e(u^\e),\quad \forall\; u\in \mathcal{K}.
\end{equation}

 \medskip

\noindent {\bf Step 2.} For any $\xi^\e
=\big(\xi^\e_1, \xi^\e_2\big)^\top \in\mathcal
V(u^\e)\subseteq \mbox{\rm Var}_{\mathcal{K}}(f_0,f)(u^\e)$,
it follows from  Definition \ref{llzd4} that  there
exists  a sequence $\{h^\e_k\}_{k=1}^\infty\subseteq (0,+\infty)$
with $\lim\limits_{k\rightarrow\infty}h_k^\e=0$  and a sequence
$\{u_k^\e\}_{k=1}^\infty\subseteq \mathcal{K}$ such that
\begin{equation}\label{llze1}
\bd_V(u^\e_k,u^\e) \leq h_k^\e
\end{equation}
and
\begin{equation}\label{llze1.1}
\mathcal F(u_k^\e)=\mathcal F(u^\e)+h_k^\e
\xi^\e+{\it  o}(h_k^\e),\  \mbox{ as
}k\rightarrow\infty,
\end{equation}
where ${\it  o}(h_k^\e)=\big({\it  o}_1(h_k^\e),  {\it  o}_2(h_k^\e)\big)^\top$ denotes a
higher-order infinitesimal with respect to $h_k^\e$. It follows from (\ref{llze1}) and (\ref{llze1.1}) that
\begin{eqnarray*}
&&\Phi_\e(u_k^\e)-\Phi_\e(u^\e)\\
&& = \frac{\big[ \big(f_0(u_k^\e) - f_0(\bar{u}) + \e\big)^+
\big]^2-\big[ \big(f_0(u^\e) - f_0(\bar{u}) + \e\big)^+ \big]^2}
{\Phi_\e(u_k^\e)+\Phi_\e(u^\e)}
\\[1mm]
&&\q + \frac{\big[\mbox{dist}\big(f(u^{\e}_{k}),E\big)\big]^2
- \big[\mbox{dist}\big(f(u^\e),E\big)\big]^2}
{\Phi_\e(u_k^\e)+\Phi_\e(u^\e)}\\[1mm]
&&=\frac{\big[ \big(f_0(u^\e)+h_k^\e
\xi^\e_1+{\it
o}_1(h_k^\e)-f_0(\bar{u})+\e\big)^+
\big]^2-\big[
\big(f_0(u^\e)-f_0(\bar{u})+\e\big)^+
\big]^2}{\Phi_\e(u_k^\e)+\Phi_\e(u^\e)}\\[1mm]
&&\quad +\frac{
\big[\mbox{dist}\big(f(u^\e)+h_k^\e\xi^\e_2+{\it
o}_2(h_k^\e),E\big)\big]^2
-\big[\mbox{dist}\big(f(u^\e),E\big)\big]^2}
{\Phi_\e(u_k^\e)+\Phi_\e(u^\e)}.
\end{eqnarray*}

Clearly,
\begin{eqnarray*}
&&\lim\limits_{k\rightarrow\infty}\frac{\big[ \big(f_0(u^\e)+h_k^\e
\xi^\e_1+{\it
o}_1(h_k^\e)-f_0(\bar{u})+\e\big)^+
\big]^2-\big[
\big(f_0(u^\e)-f_0(\bar{u})+\e\big)^+
\big]^2}{h_k^\e\big[\Phi_\e(u_k^\e)+\Phi_\e(u^\e)\big]} \\
&&=\frac{
\big(f_0(u^\e)-f_0(\bar{u})+\e\big)^+
\xi^\e_1
}{\Phi_\e(u^\e)}.
\end{eqnarray*}
In addition, by the Lipschitz continuity of distance functions and  Lemma \ref{4.3-prop11},
we obtain
\begin{eqnarray*}
&&\lim\limits_{k\rightarrow\infty}\frac{
 \mbox{dist}\big(f(u^\e)+h_k^\e\xi^\e_2+{\it
o}_2(h_k^\e),E\big)
- \mbox{dist}\big(f(u^\e),E\big) }
{h_k^\e }\\
&&=\lim\limits_{k\rightarrow\infty}\frac{
 \mbox{dist}\big(f(u^\e)+h_k^\e\xi^\e_2,E\big)
- \mbox{dist}\big(f(u^\e),E\big) }
{h_k^\e }\\
&&=\max\big\{\langle \psi, \xi^\e_2\rangle_{X', X}\in\mathbb R\ \big|\ \psi \in\partial \mbox{dist}\big(f(u^\e),E\big)\big\}.
\end{eqnarray*}
Then, by the  weak$^*$ compactness of $\partial \mbox{dist}\big(f(u^\e),E\big)$
 (the  subdifferential of  $\mbox{dist}\big(\cd,E\big)$ at $f(u^\e)$),
there is a $\psi_2^\e\in \partial \mbox{dist}\big(f(u^\e),E\big)$,
  such that
$$\lim\limits_{k\rightarrow\infty}\frac{
 \mbox{dist}\big(f(u^\e)+h_k^\e\xi^\e_2+{\it
o}_2(h_k^\e),E\big)
- \mbox{dist}\big(f(u^\e),E\big) }
{h_k^\e }= \langle \psi_2^\e, \xi^\e_2\rangle_{X', X}.$$
Therefore,
\begin{eqnarray*}
&&\lim\limits_{k\rightarrow\infty}\frac{
\big[\mbox{dist}\big(f(u^\e)+h_k^\e\xi^\e_2+{\it
o}_2(h_k^\e),E\big)\big]^2
-\big[\mbox{dist}\big(f(u^\e),E\big)\big]^2}
{h_k^\e\big[\Phi_\e(u_k^\e)+\Phi_\e(u^\e)\big]}  \\
&&=\lim\limits_{k\rightarrow\infty}\[\frac{
 \mbox{dist}\big(f(u^\e)+h_k^\e\xi^\e_2+{\it
o}_2(h_k^\e),E\big)
+\mbox{dist}\big(f(u^\e),E\big) }
{\Phi_\e(u_k^\e)+\Phi_\e(u^\e)}\\
&&\qq\qq\times\frac{
 \mbox{dist}\big(f(u^\e)+h_k^\e\xi^\e_2+{\it
o}_2(h_k^\e),E\big)
- \mbox{dist}\big(f(u^\e),E\big) }
{h_k^\e } \]  \\
&&=\frac{
  \mbox{dist}\big(f(u^\e),E\big)
}{\Phi_\e(u^\e)}\langle
\psi_2^\e, \xi^\e_2\rangle_{X', X}.
\end{eqnarray*}

On  the other hand, by  (\ref{66}),
$
\Phi_\e(u_k^\e)-\Phi_\e(u^\e)\geq -\sqrt{\e}h_k^\e.
$
Therefore, we have
\begin{equation}\label{2.14+}
\displaystyle
\lim\limits_{k\rightarrow\infty}\frac{\Phi_\e(u_k^\e)-\Phi_\e(u^\e)}{h_k^\e}=
\frac{
\big(f_0(u^\e)-f_0(\bar{u})+\e\big)^+
\xi^\e_1+ \mbox{dist}\big(f(u^\e),E\big)\langle
\psi_2^\e, \xi^\e_2\rangle_{X', X}
}{\Phi_\e(u^\e)}\geq -\sqrt{\e},
\end{equation}
for the above $\psi_2^\e\in \partial \mbox{dist}\big(f(u^\e),E\big)$.

Note that if $f(u^\e)\in E$,  then $\mbox{dist}\big(f(u^\e),E\big)=0$ and $\Phi_\e(u^\e)=
\big(f_0(u^\e)-f_0(\bar{u})+\e\big)^+$. In this case, we
 choose a $\psi_\e\in \partial \mbox{dist}\big(f(u^\e),E\big)$;  if $f(u^\e)\notin E$,
then by the strict convexity of $X^{\prime}$, $\partial \mbox{dist}\big(f(u^\e),E\big)$
 is a singleton,  and as the unique element
$\psi_2^\e$ in $\partial \mbox{dist}\big(f(u^\e),E\big)$,  $|\psi_2^\e|_{X'}=1$ and we define  $\psi_\e=\psi_2^\e$. Therefore, by \eqref{2.14+}, we can always find a $\psi_\e\in \partial \mbox{dist}\big(f(u^\e),E\big)$ such that
$$\frac{
\big(f_0(u^\e)-f_0(\bar{u})+\e\big)^+
\xi^\e_1+ \mbox{dist}\big(f(u^\e),E\big)\langle
\psi_\e, \xi^\e_2\rangle_{X', X}
}{\Phi_\e(u^\e)}\geq -\sqrt{\e}, \q \forall \; \big(\xi^\e_1, \xi^\e_2\big)^\top \in\mathcal V(u^\e).$$

Let
\begin{equation}\label{z0eze}
a_\e=\frac{\big(f_0(u^\e)-f_0(\bar{u})+\e\big)^+}{\Phi_\e(u^\e)}
\quad\mbox{ and }\quad
b_\e=\frac{\mbox{dist}\big(f(u^\e), E\big)\psi_\e}{\Phi_\e(u^\e)}.
\end{equation}
Then,
\begin{equation}\label{0010}
 a_\e \xi^\e_1+\langle b_\e, \xi^\e_2\rangle_{X', X}\geq -\sqrt{\e}, \q \forall \; \big(\xi^\e_1, \xi^\e_2\big)^\top \in\mathcal V(u^\e).
\end{equation}
Noting that $\psi_\e\in \partial \mbox{dist}\big(f(u^\e),E\big)$,  by the
definition of $b_\varepsilon$ in  (\ref{z0eze}), we have
\begin{equation}\label{0011}
\langle b_\e, x-f(u^\e)\rangle_{X', X}\leq 0,\quad \forall\; x\in E.
\end{equation}

Moreover, by the definitions of $\Phi_\e(u^\e)$, $a_\e$ and $b_\e$, we have
$a_\e^2+|b_\e|_{X'}^2=1$.
Consequently,  there exists  a  subsequence
$\{(a_{\e_k}, b_{\e_k})\}_{k=1}^\infty$ of $\{(a_\e,
b_\e)\}_{\e>0}$ and  a pair  $(z_0, z)^\top\in
[0,+\infty)\times X'$,  such that $a_{\e_k}\rightarrow
z_0$ in $\dbR$ and $b_{\e_k}\rightarrow z$
weakly$^*$ in $X'$,  as $k\rightarrow \infty$.
Meanwhile, for the subsequence $\{u^{\e_k}\}_{k=1}^\infty$ of  $\{u^{\e}\}_{\e>0}$,  it holds that
$\bd_V(u^{\e_k}, \bar u)\rightarrow 0$ as $k\rightarrow\infty$.

\medskip

Furthermore, by ${\bf (H_1)}$,  for any $(\xi_1,
\xi_2)^\top\in \mathcal V(\bar{u})$, there
exists  a  pair $(\xi^{\e_k}_1, \xi^{\e_k}_2)^\top\in
\mathcal V(u^{\e_k})$  such that
$\xi^{\e_k}_1\rightarrow \xi_1$ in $\dbR$ as $k\to \infty$,  and for sufficiently large
$k\in\mathbb N$,
\begin{equation}\label{new1}
\big|\xi_2-\xi^{\e_k}_2\big|_X\leq \ell\big(\bd_V(\bar  u, u^{\e_k})\big).
\end{equation}
Then, by (\ref{0010})  and (\ref{0011}), we obtain that
$$
z_0\xi_1+\langle z,  \xi_2\rangle_{X', X}\geq 0,
\q  \forall\; (\xi_1, \xi_2)^\top\in \mathcal
V(\bar{u})
$$
and
$$
\langle z, x-\bar{e}\rangle_{X', X}\leq 0,\quad \forall\; x\in E.
$$
This, together with Definition
\ref{llzd3}, implies that $z\in
\cN_E(\bar{e})$.

\medskip

\noindent {\bf Step 3. }    In this step,  we prove  that  the pair  $(z_0,  z)$
is not zero.
If $z_0\neq 0$, the desired result holds naturally.
Otherwise, we assume that $z_0=0$.
Since $a_{\e_k}^2+|b_{\e_k}|_{X'}^2=1$ and
$a_{\e_k}\rightarrow 0$ as $k\rightarrow\infty$, it holds that
$|b_{\e_k}|_{X'}\rightarrow 1$ as $k\rightarrow\infty$.
Hence, for sufficiently large $k$,  $|b_{\e_k}|_{X'}\geq 1/2$.
Furthermore, for the above $u^{\e_k}\in\mathcal{K}$,
by $({\bf H_2})$ and $({\bf H_1})$,
there exists a pair $(\xi^{\e_k}_1, \xi^{\e_k}_2)^\top\in \mathcal V(u^{\e_k})$,  such that for sufficiently large $k$,
\begin{equation}\label{add**1}
|\xi_1^{\e_k}|\leq C \quad\mbox{ and }\quad
|\xi_2-\xi^{\e_k}_2|_{X}\leq \ell\big(\bd_V(\bar  u,  u^{\e_k})\big).
\end{equation}
Here, $C$ is a positive constant independent of $k$,  and  $\ell(\cdot)$ is the modulus of continuity in  $({\bf H_1})$.

 Next,    by (\ref{0011}) and (\ref{0010}),   for any
$x\in E$, $\xi_2\in \pi_{X}(\mathcal{V}(\bar u))$ and the above  $\big(\xi^{\e_k}_1, \xi^{\e_k}_2\big)^\top \in\mathcal V(u^{\e_k})$, it follows that
\begin{eqnarray}\label{add*1}
\begin{array}{rl}
&\displaystyle\langle b_{\e_k}, \xi_2-(x-\bar{e})\rangle_{X', X}\\[3mm]
&\displaystyle=\langle b_{\e_k}, \xi_2\rangle_{X', X} +\langle
b_{\e_k}, f(u^{\e_k})- x\rangle_{X', X}
-\langle b_{\e_k}, f(u^{\e_k})- f(\bar{u}) \rangle_{X', X}\\[3mm]
&\displaystyle\geq \langle b_{\e_k}, \xi^{\e_k}_2\rangle_{X', X}
+\langle b_{\e_k}, \xi_2-\xi^{\e_k}_2\rangle_{X', X}
-\langle b_{\e_k}, f(u^{\e_k})- f(\bar{u})\rangle_{X', X}\\[3mm]
&\displaystyle\geq-\sqrt{{\e_k}}-a_{{\e_k}}\xi^{{\e_k}}_{1} +\langle
b_{\e_k}, \xi_2-\xi^{\e_k}_2\rangle_{X', X}
-\langle b_{\e_k}, f(u^{\e_k})- f(\bar{u})\rangle_{X', X},
\end{array}
\end{eqnarray}
where
$$
a_{\e_k}=\frac{\big[f_0(u^{\e_k})-f_0(\bar{u})+\e_k\big]^+}{\Phi_{\e_k}(u^{\e_k})}\ \mbox{ and }\
b_{\e_k}=\frac{\mbox{dist}\big(f(u^{\e_k}), E\big)\psi_{\e_k}}{\Phi_{\e_k}(u^{\e_k})}.
$$

Combining  (\ref{add*1}) with (\ref{add**1}),  we get that
\begin{eqnarray}\label{add*2}
\begin{array}{rl}
&\displaystyle\langle b_{\e_k}, \xi_2-(x-\bar{e})\rangle_{X', X}\\[3mm]
&\displaystyle\ge -\sqrt{{\e_k}}-C a_{{\e_k}}
-|\xi_2-\xi^{\e_k}_2|_{X} -|f(u^{\e_k})-
f(\bar{u})|_{X}\\[3mm]
&\displaystyle\ge -\sqrt{{\e_k}}-C a_{{\e_k}}
-\ell\big(\bd_V(\bar  u,  u^{\e_k})\big) -|f(u^{\e_k})-
f(\bar{u})|_{X}.
\end{array}
\end{eqnarray}
Let $\delta_{\e_k}=\sqrt{{\e_k}}+C a_{{\e_k}}
+\ell\big(\bd_V(\bar  u,  u^{\e_k})\big) +|f(u^{\e_k})-
f(\bar{u})|_{X}$. Then,  $\delta_{\e_k}>0$  and
$\lim\limits_{  k \to
\infty}\delta_{\e_k}=0$.  By Lemma
\ref{codim converg}  and the finite
codimensionality of
$\pi_X\big(\mathcal{V}(\bar u)\big)-E$ from $({\bf H_3})$,
we obtain that $z\neq 0$.

\medskip

\noindent {\bf Step 4. }  Finally, we prove the enhanced results  (i) and (ii). Clearly,
$z_0\ge 0$. If  $z_0\neq 0$, the conclusion (i) follows
immediately by   dividing  both sides of the inequality (\ref{FJ condition})  by $z_0$.

Now, consider the case where $z\neq 0$. Let $\{u^{\e_k}\}_{k=1}^\infty$,
$\{a_{\e_k}\}_{k=1}^\infty$,
$\{b_{\e_k}\}_{k=1}^\infty$ and
$\{\e_{k}\}_{k=1}^\infty$  be  the
subsequences obtained  in Steps 1 and  2, such that
$
\bd_{V}(\bar u, u^{\e_k})\le \sqrt{\e_k},
$
and as $
k \to \infty$,
$$
\begin{cases}
\ds \ a_{\e_k}\to z_0 \ \mbox{ in }\ \mathbb R, \\
\ns\ds
\ b_{\e_k}\rightarrow z \ \mbox{ weakly}^* \mbox{ in
}\ X', \\
\ns\ds\ \e_{k}\to 0  \ \mbox{ in }\ \mathbb R.
\end{cases}
$$
By the continuity of $f_0$ and
$f$, we have
$
f_{0}(u^{\e_k})\to f_{0}(\bar u)$   in
$\mathbb R$    and   $f(u^{\e_k})\to
f(\bar u)$    in  $X$ as $
k \to \infty$.
Since $z\neq 0$,
 the weakly* lower
semicontinuity of the norm $|\cdot|_{X^{\prime}}$ yields
$$\liminf_{k\to \infty}|b_{\e_k}|_{X^{\prime}}=\liminf_{k\to \infty}\Big| \frac{\mbox{dist}\big(f(u^{\e_k}),E
\big)\psi_{\e_k}}{\Phi_{\e_{k}}(u^{\e_k})}\Big|_{X^{\prime}}
\ge |z|_{X^{\prime}}>0,$$
which implies that $\mbox{dist}\big(f(u^{\e_k}), E\big)> 0$ for  sufficiently large $k\in\mathbb N$. Therefore,
 $f(u^{\e_k})\notin E$  for large $k$. This, combining with $\Phi_{\e_k}(u^{\e_k})\leq
\Phi_{\e_k}(\bar{u})=\e_k$, implies that $f_0(u^{\e_k})< f_0(\bar u)$.

Furthermore, assume that $X$ is a Hilbert space. Then, $\cP_{E}(f(u^{\e_k}))$ is a singleton, and
\begin{eqnarray*}
\lim_{k\to \infty}|\cP_{E}(f(u^{\e_k}))-f(\bar u)|_{X}&\leq& \lim_{k\to \infty}|\cP_{E}(f(u^{\e_k}))-f(u^{\e_k})|_{X} + \lim_{k\to \infty}|f(u^{\e_k})-f(\bar u)|_{X}\\[2mm]
&=&\lim_{k\to \infty}\mbox{dist}(f(u^{\e_k}), E) + \lim_{k\to \infty}|f(u^{\e_k})-f(\bar u)|_{X}= 0.
\end{eqnarray*}
Define
$$
\hat b_{\e_k}\deq\frac{\mbox{dist}\big(f(u^{\e_k}),E\big)}{\Phi_{\e_k}(u^{\e_k})}
\times \frac{
f(u^{\e_k})-\cP_{E}(f(u^{\e_k}))}{|f(u^{\e_k})-\cP_{E}(f(u^{\e_k}))|_{X}}=
\frac{f(u^{\e_k})-\cP_{E}(f(u^{\e_k}))}{\Phi_{\e_{k}}(u^{\e_k})}.
$$
By (\ref{gradent dist}),
 $\hat b_{\e_k}\in X$ is the element
corresponding to $b_{\e_k}\in X'$ via the
Riesz-Fr\'echet isomorphism. Therefore,
%
$\hat
b_{\e_k}\to \hat z \mbox{ weakly in } X,  \mbox{  as }
k\to\infty,$ where $\hat z\in X$ is the element
corresponding to $z\in X'$ via the
Riesz-Fr\'echet isomorphism.
%
Moreover,
\begin{eqnarray*}\label{12.21-eq2}
\begin{array}{ll}\ds
\lim_{k\to \infty}\Big\langle z,
\frac{f(u^{\e_k})
-\cP_{E}(f(u^{\e_k}))}{\Phi_{\e_{k}}(u^{\e_k})}\Big\rangle_{
X',X}\3n&\ds = \lim_{k\to \infty}\Big( \hat z,
\frac{f(u^{\e_k})
-\cP_{E}(f(u^{\e_k}))}{\Phi_{\e_{k}}(u^{\e_k})}\Big)_{
X}\\
\ns&\ds =\lim_{k\to \infty} (\hat z, \hat b_{\e_k})_{ X}=
|\hat z|_{X}^{2}=|z|^2_{X'}>0.
\end{array}
\end{eqnarray*}
 This  demonstrates that for sufficiently large $k\in\dbN$,
$$\Big(\hat z, \frac{f(u^{\e_k})- \cP_{E}(f(u^{\e_k}))}{\Phi_{\e_{k}}(u^{\e_k})}\Big)_{ X}>0.$$
Since $\Phi_{\e_{k}}(u^{\e_k})>0$ for all $k\in\dbN$,
 the conclusion (\ref{enhanced FJ eq3}) immediately follows. This completes  the proof
of Theorem \ref{tt1}.
\endpf

\medskip

In the following discussion,   we compare
  the conclusions of  Theorem \ref{tt1} with
   known results for the optimization problem {\bf (P)}.
\begin{remark}
Assume that $V$ is a reflexive Banach space,  $\mathcal{K}$ is a
nonempty closed convex subset of $V$,
 $f_0$ and $f$ are continuously differentiable, and
 all conditions in Theorem $\ref{tt1}$ hold  for
  $\mathcal{V}(\cdot)=\mbox{\rm Var}_{\mathcal{K}} (f_0, f)(\cdot)$.
 Then,  by Example $\ref{remark2.2}$, for any $u\in \mathcal{K}$,
$$\mbox{\rm Var}_{\mathcal{K}}(f_0, f)(u) =
\Big\{\big(f^{\prime}_{0}(u)v,f^{\prime}(u)v\big)^\top\in \dbR\times X\ \Big|\ v\in \cT_{\mathcal{K}}(u)\cap B_V(0, 1)\Big\},$$
and  the first-order necessary condition in Theorem $\ref{tt1}$ can be expressed as
$$
\langle z_0 f^{\prime}_{0}(\bar u)+f^{\prime}(\bar u)^*z, v\rangle_{V', V}\ge 0, \q \forall\;  v\in \cT_{\mathcal{K}}(\bar u)
\cap B_V(0,1),
$$
or equivalently,
$
-z_0 f^{\prime}_{0}(\bar u)-f^{\prime}(\bar u)^*z\in \cN_{\mathcal{K}} (\bar u).
$

\smallskip

In the special case where  $\mathcal{K}=V$, we have
  $\cT_{\mathcal{K}}(u)=V$ and  $\cN_{\mathcal{K}}(u)=\{0\}$ for any $u\in V$, and  the first-order necessary condition for $\bar u$
reduces to
$z_0f^{\prime}_{0}(\bar u)+f^{\prime}(\bar u)^*z=0.$
This  demonstrates that the first-order necessary condition in Theorem $\ref{tt1}$
extends the classical Fritz John condition $($e.g., $\cite[\mbox{ Page } 153]{BS})$
 to the setting of complete metric spaces.
Furthermore, the enhanced results $(i)$  and $(ii)$  in  Theorem $\ref{tt1}$
 extend  known enhanced Fritz John conditions $($e.g.,
$\cite{Bertsekas2006,guoYezhang2013})$ to  cases
 where the range of  the constraint map $f$ is   infinite dimensional.

\smallskip

Note that in many optimization problems,  such as
 optimal control problems where the control regions are
   genuine  metric spaces lacking a linear structure,
    it is impossible to define standard  directional derivatives or
    even more generalized
     directional derivatives from non-smooth analysis for $f_0$ and $f$.
In such cases,   the variations of $f_0$
and $f$ remain well-defined,
and  the first-order necessary condition
formulated via   these variations  in Theorem  $\ref{tt1}$
offers  broader applicability.
It is also worth noting that precisely characterizing the set-valued map
$\mbox{\rm Var}_{\mathcal{K}} (f_0, f)(\cdot)$ can
 be quite challenging in  concrete problems.
 Therefore, replacing it with a suitable
  set-valued map $\mathcal{V}(\cdot)\subseteq
  \mbox{\rm Var}_{\mathcal{K}}(f_0,f)(\cdot))$
   provides greater flexibility and convenience.
    This approach is particularly  advantageous
    when  constructing  $\mathcal{V}(\cdot)$ is more straightforward
    than fully  determining   $\mbox{\rm Var}_{\mathcal{K}} (f_0, f)(\cdot)$.
    An illustrative example will be
   presented in
 Subsection $\ref{section 2.3}$.
\end{remark}

In what follows,  we present several examples and  remarks
   to further clarify   the conditions
${\bf  (H_1)}$--${\bf  (H_3)}$ of Theorem \ref{tt1}.
The condition ${\bf (H_1)}$  plays a particularly  crucial role.
As demonstrated  below,
a counterexample
shows  that the Fritz John  condition may fail in the absence of ${\bf (H_1)}$,
even when   both ${\bf (H_2)}$ and ${\bf (H_3)}$   are satisfied.
\begin{example}\label{example 2.3}
Let us first recall $\cite[Example\  3.5]{Borgens2020}$. Let
 $V=\mathbb R\times L^2(0, 1)$,
 $X=L^2(0, 1)$, $E=\{0\}$, and
 $\mathcal{K}=\mathbb R\times\big\{ v\in L^2(0, 1)\ \big|\ -1\leq v(t)\leq 1,
 \mbox{ a.e. } t\in (0,1)\big\}$.
 Let  $f_0(\alpha, v)=-\alpha \mbox{  and  } f(\alpha, v)=
 \alpha q-v$ for any $(\alpha, v)\in V,$  where   $q\in L^2(0, 1)\setminus L^\infty(0, 1)$   is a given function.

 Clearly, $(0,0)^\top\in \mathbb R\times L^2(0, 1)$ is the unique feasible point;
   that is, $f(\alpha, v)\in
 E$  if and only if $(\alpha, v)=(0, 0)$.  Therefore, $\bar u=(0, 0)^\top$ is the unique solution to the associated  optimization problem {\bf(P)}.

In this example, $f_0$ and $f$ are continuously differentiable.  By Example $\ref{remark2.2}$, for any $(\alpha, v)^\top\in \mathcal{K}$,
\begin{equation}\label{addre1}
\pi_{\mathbb{R}}(\mbox{\rm Var}_{\mathcal{K}}(f_0,f)(\alpha, v))=
\big\{\!-\beta\in \mathbb{R}\ \big|\ (\beta, w)\in
 \cT_{\mathcal{K}}(\alpha, v),\ |\beta| +|w|_{L^2(0, 1)}\leq 1\big\}
\end{equation}
and
\begin{eqnarray}\label{new11}
\begin{array}{ll}
&\displaystyle\pi_{X}(\mbox{\rm Var}_{\mathcal{K}}(f_0,f)(\alpha, v))
\displaystyle=
\Big\{ \beta q-w\in  L^2(0, 1)\ \Big|\
(\beta, w)\in \cT_{\mathcal{K}}(\alpha, v), \ |\beta|
+|w|_{L^2(0, 1)}\leq 1\Big\}.
\end{array}
\end{eqnarray}
Furthermore, by Definition $\ref{llzd1}$, it is observed that
$\cT_{\mathcal{K}}(0,0) =V=\mathbb{R}\times L^2(0, 1)$,  and therefore,
\begin{equation}\label{new2}
\pi_{X}(\mbox{\rm Var}_{\mathcal{K}}(f_0,f)(\bar u))=\pi_{X}(\mbox{\rm Var}_{\mathcal{K}}(f_0,f)(0, 0))=
\Big\{ \beta q-w\in L^2(0, 1)\ \Big|\ |\beta| +|w|_{L^2(0, 1)}\leq 1\Big\}.
\end{equation}

In Theorem $\ref{tt1}$,
 choose $\mathcal{V}(\cdot)= \mbox{\rm Var}_{\mathcal{K}}(f_0,f)(\cdot)$. Then,
$(\ref{addre1})$ implies the local boundedness of $\pi_\mathbb R(\mbox{\rm Var}_{\mathcal{K}}(f_0,f)(\cdot))$
at $\bar u$, and  $(\ref{new2})$ indicates the finite codimensionality of  $\pi_X(\mbox{\rm Var}_{\mathcal{K}}(f_0,f)(\bar u))-E$ in $X$. Hence,  both ${\bf (H_2)}$   and
${\bf (H_3)}$  are satisfied. However, it was shown in  $\cite{Borgens2020}$
 that the Fritz John condition fails  for this example.

\smallskip

Now,  we  demonstrate  that the condition $(\ref{modulus-contin})$ in  ${\bf (H_1)}$  fails. It suffices to  show that, for sufficiently small  $\varepsilon>0$ and any modulus of continuity $\ell: [0, +\infty) \to [0, +\infty)$ with $\lim\limits_{s\rightarrow 0^+}\ell(s)=0$, there exist  $(\alpha_{\varepsilon},v_\varepsilon)^\top\in B_V(\bar u, \varepsilon)\cap \mathcal{K}$ and $  w_\varepsilon \in \pi_{X}(\mbox{\rm Var}_{\mathcal{K}}(f_0,f)(\bar u))$, such that
\begin{equation}\label{add}
 w_\varepsilon \notin \pi_{X}(\mbox{\rm Var}_{\mathcal{K}}(f_0,f)(\alpha_{\varepsilon},v_\varepsilon))
+\ell\big( |\alpha_{\varepsilon}|  +|v_\varepsilon|_{L^2(0, 1)}\big) B_{L^2(0, 1)}(0, 1).
\end{equation}

To this end, we choose $\alpha_\varepsilon=0$ and
\begin{equation}\label{v_eps}
v_\varepsilon(t)=\left\{
\begin{array}{ll}\ds
0, &\quad t\in [0, 1]\setminus E_\varepsilon,\\
\ns\ds 1, &\quad t\in E_\varepsilon,
\end{array}
\right.
\end{equation}
where $E_{\varepsilon}\subseteq [0, 1]$ is a measurable set with  Lebesgue
measure $\mbox{mes}(E_\varepsilon)=\varepsilon^2$.
It follows that  $(\alpha_\varepsilon,  v_\varepsilon)^\top\in B_V(\bar u, \e)\cap \mathcal{K}$  and
$$|(\alpha_\varepsilon, v_\varepsilon)-\bar u|_V=
|(0, v_\varepsilon)-(0, 0)|_V=|v_\varepsilon|_{L^2(0, 1)}=\varepsilon\rightarrow 0,
\mbox{ as }\varepsilon\rightarrow 0^+.
$$
Furthermore,  by Definition $\ref{llzd1}$,
\begin{eqnarray*}
&&\cT_{\mathcal{K}}(\alpha_\varepsilon, v_\varepsilon)\\
&&=\mathbb R\times \overline{\bigcup_{h\geq 0}\Big\{
h(v-v_\varepsilon)\in L^2(0, 1)\ \Big|\ v(\cdot)\in L^2(0, 1) \mbox{ with } v(t)\in [-1, 1] \Big\}}\\
&&=\mathbb R\times\overline{\Big\{
w\in L^\infty(0, 1)\ \Big|\ w(t)\leq 0, \ a.e.\ t\in E_\varepsilon \Big\}}
=\mathbb R\times \Big\{
w\in L^2(0, 1)\ \Big|\ w(t)\leq 0, \ a.e.\ t\in E_\varepsilon \Big\}.
\end{eqnarray*}
Therefore, by  $(\ref{new11})$ and the  equality above,  we have
\begin{eqnarray*}
&&\pi_{X}(\mbox{\rm Var}_{\mathcal{K}}(f_0,f)(\alpha_\varepsilon, v_\varepsilon))\\[2mm]
&&=\Big\{\beta q-w\in L^2(0, 1)\ \Big|\
|\beta|+|w|_{L^2(0, 1)}\leq 1,\ w(\cdot)\!\in\! L^2(0, 1) \mbox{ with }
w(t)\leq 0, \, a.e.\, t\in E_\varepsilon \Big\}.
\end{eqnarray*}

On the other hand, we  define
\begin{equation}\label{w_eps}
w_\varepsilon(t)=\left\{
\begin{array}{ll}\ds
0, &\quad t\in [0, 1]\setminus E_\varepsilon,\\
\ns\ds -\frac{1}{\varepsilon}, &\quad t\in E_\varepsilon.
\end{array}
\right.
\end{equation}
It is evident  that $|w_\varepsilon|_{L^2(0, 1)}=1$,  and
by $(\ref{new2})$,  $w_\varepsilon\in \pi_{X}(\mbox{\rm Var}_{\mathcal{K}}(f_0,f)(\bar u))$.
Moreover,  for any $\tilde u_\varepsilon\in \pi_{X}(\mbox{\rm Var}_{\mathcal{K}}(f_0,f)(\alpha_\varepsilon, v_\varepsilon))$,
 there exists a  $\beta\in\mathbb R$ and $w\in L^2(0, 1)$, such  that
$|\beta| +|w|_{L^2(0, 1)}\leq 1$, $w(t)\leq 0$ a.e. in $E_\varepsilon$, and $\tilde u_\varepsilon
=\beta q-w$.
Consequently,
\begin{eqnarray*}
&&|w_\varepsilon-\tilde u_\varepsilon|^2_{L^2(0, 1)}
=
\int^1_0 |w_\varepsilon(t)-\tilde u_\varepsilon(t)|^2dt \\
&&=
\int_{E_\varepsilon} |w_\varepsilon(t)-\tilde u_\varepsilon(t)|^2dt
+\int_{[0, 1]\setminus E_\varepsilon} |w_\varepsilon(t)-\tilde u_\varepsilon(t)|^2dt\\
&&\geq \int_{E_\varepsilon} \big|w_\varepsilon(t)-[\beta q(t)-w(t)]\big|^2dt \geq \frac{1}{2} \int_{E_\varepsilon} |w_\varepsilon(t)+ w(t) |^2dt
-\int_{E_\varepsilon} |q(t)|^2dt.
\end{eqnarray*}
Notice that $w(t)\leq 0$ and $w_\varepsilon(t)=-1/\varepsilon$ in $E_\varepsilon$.
Hence,  the above inequality implies
$$
|w_\varepsilon-\tilde u_\varepsilon|^2_{L^2(0, 1)}\geq
\displaystyle\frac{1}{2}\int_{E_\varepsilon} \frac{1}{\varepsilon^2}dt-
\int_{E_\varepsilon}|q(t)|^2dt=\frac{1}{2}-\int_{E_\varepsilon}|q(t)|^2dt.
$$
This means that $|w_\varepsilon-\tilde u_\varepsilon|^2_{L^2(0, 1)}\geq 1/4$ for sufficiently small $\varepsilon$.  However,  for any modulus of continuity $\ell(\cdot)$ in ${\bf (H_1)}$,
we have
$$
\ell(|\alpha_\varepsilon|+|v_\varepsilon|_{L^2(0, 1)})
=\ell(|v_\varepsilon|_{L^2(0, 1)})=\ell(\varepsilon)\rightarrow 0,\quad\mbox{ as }\varepsilon\rightarrow 0^+.
$$
Consequently,  $(\ref{add})$
holds, which implies that   the condition ${\bf (H_1)}$ fails in this example.
\end{example}

Although Example \ref{example 2.3}
shows that
${\bf (H_1)}$ may fail for  general optimization problems,
the following example demonstrates that, under appropriate conditions,
both ${\bf(H_1)}$ and ${\bf(H_2)}$ in Theorem \ref{tt1} are
 readily  satisfied   for $\mathcal{V}(\cdot)=\mbox{Var}_{\mathcal{K}} (f_0, f)(\cdot)$ in
the problem {\bf(P)}.

\begin{example}\label{remark 2.4}
Assume that $V$ is a reflexive Banach space, $\mathcal{K}$  is
a nonempty closed convex subset of $V$, and  $f_0$ and $f$ are continuously
differentiable. Moreover,   assume that
  there exists  a positive constant  $\delta$ and  a   modulus of continuity $\ell_{i}: [0, +\infty) \to [0, +\infty)$ with $\lim\limits_{s\rightarrow 0^+}\ell_{i}(s) = 0$ $(i = 1, 2)$, such that
\begin{equation}\label{2.2-eq1}
\begin{array}{ll}
	\cT_{\mathcal{K}}(\bar u)\cap B_V(0, 1)\\[2mm]
	\subseteq   \cT_{\mathcal{K}}(u)\cap B_V(0, 1)
	+\ell_{1}\big(|\bar  u-u|_V\big) B_V(0, 1), \q \forall\;  u\in B_V(\bar u, \delta)\cap
	\mathcal{K}
\end{array}
\end{equation}
and
\begin{equation}\label{2.2-eq2}
	|f^{\prime}(u_1)-f^{\prime}(u_2)|_{\mathcal{L}(V;X)}\le \ell_{2}\big(|u_1-u_2|_V\big),\q
	\forall\; u_1,u_2\in B_V(\bar u, \delta)\cap \mathcal{K},
\end{equation}
where $\bar u\in\mathcal{K}$ is a solution to the problem {\bf (P)}.
Then it will be shown that  the  conditions  $({\bf H_1})$--$({\bf H_2})$
in Theorem $\ref{tt1}$  hold true  for $\mathcal{V}(\cdot)=\mbox{\rm Var}_{\mathcal{K}}(f_0, f)(\cdot)$.

\smallskip

Indeed,  since $f_0$ and $f$  are continuously
differentiable,
$f_{0}^{\prime}(\cdot)$ and $f^{\prime}(\cdot)$ are locally bounded at $\bar u$,
and
$f_0(\cdot)$ and $f(\cdot)$ are
 locally Lipschitz continuous at
$\bar u$.

For any $(\xi_1,\xi_2)^\top\in \mbox{\rm Var}_{\mathcal{K}}(f_0, f)(\bar u)$, by Example $\ref{remark2.2}$,
there is a $v\in  \cT_{\mathcal{K}}(\bar u)\cap B_V(0, 1)$,  such that
$(\xi_1,\xi_2)=\big(f_{0}^{\prime}(\bar u)v,f^{\prime}(\bar u)v \big). $
For any $\e>0$ and any $u^{\e} \in\mathcal{K}$ satisfying
 $u^{\e} \to \bar u$ in $V$ as $\e\to 0^+$,   the condition  $(\ref{2.2-eq1})$
 implies that  for sufficiently small $\varepsilon$,
$$v\in \cT_{\mathcal{K}}(u^{\e})\cap B_V(0, 1)
+\ell_{1}\big(|\bar  u-u^{\e}|_V\big) B_V(0, 1). $$
This implies that there exists a
 $v^{\e}\in \cT_{\mathcal{K}}(u^{\e})\cap B_V(0, 1)$ such that
\begin{equation}\label{2.2-eq3}
|v-v^{\e}|_{V}\le \ell_{1}\big(|\bar  u-u^{\e}|_V\big).
\end{equation}

Next, define
$$ (\xi_{1}^{\e},\xi_{2}^{\e})
\deq\big(f_{0}^{\prime}(u^{\e})v^{\e},f^{\prime}(u^{\e})v^{\e} \big).$$
Then, $(\xi_{1}^{\e},\xi_{2}^{\e})^\top\in  \mbox{\rm Var}_{\mathcal{K}}(f_0, f)(u^{\e})$.
By the local boundedness of $f^{\prime}_0(\cdot)$ and $f^{\prime}(\cdot)$ at $\bar u$,  there
exists an $L>0$, such that
\begin{equation}\label{LL*}
|f^{\prime}(u^{\e})|_{\mathcal{L}(V;X)}+|f^{\prime}_0(u^{\e})|_{V'}\leq L.
\end{equation}
It follows that
\begin{eqnarray*}
&&|\xi_1-\xi_1^\e|=|f^{\prime}_0(\bar u)v-f^{\prime}_0(u^{\e})v^{\e}|\\[+0.2em]
&&\quad\quad\quad\ \;\,\le|f^{\prime}_0(\bar u)-f^{\prime}_0(u^{\e})|_{V'}|v|_{V}
+|f^{\prime}_0(u^{\e})|_{V'}|v^{\e}-v|_{V}\\[+0.2em]
&&\quad\quad\quad\ \;\,\le|f^{\prime}_0(\bar u)-f^{\prime}_0(u^{\e})|_{V'}
+L|v^{\e}-v|_{V}.
\end{eqnarray*}
This indicates that
$\xi_{1}^{\e}\to \xi_{1}$ in $\dbR$ as $\e\to 0^+$. Meanwhile,
by $(\ref{2.2-eq2})$ and $(\ref{2.2-eq3})$,
\begin{eqnarray*}
&&|\xi_2-\xi_{2}^{\e}|_{X}=|f^{\prime}(\bar u)v-f^{\prime}(u^{\e})v^{\e}|_{X}\\[+0.2em]
&&\quad\quad\quad\quad\;\,\le|f^{\prime}(\bar u)-f^{\prime}(u^{\e})|_{\mathcal{L}(V;X)}|v|_{V}
+|f^{\prime}(u^{\e})|_{\mathcal{L}(V;X)}|v^{\e}-v|_{V}\\[+0.2em]
&&\quad\quad\quad\quad\;\,\le \ell_{2}\big(|\bar  u-u^{\e}|_V\big)+L\ell_1\big(|\bar  u-u^{\e}|_V\big)
\deq \ell\big(|\bar  u-u^{\e}|_V\big).
\end{eqnarray*}
This implies that $({\bf H_1})$ holds for the modulus of continuity $\ell(\cdot)$.
Clearly,
$$\pi_\dbR\big(\mbox{\rm Var}_{\mathcal{K}}(f_0, f)(u)\big)=\big\{ f_{0}^{\prime}(u)v\in X\;    \big|\;  v\in
 \cT_{\mathcal{K}}(u)\cap B_V(0, 1)\big\}.$$
Then, $({\bf H_2})$  follows directly  from  the local boundedness of $f_{0}^{\prime}(\cdot)$
at $\bar u$.

\smallskip

Notice that if  $\bar{u}$ is an interior point of
$\mathcal{K}$,  $(\ref{2.2-eq1})$ is satisfied trivially.
Besides, in the following two typical cases, the condition  $(\ref{2.2-eq1})$
also holds.
\begin{itemize}
\item[$\bullet$] {\bf Case $1$:} Suppose that $\mathcal{K}$ is a polyhedral set  represented by
$$\mathcal{K}=
\big\{u\in V\ \big|\ \langle a_i, u\rangle_{V^{\prime}, V}\le b_i,\ i=1,2,\cdots,m\big\},$$
where $a_i\in V^{\prime}$ and
$b_i\in \dbR$, for $i=1,2,\cdots, m$ and some $m\in \dbN$.

In this situation,
$$\cT_{\mathcal{K}}(u)=\big\{v\in V\ \big|\ \langle a_i, v\rangle_{V^{\prime}, V}\le 0,\ i=I(u)\big\},$$
where $I(u)=\{i\in\{1,2,\cdots,m\}\ |\ \langle a_i, u\rangle_{V^{\prime}, V}= b_i\}$.
Clearly, there exists a $\delta>0$, such that  $I(u) \subseteq I(\bar u)$ for any $u\in B_{V}(\bar u, \delta)\cap \mathcal{K}$. Consequently,
$\cT_{\mathcal{K}}(\bar u)\subseteq \cT_{\mathcal{K}}(u)$ for all $u\in B_V(\bar u, \delta)\cap \mathcal{K},$
and the condition \eqref{2.2-eq1} holds true.

\item[$\bullet$] {\bf Case $2$:} Let $V$ be a Hilbert space, and let  $g: V\to \dbR$
be   continuously differentiable and convex.
Define $$ \mathcal{K}= \big\{u\in V\ \big|\ g(u)\le 0 \big\}.$$
Let $\bar{u}\in  \mathcal{K}$ such that
$g^\prime(\bar{u})\neq 0$ and $g(\bar{u})=0$,
and assume that  there exists a $\delta>0$ and  modulus of continuity
$\ell: [0,\infty)\to [0,\infty)$ so that
$|g^{\prime}(\bar u)-g^{\prime}(u)|_{V'}\le \ell(|\bar u-u|_{V})$ for all $u\in B_{V}(\bar u,\delta).$
Then,
 the condition \eqref{2.2-eq1} is satisfied; the proof of this fact is given in Appendix D.
\end{itemize}
\end{example}

The above Example \ref{remark 2.4} demonstrates  that
 when  $f_0$ and  $f$ are continuously differentiable, and the tangent cone of
   $\mathcal{K}$ exhibits  appropriate continuity near
$\bar u$, the conditions ${\bf(H_1)}$--${\bf(H_2)}$ are satisfied. In Section $\ref{sec-4}$,
 we shall see that  these  requirements  are naturally
 met  for   optimal control problems under
 mild assumptions.

\medskip

It is worth noting
 that in \cite[Corollary  6.4.5, Page  267]{FF},
 the penalty functional method is also adopted to prove the
  Fritz John condition for the problem {\bf (P)}. To ensure the nontriviality of
   the multiplier pair $(z_0, z)^\top\in\mathbb R\times X^{\prime}$, a specific
    constraint qualification (see ${\bf(H_4)}$ below) was introduced.
In the following remark,  we compare   this condition  from  \cite{FF}
with  ${\bf(H_3)}$ of the present  paper.

\begin{remark}\label{propsition 2.1}
In  $\cite[Corollary\ 6.4.5, Page\  267]{FF}$, the
Fritz John   condition for the problem {\bf (P)}
 was studied in the case where  $\mathcal{K}=\mathcal{D}(f^0)\cap\mathcal{D}(f)$.
To guarantee the nontriviality of the multiplier pair $(z_0, z)^\top\in[0, +\infty)\times X^{\prime}$,
the following condition was imposed:
\begin{enumerate}
\item [${\bf  (H_4)}$] There exists a $k_0\in\dbN$,  a $\rho>0$ and  a precompact
sequence of sets
$\{Q_k\}_{k=1}^\infty\subseteq X$ $($i.e.,  any
sequence $\{q_k\}_{k=1}^\infty\subseteq X$   with
$q_k\in Q_k$  has a convergent
subsequence$)$ such that
\begin{equation}
{\rm Int}\Big[\bigcap_{k=k_0}^{\infty}
(\Delta_k+Q_k)\Big]\neq \emptyset,
\end{equation}
where $\Delta_k=\pi_{X}\big(\mbox{\rm Var}_{\mathcal{K}}(f_0,f
)(u_k)\big)-\big(\cT_E(e_k)\cap B_X(0,\rho)\big)$, and
$\{u_k\}_{k=1}^{\infty}\subseteq \mathcal{K}$  and
$\{e_k\}_{k=1}^{\infty}\subseteq E$ satisfy the
conditions in $\cite[Theorem\ 6.4.2, Page\ 265]{FF}$.
\end{enumerate}

Generally, verifying the condition ${\bf(H_4)}$ directly is challenging because it relies on the selection of sequences  $\{u_k\}_{k=1}^{\infty}$, $\{e_k\}_{k=1}^{\infty}$ and $\{Q_k\}_{k=1}^\infty$. In contrast, the condition   ${\bf(H_3)}$ depends  only on the solution $\bar u$.
As  we shall see  in Section $\ref{sec-4}$,
 it can be characterized  by  equivalent {\it a priori} estimates,
 making  it much easier to  verify in  concrete problems.

In what follows,  we demonstrate that,
  under suitable conditions, the condition ${\bf(H_3)}$ indeed  implies
  the condition ${\bf(H_4)}$.

Assume that  $V$ is a reflexive Banach space,  $\mathcal{K}$  is  a nonempty closed convex subset of $V$, $E=\{0\}$,
 $f_0$ and $f$ are continuously differentiable,
 and $\mathcal{V}(\cdot)=\mbox{\rm Var}_{\mathcal{K}}(f_0, f)(\cdot)$.
 Suppose that  there exists a  $\delta>0$ and a modulus of continuity
     $\ell_{1}: [0, +\infty) \to [0, +\infty)$ with $\lim\limits_{s\rightarrow 0^+}\ell_{1}(s)=0$
         such that $(\ref{2.2-eq1})$  is satisfied. Then
      the condition ${\bf  (H_3)}$ implies  the condition ${\bf  (H_4)}$, where
              $\{u_k\}_{k=1}^{\infty}\subseteq\mathcal{K}$ can be any sequence satisfying
     $u_k\to\bar u$ in $V$ as $k\to \infty$.

\smallskip

Indeed, under these conditions, for any
$u\in \mathcal{K}$, it follows from Example $\ref{remark2.2}$ that
\begin{equation}\label{2.2-eqzx}
\pi_{X}(\mbox{\rm Var}_{\mathcal{K}}(f_0, f)(u))=\big\{f^{\prime}(u)v\in X\ \big|\ v\in \cT_{\mathcal{K}}(u)\cap B_V(0, 1)\big\},
\end{equation}
the condition $({\bf H_3})$
with $\mathcal V(\cdot)=\mbox{\rm Var}_{\mathcal{K}}(f_0, f)(\cdot)$
reduces to the requirement that
$\pi_{X}(\mbox{\rm Var}_{\mathcal{K}}(f_0, f)(\bar
u))$
is finite codimensional in $X$.

Clearly,  $\pi_{X}(\mbox{\rm Var}_{\mathcal{K}}(f_0, f)(\bar
u))$ is a closed convex set in $X$.
By  $\cite[Corollary\  3.3, Page\ 140]{2}$ and the  finite codimensionality of this set,
there exists a compact subset $\mathcal{Q}\subseteq X$,
such that
$$\mbox{Int}\big(\pi_{X}(\mbox{Var}_{\mathcal{K}}(f_0, f)(\bar
u))+\mathcal{Q}\big)\neq \emptyset.$$
Therefore, there exists an $x_0\in X$ and a $\gamma_0>0$, such that
$
B_X(x_0, \gamma_0)\subseteq \pi_{X}(\mbox{\rm Var}_{\mathcal{K}}(f_0, f)(\bar
u))+\mathcal{Q}.
$
It follows that for any $x\in B_X(0, \gamma_0)$,
we can find a $\bar v\in \cT_{\mathcal{K}}(\bar u)\cap B_V(0, 1)$ and a $q\in \mathcal{Q}$
 such that
\begin{equation}\label{2.2-eq4}
x_0+x=f'(\bar u)\bar v+q.
\end{equation}

Furthermore, for any $\varepsilon\in (0, \gamma_0)$, there exist two constants $M_1\ge 0$ and $\eta\in(0,  \delta)$,  with $\delta$ being the constant from
$(\ref{2.2-eq1})$, such that for any
$u\in B_V(\bar u, \eta)\cap\mathcal{K}$,
$$ |f'(\bar u)-f'(u)|_{\mathcal{L}(V; X)}\le  \frac{\varepsilon}{2},\q |f'(u)|_{\mathcal{L}(V; X)}\le M_1
\q \mbox{and} \q \ell_{1}(|\bar u-u|_{V})\le \frac{\varepsilon}{2 M_1}.$$
By $(\ref{2.2-eq1})$, for any
$u\in B_V(\bar u, \eta)\cap\mathcal{K}$, there is a
 $ v\in \cT_{\mathcal{K}}( u)\cap B_V(0, 1)$, such that
$|v-\bar v|_{V}\le \ell_{1}(|\bar u-u|_{V}).$
Consequently,
\begin{eqnarray}\label{2.2-eq5}
\begin{array}{rcl}
|f'(\bar u)\bar v-f'(u)v|_X\!\!\!&\leq&\!\!\! |f'(\bar u)-f'(u)|_{\mathcal{L}(V; X)}
|\bar v|_{V}+|f'(u)|_{\mathcal{L}(V; X)}|\bar v-v|_{V}\\
\!\!\!&\leq&\!\!\! |f'(\bar u)-f'(u)|_{\mathcal{L}(V; X)}+M_1\ell_{1}(|\bar u-u|_{V})\leq \varepsilon.
\end{array}
\end{eqnarray}
By $(\ref{2.2-eq4})$, we have
\begin{equation}\label{2.2-eq5+}
x=-x_0+f'(u)v+\big[f'(\bar u)\bar v-f'(u)v\big]+q.
\end{equation}
Furthermore, by \eqref{2.2-eqzx},   \eqref{2.2-eq5} and \eqref{2.2-eq5+}, we obtain
\begin{eqnarray}\label{2.2-eq6}
\begin{array}{rcl}
B_X(0, \gamma_0)\!\!\!&\subseteq&\!\!\! -x_0+\pi_{X}(\mbox{\rm Var}_{\mathcal{K}}(f_0, f)(u))+
B_X(0, \varepsilon)+\mathcal{Q}\\
\!\!\!&\subseteq&\!\!\! -x_0+\pi_{X}(\mbox{\rm Var}_{\mathcal{K}}(f_0, f)(u))+B_X(0, \varepsilon)+\overline{\mbox{co}}\mathcal{Q}.
\end{array}
\end{eqnarray}
Since $\mathcal{Q}$ is compact, $\overline{\mbox{co}} \mathcal{Q}$ is also compact,
which implies that
$\pi_X(\mbox{\rm Var}_{\mathcal{K}}(f_0, f)(u))+\overline{\mbox{co}}\mathcal{Q}$ is closed and convex. This implies that
$
B_X(0, \gamma_0-\varepsilon)\subseteq -x_0+\pi_{X}(\mbox{\rm Var}_{\mathcal{K}}(f_0, f)(u))+\overline{\mbox{co}}\mathcal{Q}.
$
Due to the arbitrariness of $u$ in $B_V(\bar u, \eta)\cap\mathcal{K}$,
we deduce that
$$
B_X(x_0, \gamma_0-\varepsilon)\subseteq \bigcap\limits_{u\in B_V(\bar u, \eta)\cap\mathcal{K}}\Big[\pi_{X}(\mbox{\rm Var}_{\mathcal{K}}(f_0, f)(u))+\overline{\mbox{co}}\mathcal{Q}\Big].
$$
In ${\bf(H_4)}$, by  choosing $\Delta_k=\pi_{X}(\mbox{\rm Var}_{\mathcal{K}}(f_0, f)(u_k))$, $e_k=0$ and
$Q_k=\overline{\mbox{co}}\mathcal{Q}$,
we see that the condition ${\bf(H_4)}$ holds true for a
 sufficiently large $k_0\in\mathbb N$   such  that $|\bar u-u_k|_V\le \eta$ for all
  $k\ge k_0$.
This demonstrates that, under certain   special circumstances,
 the condition  ${\bf(H_4)}$ can be verified  via  the condition  ${\bf(H_3)}$.
\end{remark}

The following result is an immediate  consequence of Theorem \ref{tt1}.
\begin{corollary}\label{coro2.1}
Let $X$ be a reflexive Banach space whose dual
  $X^{\prime}$ is strictly convex, let $\bar u\in\mathcal{K}$ be a solution to the problem {\bf
(P)},  and let  $\bar{e}=f(\bar{u})$.
Assume that the conditions ${\bf (H_1)}$--${\bf (H_2)}$ are satisfied.  If the condition
\begin{equation}\label{cqKKT}
0\in\mbox{\rm  Int}\Big(
\pi_{X}
\big(\mathcal{V}(\bar u)\big)-[E-f(\bar u)]\Big)
\end{equation}
holds true, then there exists a $z\in X'$ with $z\in
\cN_E(\bar{e})$,   such that
\begin{equation}\label{KKT}
\xi_1+\langle z,  \xi_2\rangle_{X', X}\geq 0, \q
\forall\; (\xi_1, \xi_2)^\top\in \mathcal{V}(\bar u).
\end{equation}
Moreover, if $z\neq 0$, then there exists  a sequence $\{
u^k\}_{k=1}^\infty\subseteq \mathcal{K}$,
converging to $\bar u$ as  $k\to \infty$, such
that the conclusions in $(\ref{enhanced FJ eq1})$ hold.
In particular,  when $X$ is a Hilbert space,  $\cP_{E}(f(u^{k}))\to f(\bar u)$ as $k\to \infty$,  and the condition
$(\ref{enhanced FJ eq3})$ holds true.
\end{corollary}
\noindent  {\bf  Proof. }
It is obvious that (\ref{cqKKT}) implies the condition ${\bf(H_3)}$.
Thus, by Theorem \ref{tt1},  there exists  a non-zero pair
$(z_0, z)^\top\in [0,+\infty)\times X'$ with $z\in
\cN_E(\bar{e})$,  such that
$$
z_0\xi_1+\langle z,  \xi_2\rangle_{X', X}\geq
0, \q  \forall\; (\xi_1, \xi_2)^\top\in
\mathcal{V}(\bar u).
$$

Next,
we argue  by contradiction  to prove that $z_0\neq 0$.  If $z_0=0$, then  by
(\ref{add*2}) in  Step 3 of the proof of Theorem  \ref{tt1}, we obtain
$$\langle z, \xi_2-(x-f(\bar u))\rangle_{X', X}
\geq 0,\q \forall\; \xi_2\in
\pi_X\big(\mathcal{V}(\bar u)\big)\mbox{ and } x\in E.
$$
From the condition (\ref{cqKKT}),  this implies that $z=0$. However,  this
contradicts  the fact that  $(z_0,
z)\neq (0,0)$. Hence,  $z_0\neq 0$.

The remaining  conclusions of this corollary follow  immediately from Theorem \ref{tt1}.
\endpf

\medskip

The first-order necessary condition \eqref{KKT} in Corollary \ref{coro2.1} constitutes
 a KKT condition for the problem {\bf(P)}. In the following remark,  we
 compare the constraint qualification
  (\ref{cqKKT}) with the Robinson constraint qualification,
   which is a classical and frequently used condition in the study of
    the KKT conditions for optimization problems.

\begin{remark}
Assume that $V$ is a reflexive Banach space, $\mathcal{K}$  is
 a nonempty closed convex subset of $V$,  $f_0$ and $f$  are  continuously differentiable,  and $\mathcal V(\cdot)=\mbox{\rm Var}_{\mathcal{K}}(f_0, f)(\cdot)$.  Under these conditions,
    \eqref{cqKKT} can equivalently  be expressed  as
\begin{equation}\label{cqKKT2}
0\in {\rm Int}\big(\big\{ f'(\bar u)v- [e-f(\bar
u)]\in X \  \big|\  v\in \cT_{\mathcal{K}}(\bar u)\cap B_V(0, 1) \mbox{ and }
e\in E \big\}\big).
\end{equation}

In the following, we compare the condition $(\ref{cqKKT2})$ with
  the Robinson constraint qualification
\begin{equation}\label{??!!}
0\in {\rm Int}\big(\big\{   f'(\bar u)(u-\bar u)-
[e-f(\bar u)]\in X \  \big|\  u\in \mathcal{K}\mbox{ and
}e\in E \big\}\big).
\end{equation}
First,  we show  that
$(\ref{??!!})$ implies  $(\ref{cqKKT2})$.
Indeed,  \eqref{??!!} indicates   that
$$0\in  {\rm Int}\mathcal{Z}(\bar u),\ \mbox{ where }\ \mathcal{Z}(\bar u)\deq f^\prime(\bar
u)\big(\mathcal{R}_{\mathcal{K}}(\bar u)\big)-
\mathcal{R}_{E}(f(\bar u)).$$
Since $\mathcal{Z}(\bar u)$ is a cone, we have  $X=  \mathcal{Z}(\bar u).$
By $\cite[Theorem\ 2.1]{Kurcyusz1979}$, it follows that
$$
0\in {\rm Int}\big(\big\{ f'(\bar u)\tilde v- \tilde
w\in X \,  \big|\,  \tilde v\in  (\mathcal{K}- \bar u)\cap B_V(0, 1)\mbox{ and } \tilde
w\in (E-f(\bar u))\cap B_X(0, 1)
\big\}\big).
$$
Since $\mathcal{K}- \bar u\subseteq
 \cT_{\mathcal{K}}(\bar u)$ and
  $(E-f(\bar u))\cap B_X(0, 1)\subseteq E-f(\bar u)$, we obtain that \eqref{cqKKT2} holds.

On the other hand, if we further assume that $E=\{0\}$ and
that  $(\ref{2.2-eq1})$  from Example $\ref{remark 2.4}$ is satisfied,
 then \eqref{cqKKT2}  implies  $(\ref{??!!})$. Indeed, in this special case,
  \eqref{cqKKT2}  reduces to
\begin{equation}\label{cqKKT2+}
	0\in {\rm Int}\big(\big\{ f'(\bar u)v\in X \  \big|\  v\in
	\cT_{\mathcal{K}}(\bar u)\cap B_V(0, 1) \big\}\big).
\end{equation}
Similar to the proof of $(\ref{2.2-eq6})$ in Remark $\ref{propsition 2.1}$, we can find positive constants $\gamma_0$, $\varepsilon_0$ and $\eta_0$ with $\varepsilon_0<\gamma_0$, such that for any $u\in B_V(\bar u, \eta_0)\cap\mathcal{K}$,
 $$
 B_X(0, \gamma_0) \subseteq \big\{ f'(u)v\in X \, \big|\, v\in \cT_{\mathcal{K}}(u)\cap B_V(0, 1) \big\}+B_X(0, \varepsilon_0).
 $$
By $\cite[Theorem\ 3.4.5, Page\ 96]{00}$,  the equation
 $
 f(u)=0  \text{ with }   u\in \mathcal{K}
 $
enjoys the following metric regularity at $\bar u\in f^{-1}(0)\cap \mathcal{K}$:
$$
\text{\rm
dist}\big(u,f^{-1}(x)\cap \mathcal{K}\big)
\le M_2|f(u)-x|_{X}, \q \forall\; (u,x)\in \left(B_V(\bar u, \eta_1)
\cap\mathcal{K}\right)\times B_X(0, \gamma_1),
$$
for certain positive constants $\gamma_1$, $\eta_1$ and $M_2$.
Consequently,  as shown in $\cite[Corollary\ 2.2]{Cominetti}$,
 this metric regularity is equivalent to the Robinson constraint qualification
$(\ref{??!!})$, which yields  the desired conclusion.

From the above arguments,
 the condition \eqref{cqKKT}  is generally  weaker than the
  Robinson constraint qualification $(\ref{??!!})$. Furthermore,
  the first-order necessary  condition in Corollary $\ref{coro2.1}$
  serves as a generalization of both
   the classical KKT condition and the enhanced KKT condition
    for constrained optimization problems in metric spaces,  even
    when   the range of the constraint map $f$ is  infinite dimensional.
\end{remark}

To conclude this subsection, inspired by the work in \cite{KK1951},
we present an example illustrating the
effectiveness of the Fritz John condition alongside
 certain
  limitations of the KKT condition. Moreover, this example highlights
  how the enhanced Fritz John condition offers more precise information
  regarding  solutions to the optimization problem.

\begin{example}
Let $V=X\deq\ell^{2}=\Big\{
(u_1,u_2,...)^{\top}\;\Big|\;
u_i\in\dbR,  \ i\in \dbN  \mbox{ and  }\sum\limits_{i=1}^{\infty}
u^{2}_{i}< +\infty \Big\}$. Define a map
$f: V\to X$ as  follows:
$$
f(u)=\big(u_2+(u_1-1)^3,-u_2+(u_1-1)^3, 0, u_4,
u_5,\cds \big)^{\top}, \q \forall\;
u=(u_1,u_2,\cdots)^\top\in \ell^2,
$$
and consider the nonlinear optimization problem:
\begin{equation}\label{ex P}
\mbox{Minimize }  f_0(u)=u_1  \q\mbox{ subject to
}  f(u)=0  \mbox{ for }u=(u_1,u_2,\cdots)^\top\in\ell^2.
\end{equation}
Clearly, $ \mathcal{K}=\ell^2$, $E=\{0\}$, and  for any $u_3\in \dbR$,
\begin{equation}\label{ex eq1}
\bar u=(1,0, u_3, 0, 0,\cds)^{\top}
\end{equation}
solves the optimization problem $(\ref{ex P})$.

Note that  $f_0$ and $f$ are continuously  differentiable on
$\ell^2$,
\begin{eqnarray*}
&&f_{0}^{\prime}(\bar u) v= v_1,\
 f^{\prime}(\bar u)v= (v_2, -v_2, 0, v_4, v_5, ...)^{\top},
 \q \forall\; v=(v_1, v_2,  \cdots)^{\top}\in \ell^2,\\[2mm]
&&
f'(\bar{u})^*z=
(0, z_1-z_2, 0, z_4, z_5, \cdots)^\top,\q\forall\;  z=(z_1, z_2, \cdots)^{\top}\in \ell^2,
\end{eqnarray*}
and for any $u=(u_1, u_2, \cdots)^\top\in\ell^2$ and  $v=(v_1, v_2, \cdots)^\top\in\ell^2$,
$$
f'(u)v=
\big(v_2+3(u_1-1)^2 v_1,  -v_2+3(u_1-1)^2 v_1,  0, v_4, v_5, \cdots\big)^\top.
$$

Choose
$\mathcal{V}(\cdot)=\mbox{\rm Var}_{\mathcal{K}}(f_0, f)(\cdot)$.
 By Example $\ref{remark 2.4}$,   it  is easy to verify  that
$(\ref{2.2-eq2})$ holds and that   $\bar u$ is an
 interior point of $V$.  Therefore,  ${\bf (H_1)}$ and ${\bf (H_2)}$ are satisfied.
Furthermore, it holds that
$$
|\phi|_{\ell^2}\leq |f'(\bar{u})^*\phi|_{\ell^2}, \quad\forall\; \phi=(0, 0, 0, z_4, z_5, \cdots)\in\ell^2.
$$
By Theorem $\ref{t2}$ in the next section,
 $\mbox{\rm Var}_{\mathcal{K}}f(\bar u)$ is
finite codimensional  in $X$, meaning  the condition
 ${\bf (H_3)}$ holds.
Hence,   any  solution
$\bar u$ defined by $(\ref{ex eq1})$ satisfies the  classical
Fritz John condition:
$$
z_0f_{0}^{\prime}(\bar u)+ f^{\prime}(\bar u)^*z=0,
$$
with $z_0=0$ and $z=(\alpha, \alpha, \beta,
0,0,...)^{\top}$, where   $\alpha, \beta\in \dbR$ satisfy
  either  $\alpha\neq 0$, or $\alpha=0$ and $\beta\neq 0$.

\smallskip

Since  $f_0'(\bar u)\notin {\rm Im} (f^{\prime}(\bar u)^*)$, the KKT
condition for the solution $\bar u$ is not satisfied.

\smallskip

Furthermore, for any  $z=(\alpha, \alpha, \beta, 0, 0, \cdots)^{\top}$ with
$\alpha, \beta\in \dbR$,  since  $E=\{0\}$,   we have
$$ \big(z, f(u)-\cP_{E}(f(u))\big)_{\ell^2 }
=(z, f(u))_{\ell^2 }=2\alpha (u_1-1)^3,\quad \forall\; u\in \ell^2.$$
Consequently,  $z_0=0$ and $\hat z=(0, 0, \beta, 0,0,...)^{\top}$ with $\beta\neq 0$ do not satisfy the enhanced Fritz John condition $(\ref{enhanced FJ eq3})$.
In addition,
by $(ii)$ of Theorem $\ref{tt1}$, for sufficiently large $k\in\dbN$,
 $f_0(u^k)-f_0(\bar{u})=u^{k}_{1}-1< 0$ for $u^k=(u^k_1, u^k_2, \cdots)^\top$. Furthermore,
by the enhanced Fritz John condition $(\ref{enhanced FJ eq3})$,
 $ \big(z, f(u^k)-\cP_{E}(f(u^k))\big)_{\ell^2 }
=2\alpha (u_1^k-1)^3> 0.$ Hence, $\alpha<0$, and
 the multiplier pair in the enhanced Fritz John condition
  satisfies:
$$
\begin{cases}\ \ds z_0=0, \\
\ns\ds\  z=(\alpha, \alpha, \beta, 0, 0, \cdots)^{\top}\ \mbox{ with }\alpha<0 \mbox{ and } \beta\in \dbR.
\end{cases}
$$
\end{example}

\subsection{An application of  Theorem \ref{tt1}}\label{section 2.3}

 In this section, we  apply the conclusions of
  Theorem \ref{tt1}
to study an optimal control problem with constraints in infinite dimensions.

Let $T>0$,  let $U$ be  a separable complete metric space,
and let $\mathbb X$ be a Hilbert space. Assume that $A:
\mathcal{D}(A)\subseteq \mathbb X\rightarrow \mathbb X$ generates a
$C_0$-semigroup  on $\mathbb X$. Set

\begin{equation}\label{2.21}
\mathcal{U}_1(0,T)\deq\big\{ u: (0,
T)\rightarrow U\ \big|\  u(\cdot)\mbox{ is a
measurable function} \big\}.
\end{equation}
Clearly,
$\mathcal{U}_1(0,T)$
is a complete metric space equipped  with the Ekeland metric (see
\cite[Proposition 3.10 on Page 145]{2}):
$$
{\mathbf d}_{\mathcal{U}_1(0,T)}(u,v)\deq \mbox{mes}\big(\{t\in (0, T)\ |\
u(t)\neq v(t)\}\big),\q \forall\; u, v\in \mathcal{U}_1(0,T).
$$
Let $y_0\in \mathbb X$, and let  $F: [0, T]\times\mathbb X\times U\rightarrow
\mathbb X$ and $g_0:  [0, T]\times \mathbb X\times
U\rightarrow\dbR$ be given.  Consider the following
controlled evolution equation:
\begin{eqnarray}\label{LL61}
\left\{
\begin{array}{ll}
\ds y_t(t)=Ay(t)+F(t, y(t), u(t)),  \q  t\in [0,
T],&\\  \ns
y(0)=y_0,&
\end{array}\right.
\end{eqnarray}
with a cost functional
\begin{equation}\label{cost fun 2.29}
J\big(u(\cdot)\big)=\displaystyle\int^T_0 g_0(t, y(t),
u(t))dt.
\end{equation}
Here $F$ and $g_0$  satisfy suitable
conditions (stated below) such that for any  control $u(\cdot)\in \mathcal{U}_1(0,T)$,
  the  equation
(\ref{LL61}) admits a unique mild solution.
Furthermore,
 for any $u(\cdot)\in \mathcal{U}_1(0,T)$ and the
corresponding solution $y(\cdot)=y(\cdot; u)$ to (\ref{LL61}),  it holds that $g_0(\cdot, y(\cdot; u),  u(\cdot))\in
L^1(0, T)$.

Let  $E$ be a closed  convex subset of $\mathbb X$. Define the  set
$\mathcal{U}_{ad}=\big\{ u(\cdot)\in \mathcal{U}_1(0,T)\ \big|\  y(T;  u(\cdot))\in
E  \big\}$, and consider the following optimal control problem:
\begin{equation}\label{ocp1}
\text{Minimize } J(u(\cdot))   \quad
\text{subject to } u(\cdot)\in
\mathcal{U}_{ad}.
\end{equation}

Let
$$ V=\mathcal{K}=\mathcal{U}_1(0,T),\q X=\mathbb
X, \q f_0(u)=J(u(\cdot))  \q \mbox{and} \q f(u)=y(T;  u).
$$
Then,  the optimal control  problem \eqref{ocp1} can be cast as a
special case of  the problem  {\bf(P)}.
The optimal control  $\bar u(\cdot)\in \mathcal{U}_{ad}$  is a solution
to the problem  {\bf(P)},  and
the corresponding optimal state is  $\bar y(\cdot)=y(\cdot; \bar u)$.
To derive  first-order necessary conditions  for the optimal
control $\bar u(\cdot)$ via Theorem \ref{tt1},  we need
to construct a   set-valued map
 $\mathcal  V(\cdot)$ satisfying
${\bf (H_1)}$--${\bf (H_3)}$. Note
 that the assumption that   $U$ is  a  separable complete
  metric space is motivated by
 engineering applications,  a
typical example being  $U=\{-1, 1\}$,
which corresponds to the classical ``bang-bang" controls.

Let us begin with the following hypotheses on
$F$ and $g_0$:

\medskip

\noindent ${\bf (A_{1})}$ {\it   $F: [0, T]\times
\mathbb X\times U\rightarrow \mathbb X$
 and $g_0: [0, T]\times \mathbb X\times U\rightarrow \dbR$ are
  strongly measurable with respect
to $t$ in $(0, T)$, and continuously
differentiable
 with respect to  $y$ in $\mathbb X$.
The maps $F(t, \cdot, \cdot)$,  $F_y(t, \cdot, \cdot)$, $g_0(t, \cdot, \cdot)$ and
$g_{0, y}(t, \cdot, \cdot)$ are continuous, and  there exist
positive constants $\delta_1^*$, $L_1$ and $L_2$, such  that
\begin{eqnarray}\label{newequation1}
\begin{array}{ll}
|F_y(t, y, u)|_{\mathcal{L}(\mathbb X)}+|g_{0, y}(t,
y, u)|_{\mathbb X'}+|F(t, 0, u)|_{\mathbb
X}+|g_0(t, 0, u)|\leq L_1,\\
\ns\ds\qq\qq\qq\qq\qq\qq\qq\qq\   \forall\; (t, y, u)\in [0, T]\times
\mathbb X\times U,
\end{array}
\end{eqnarray}
and
\begin{eqnarray}\label{newequation2}
\begin{array}{ll}
|F_y(t, y_1, u) - F_y(t, y_2, u)|_{\mathcal{L}(\mathbb X)} \leq
L_2|y_1-y_2|_{\mathbb X}, \\
\ns\ds\qq\qq   \forall \;   t\in [0, T],  u\in U \mbox{  and }
y_1, y_2\in\mathbb X \mbox{ with }|y_1-y_2|_\mathbb  X\leq \delta_1^*.
\end{array}
\end{eqnarray}
}

Under  the condition ${\bf (A_{1})}$, the  equation
(\ref{LL61}) is  well-posed for any $u(\cdot)\in \mathcal{U}_1(0,T)$. Let
 $y(\cdot; u)$ denote
  the mild solution to (\ref{LL61}) corresponding to the control $u(\cdot)\in \mathcal{U}_1(0,T)$.
Furthermore,
  we consider the linearized  control system:
\begin{eqnarray}\label{LLL62**}\left\{\!\!
\begin{array}{ll}
\ds \xi_t(t)=A\xi(t)+F_y(t, y(t; u),  u(t))\xi(t)+
F(t, y(t; u), v(t))-F(t, y(t; u),  u(t)), \q  t\in [0,
T],&\\  \ns\ds
\xi(0)=0,&
\end{array}\right.
\end{eqnarray}
where $v(\cdot)  \in \mathcal{U}_1(0,T) $ is the control variable.
We  define  a set-valued map by
\begin{eqnarray*}\label{VS}
&&\mathcal V(u)\!=\!\Big\{ \displaystyle\frac{1}{T}(\xi^0 , \xi(T; v))^\top \in \dbR \times\mathbb{X} \  \Big|\  \xi(\cd;v)
\mbox{ is the solution to }
\eqref{LLL62**}\mbox{ and }\nonumber\\
&&\displaystyle\q\q\q\ \   \xi^0=\!\!\int^T_0\!\! \big[\langle g_{0,y}(t,
y(t; u),  u(t)),\xi(t;v)\rangle_{\mathbb X',\mathbb X}+
g_0(t, y(t; u), v(t))-g_0(t, y(t; u),  u(t))\big]dt, \nonumber\\
&&\displaystyle\q\q\q\q\q\q\q\q\q\q\q\q\q\q\q\ \   \mbox{for some }
 v(\cdot)\in \mathcal{U}_1(0,T)  \Big\},\q \forall\; u\in
\mathcal{U}_1(0,T).
\end{eqnarray*}

By the spike variation  technique, one can prove that (e.g., \cite{LLZ})
$$\mathcal V(u)\subseteq
\mbox{Var}_{\mathcal{U}_1(0,T)}(f_0,f )(u), \q \forall \; u\in \mathcal{U}_1(0,T).$$
In particular,  for the optimal pair $(\bar{u}(\cdot), \bar{y}(\cdot))$, letting $\hat\xi(\cdot;v)
$ be the solution to
the equation:
\begin{eqnarray}\label{LLL62}
\left\{
\begin{array}{ll}
\ds \hat\xi_t(t)=A \hat\xi(t)+F_y(t, \bar{y}(t),
\bar{u}(t))\hat\xi(t)+ F(t, \bar{y}(t),
v(t))-F(t, \bar{y}(t), \bar{u}(t)), &
t\in [0, T], \\  \ns\ds \hat\xi(0)=0,
\end{array}\right.
\end{eqnarray}
 we have
\begin{eqnarray}\label{VSbar}
\begin{array}{ll}
&\mathcal V(\bar u)=\Big\{ \displaystyle\frac{1}{T}(\hat\xi^0 , \hat\xi(T;v))^\top \in \dbR \times\mathbb{X} \  \Big|\  \hat\xi(\cd;v)
\mbox{ is the solution to }
\eqref{LLL62}\mbox{ and }\\[3mm]
&\displaystyle\q\q\q\q\hat\xi^0=\int^T_0 \big[\langle g_{0,y}(t,
\bar y(t),  \bar u(t)),\hat\xi(t; v)\rangle_{\mathbb X',\mathbb X}+
g_0(t, \bar y(t), v(t))-g_0(t, \bar y(t),  \bar u(t))\big]dt,   \\[3mm]
&\displaystyle\q\q\q\q\q\q\q\q\q\q\q\q\q\q\q\ \   \mbox{for some }
 v(\cdot)\in \mathcal{U}_1(0,T)  \Big\}.
\end{array}
\end{eqnarray}
Then,
$$
\pi_X\big(\mathcal
V(\bar  u)\big)=\Big\{
\displaystyle\frac{1}{T}\hat\xi(T;v)
\in \mathbb{X}\ \Big|\ \hat\xi(\cd;v)
\mbox{ solves }
\eqref{LLL62} \mbox{ for }
v(\cdot) \in \mathcal{U}_1(0,T)\Big\}
$$
and
\begin{eqnarray*}
&&\3n\3n\3n\pi_\dbR\big(\mathcal V(\bar  u)\big)\!=\!\Big\{  \displaystyle\frac{1}{T}\hat\xi^0\!\in\!\dbR\; \Big|\;
\hat\xi^0\!=\!\int^T_0 \!\!\big[\langle g_{0,y}(t,
\bar{y}(t), \bar{u}(t)),
\hat\xi(t;v)\rangle_{\mathbb X',\mathbb X}\\
&&\qq\qq\q  +
g_0(t, \bar{y}(t), v(t))-g_0(t, \bar{y}(t), \bar{u}(t))\big]dt,\mbox{ for  some }v(\cdot)\in
\mathcal{U}_1(0,T) \Big\}.
\end{eqnarray*}

Now,  we prove that ${\bf (H_1)}$ in
Theorem \ref{tt1} holds.
Let $\xi(\cdot; v)$ and $\hat\xi(\cdot; v)$ be respectively the solutions to (\ref{LLL62**}) and (\ref{LLL62}) corresponding to  $v(\cdot)\in
\mathcal{U}_1(0,T) $. Setting
 $\eta(\cdot)=\xi(\cdot; v)-\hat\xi(\cdot; v)$, we find that $\eta(\cdot)$ satisfies
\begin{eqnarray}\label{new3}
\left\{
\begin{array}{ll}
\ds \eta_t=A\eta+F_y(t, y(t;u), u(t))\eta\!+\!
\left[F_y(t, y(t;u), u(t))-F_y(t, \bar y(t), \bar u(t))\right]\hat\xi\\
\ns\ds \quad   \quad\quad+\left[F(t, y(t;u),  v(t))-F(t, \bar y(t), v(t))\right]\\
\ns\ds   \quad\quad\quad-\left[F(t, y(t;u),  u(t))-F(t, \bar y(t), \bar  u(t))\right], \quad\quad\quad t\in[0, T],\\
\ns\ds \eta(0)=0.
\end{array}
\right.
\end{eqnarray}

In what follows,  we denote by $\mathcal{C}$ a generic
 positive  constant  independent of $u(\cdot)\in \mathcal{U}_1(0,T)$,
 which may vary from line to line.
 When necessary  to emphasize independence,
 we use  $\cC_1$, $\cC_2$, etc.
By  classical estimates for  solutions to evolution equations
 (see \cite[\mbox{Lemma } 4.1, \mbox{Page} 151]{2} for example), we obtain
\begin{equation}\label{new3-1}
|y(\cdot; u)|_{C([0, T]; \mathbb X)}\leq \mathcal{C},
\end{equation}
\begin{equation}\label{new3-2}
|y(\cdot;  u)-\bar y(\cdot)|_{C([0, T]; \mathbb X)}\leq \mathcal{C}_1
\bd_{\mathcal{U}_1(0,T)}(\bar u, u),
\end{equation}
and
\begin{equation}\label{new3-3}
|\hat\xi(\cdot; v)|_{C([0, T]; \mathbb X)}\leq \mathcal{C}.
\end{equation}

By \eqref{newequation1}, for the mild solution  $\eta(\cdot)$ to (\ref{new3}),  we have
\begin{eqnarray}\label{2.52+}
\begin{array}{rl}
&\displaystyle|\eta(t)|_{\mathbb X}\leq
\mathcal{C}\int^t_0 |F_y(s, y(s;u),  u(s))|_{\mathcal{L}(\mathbb X)}
|\eta(s)|_{\mathbb X} ds\\[4mm]
&\displaystyle\quad\quad\quad\q
+\mathcal{C}\int_0^t
|F_y(s,  y(s; u), u(s))-F_y(s, \bar y(s),  u(s))|
_{\mathcal{L}(\mathbb X)}|\hat \xi(s)|_{\mathbb X} ds\\[4mm]
&\displaystyle\quad\quad\quad\q
+\mathcal{C}\int_{\{ s\in[0, t]\ |\ u(s)\neq\bar u(s) \}}
|F_y(s, \bar y(s), u(s))-F_y(s, \bar y(s),  \bar u(s))|
_{\mathcal{L}(\mathbb X)}|\hat \xi(s)|_{\mathbb X} ds\\[4mm]
&\displaystyle\quad\quad\quad\q
+\mathcal{C}\int^t_0 2L_1
|y(s;u)-\bar y(s)|_{\mathbb X}ds\\[4mm]
&\displaystyle\quad\quad\quad\q
+\mathcal{C}\int_{\{ s\in[0, t]\  |\  u(s)\neq\bar u(s) \}}
|F(s, \bar y(s),  u(s))-F(s, \bar y(s),  \bar u(s))|
_{\mathbb X} ds.
\end{array}
\end{eqnarray}
When $\bd_{\mathcal{U}_1(0,T)}(\bar u, u)\leq \displaystyle\frac{\delta_1^*}{\cC_1}$,
it follows from \eqref{new3-2}   that $\sup\limits_{t\in [0,T]}|y(t;  u)-\bar y(t)|_{\mathbb X}\leq \delta_1^*$. Then, by \eqref{newequation2}, for any $s\in [0,t]$,
\begin{equation}\label{2.53+}
|F_y(s,  y(s; u), u(s))-F_y(s, \bar y(s),  u(s))|
_{\mathcal{L}(\mathbb X)}\leq L_2 |y(s; u)-\bar y(s)|_{\mathbb X}.
\end{equation}
Combining (\ref{2.52+})--(\ref{2.53+})  with (\ref{new3-1}), (\ref{new3-3}) and \eqref{newequation1},  we obtain
\begin{eqnarray*}
&&|\eta(t)|_{\mathbb X}\leq
\mathcal{C}\int^t_0 |\eta(s)|_{\mathbb X} ds
+\mathcal{C}
\bd_{\mathcal{U}_1(0,T)}(\bar u, u) \\
&&\quad\quad\quad\q
+\mathcal{C}\int_{\{ s\in[0, t]\  |\  u(s)\neq\bar u(s) \}}
|F(s, \bar y(s),  u(s))-F(s, 0,  u(s))|
_{\mathbb X} ds\\
&&\quad\quad\quad\q
+\mathcal{C}\int_{\{ s\in[0, t]\  |\  u(s)\neq\bar u(s) \}}
|F(s, \bar y(s), \bar u(s))-F(s, 0,  \bar u(s))|
_{\mathbb X} ds\\
&&\quad\quad\quad\q
+\mathcal{C}\int_{\{ s\in[0, t]\  |\  u(s)\neq\bar u(s) \}}
|F(s, 0,  u(s))-F(s, 0,  \bar u(s))|
_{\mathbb X} ds\\
&&\quad\quad\quad\leq \mathcal{C}\int^t_0 |\eta(s)|_{\mathbb X} ds
+\mathcal{C}
\bd_{\mathcal{U}_1(0,T)}(\bar u, u).
\end{eqnarray*}
This, together with   the Gronwall inequality,   implies that when $\bd_{\mathcal{U}_1(0,T)}(\bar u, u)\leq \displaystyle\frac{\delta_1^*}{\cC_1}$,   there exists  a constant  $L^*_1>0$,  such that
\begin{eqnarray}\label{ne1}
\begin{array}{ll}
	&\displaystyle\sup_{t\in [0,T]}|\xi(t; v) -\hat\xi(t; v)|_\mathbb X\\[2mm]
	&\leq L^*_1
	\bd_{\mathcal{U}_1(0,T)}(\bar u, u)\displaystyle=L^*_1 \mbox{mes}\big(
	\big\{  t\in [0, T]\  \big|\  \bar u(t)\neq u(t)   \big\}\big),\q\forall\; v(\cdot)\in \mathcal{U}_1(0,T).
	\end{array}
\end{eqnarray}

For any $\frac{1}{T}(\hat\xi^0, \hat\xi(T;v))^\top
\in\mathcal{V}(\bar{u})$ (recall (\ref{VSbar})), there exists  a control
 $v(\cdot)\in \mathcal{U}_{1}(0,T)$ such that
$\hat\xi(T;v)$ is the terminal  value of  the solution $\hat\xi(\cdot;v)$  to (\ref{LLL62}), with
$$\hat\xi^0=\int^T_0 \big[\langle g_{0,y}(t,
\bar y(t),  \bar u(t)),\hat\xi(t;v)\rangle_{\mathbb X',\mathbb X}+
g_0(t, \bar y(t), v(t))-g_0(t, \bar y(t),  \bar u(t))\big]dt.$$
Furthermore, let $\xi(T;v)$ be the terminal value of the solution $\xi(\cdot;v)$  to (\ref{LLL62**})
 corresponding  to the same  $v(\cdot)$,   and let
$$\xi^0=\int^T_0 \big[\langle g_{0,y}(t,
y(t; u),  u(t)),\xi(t;v)\rangle_{\mathbb X',\mathbb X}+
g_0(t, y(t; u), v(t))-g_0(t, y(t; u),  u(t))\big]dt.$$
Then, $\frac{1}{T}(\xi^0, \xi(T; v))^\top\in \mathcal{V}(u)$.
Moreover, by (\ref{new3-2}), (\ref{ne1}) and Lebesgue's dominated convergence theorem,
 we have
 $$\xi^0\rightarrow\hat\xi^0,\quad\mbox{ as }\  \bd_{\mathcal{U}_1(0,T)}(\bar u, u)\rightarrow 0,$$
  and when $\bd_{\mathcal{U}_1(0,T)}(\bar u, u)\leq \displaystyle\frac{\delta_1^*}{\cC_1}$,
$$|\xi(T;v) -\hat\xi(T;v)|_\mathbb X\leq L^*_1
	\bd_{\mathcal{U}_1(0,T)}(\bar u, u).$$
Consequently, the condition ${\bf (H_1)}$ is satisfied.
By the condition ${\bf (A_{1})}$, it is straightforward  to verify  that  $\pi_\dbR\big(\mathcal V(\cdot)\big)$ is  locally bounded at $\bar{u}$, meaning
 the condition  ${\bf (H_2)}$ in  Theorem \ref{tt1} also holds.

\ss

Let $\mathcal M_1\deq\pi_X\big(\mathcal
V(\bar  u)\big)$.
As a consequence of Theorem \ref{tt1},
we have the following first-order necessary optimality condition for $\bar u$.

\begin{corollary}\label{coro 2.2}
Assume that ${\bf (A_{1})}$ holds. Let $\bar u$ be the optimal control of $(\ref{ocp1})$  and $\bar y$ be the corresponding optimal state. If $\mathcal{M}_1-E$ is finite
codimensional  in $\mathbb{X}$, then  there
exists a non-zero pair $(z_0,
z)^\top\in [0,+\infty)\times \mathbb{X}'$ with
 $z\in \mathcal{N}_{E}(\bar y(T))$ such that
\begin{eqnarray}\label{9.4-eq1}
H(t, \bar{y}(t), \bar{u}(t), z_0,
\psi(t))= \max\limits_{u\in U} H(t,
\bar{y}(t), u, z_0, \psi(t)),\q \mbox{ a.e.
}t\in (0, T),
\end{eqnarray}
where the adjoint state $\psi(\cdot)$ satisfies
\begin{eqnarray}\label{LL65}
	\left\{
	\begin{array}{ll}
		\ds \psi_t(t)=-A^*\psi(t)-F_y(t,
		\bar{y}(t), \bar{u}(t))^*\psi+ z_0
		g_{0, y}(t, \bar{y}(t), \bar{u}(t)),  \q
		t\in [0, T], \\  \ns\ds \psi(T)=-z,
	\end{array}\right.
\end{eqnarray}
with the Hamiltonian defined by
\begin{equation}\label{hamiltion}
H(t, y, u, z_0,  \psi)\deq \langle \psi, F(t, y,
u)\rangle_{\mathbb X',\mathbb X} -z_0 g_0(t, y,
u), \q \forall\;  (t, y, u,z_0, \psi)^\top\in (0,
T)\times\mathbb X\times U\times \dbR\times \mathbb X'.
\end{equation}

Moreover,

\medskip

\noindent  $(i)$ If $z_0\neq 0$,
 the above results hold with $(z_0,z)$ replaced by $(1,\frac{z}{z_0})$.

\noindent  $(ii)$  If $z\neq0$,   there exists a sequence $\{
u^k(\cdot)\}_{k=1}^\infty \subseteq \mathcal{U}_1(0, T)$,    converging
  to $\bar u(\cdot)$ $($in the Ekeland metric$)$  as $k\rightarrow \infty$, such
that
$$
\begin{cases}\ds
y(T; u^k)\notin E  \text{ and }  J(u^k)<J(\bar u)  \hbox{ for sufficiently large $k$},\\
\ns\ds
\lim\limits_{k\rightarrow\infty}\mbox{\rm
dist}\big(y(T; u^k), E\big)=0, \\
\lim\limits_{k\rightarrow\infty}
J(u^k)=J(\bar{u}).
\end{cases}
$$
Furthermore,
$\lim\limits_{k\to \infty} \cP_{E}(y(T; u^k))= y(T; \bar u)  \mbox{ in } \mathbb X$
and
$
\big(\hat z,  y(T; u^k)-\cP_{E}(y(T;
u^k))\big)_{\mathbb X}>0$ for sufficiently large $k$,
where $\hat z \in \mathbb X$ is the element
corresponding to $z \in \mathbb X'$ via the
Riesz-Fr\'echet isomorphism.
\end{corollary}

\noindent  {\bf  Proof. } By Theorem \ref{tt1}, there
exists a non-zero pair $(z_0,
z)^\top\in [0,+\infty)\times \mathbb{X}'$ with $z\in \mathcal{N}_{E}(\bar y(T))$, such that
\begin{equation}\label{LLL63}
z_0 \hat\xi^0+\langle z, \hat\xi(T;v)\rangle_{\mathbb X',\mathbb X}\geq 0,\q \forall\; (\hat\xi^0, \hat\xi(T;v))^\top
\in\mathcal{V}(\bar{u}).
\end{equation}
Then, using the duality relationship between the equation (\ref{LLL62}) and the adjoint equation (\ref{LL65}), we  deduce from (\ref{LLL63}) that
$$
\int_{0}^{T} \big[H(t, \bar{y}(t), \bar{u}(t), z_0,\psi(t))-
H(t,\bar{y}(t), v(t), z_0, \psi(t)) \big] dt \ge 0,\quad \forall\; v(\cdot)\in \mathcal{U}_{1}(0,T),
$$
which   implies (\ref{9.4-eq1}) via the same arguments as in $\cite{2}$. The remaining
assertions  of Corollary \ref{coro 2.2} follow directly from   Theorem \ref{tt1}.
\endpf

\medskip

The first-order necessary condition (\ref{9.4-eq1}), also known as {\it
the Pontryagin maximum principle}, is a classical result (see \cite[Chapter
4, Theorem 1.6]{2}).
 As noted in \cite[Chapter
4]{2},  if the finite codimensionality  condition
fails,  it may be impossible to find a   non-zero
pair $(z_0,z)$ satisfying  \eqref{9.4-eq1}.

\begin{remark}
In $\cite[Page~ 130]{2}$ $($as well as $\cite[Page~ 185]{LLZ})$, only the condition  \eqref{newequation1}  was assumed
to study infinite-dimensional optimal control problems for
evolution equations with end-point constraints. However,
the condition \eqref{newequation1} alone  is insufficient
 to guarantee that
   the condition ${\bf (H_1)}$  of
   Theorem $\ref{tt1}$  holds true, i.e., the uniform continuity of the sets of variations with respect to the controls. Note that
   the condition ${\bf (H_1)}$ plays an  essential role
    in Step $3$ of the proof of Theorem $\ref{tt1}$,
      and its critical  importance is also illustrated
       in Example $\ref{example 2.3}$. Therefore,
       to  ensure the validity of the condition  ${\bf (H_1)}$, both \eqref{newequation1} and \eqref{newequation2} are assumed in the present work.
\end{remark}

\section{Verification of finite codimensionality condition and examples}\label{sec-4}

In this  section,  we  provide
 equivalent criteria for  finite codimensionality condition within
  an  abstract framework and demonstrate how to apply them through
    some specific examples.

Assume that $V$ and $X$ are two reflexive  Banach spaces,
 and ${\bf F}\in \mathcal{L}(V; X).$ Let $\mathbf
M=\{ {\bf F}(v)\in X\  |\   |v|_V\leq 1 \}$. Then it is easy to check that
\begin{equation}\label{llz53}
\ds\overline{\mbox{span}}\,\mathbf M=\overline{\Im({\bf F})} \quad \mbox{and}\quad
\overline{\mbox{co}}\, \mathbf M=\mathbf M.
\end{equation}

We have the following results.

\begin{theorem}\label{t2}
The  following  assertions are equivalent:

\medskip

\noindent  {\bf (1) }  The set $\mathbf M$   is finite codimensional in $X$.

\medskip

\noindent  {\bf (2) } The subspace  $\Im({\bf F})$ is
a  finite  codimensional closed subspace in $X$.
\medskip

\noindent {\bf (3)}  There  exists a  finite
codimensional closed subspace $\wt X$ of  $X'$ and a constant $\cC>0$
$($independent of $\phi\in \wt X)$,  such
that
\begin{equation}\label{221}
|\phi|_{X'}\leq \cC |{\bf F}^*(\phi)|_{V'},\quad\forall\;  \phi\in \wt X.
\end{equation}

\noindent {\bf (4)}   There exists a Banach space $W$, a  compact
operator $\mathcal G: X'\to W$ and a constant $\cC>0$ $($independent of $\phi\in X')$, such that
\begin{equation}\label{5.21-eq1}
|\phi|_{X'}\leq
\cC\big(|{\bf F}^*(\phi)|_{V'} +
|\mathcal G(\phi)|_{W}\big),  \quad\forall\;
\phi\in X'.
\end{equation}

\noindent {\bf (5)}  $\ker({\bf F}^*)$  is  finite-dimensional in $X'$ and $\Im({\bf F}^*)$  is closed in $V'$.

\end{theorem}

\medskip

\noindent {\bf Proof. }  First,  by \cite[Lemma 2.1 and
Remark 4.2]{D} and \cite[Lemma 3]{P},
the equivalence among  {\bf (3)},  {\bf (4)} and {\bf (5)}
can be derived directly by treating  the continuous linear mapping  in those
 known results as the bounded linear
operator ${\bf F}^*$.
Furthermore,
$\Im({\bf F}^*)$  is closed  if and only if $\Im({\bf F})$ is closed.
Thus,  $\ker({\bf F}^*)$  being   finite-dimensional and
$\Im({\bf F}^*)$  being closed is equivalent to
 $\Im({\bf F})$ being  a closed finite  codimensional subspace, which
 establishes
  the equivalence between {\bf (2)} and {\bf (5)}.

\medskip

Next,  we prove that  {\bf (2)} implies {\bf
(1)}.
For any  $k\in\dbN$, let
$
\mathbf N_k=\big\{  {\bf F}(v)\in X\ \big|\   |v|_V\leq k
\big\}.
$
Then $\mathbf N_1=\mathbf M=\overline{\co}  \mathbf M$ and
$ \bigcup\limits_{n\in\dbN} \mathbf N_k
=\Im({\bf F})$. By  {\bf (2)} and  (\ref{llz53}),
$\Im({\bf F})=\overline{\Im({\bf F})}=\overline{\mbox{span}}  \mathbf
M$ is a finite codimensional subspace of $X$. By
the Baire category theorem,  there exists
$\tilde k\in\dbN$ such that $\mathbf N_{\tilde
k}=\overline{\mathbf N_{\tilde k}}$ has at least
one interior point $x_0$ in $\overline{\mbox{span}}\mathbf
M$. Consequently,  $ x_0/\tilde k $ is an interior point
of $\overline{\co} \mathbf M$ in $\overline{\mbox{span}}\mathbf
M$, meaning  $\mathbf M$ is finite codimensional
in $X$.
\medskip

Conversely,
 we prove that    {\bf (1)} implies {\bf (2)}.
Notice that $\overline{\co}\mathbf M\subseteq
\mathbf \Im({\bf F})$.
By Proposition $\ref{add**}$,  we may choose $x_0=0$ in Definition $\ref{llzd6}$.
By
 $({\bf ii})$ of Definition
\ref{llzd6}, $\overline{\mbox{co}}\mathbf M$ has at
least  one  interior point in  the subspace
$\overline{\mbox{span}} \mathbf M=\overline{\Im({\bf F})}$. Hence, $\Im({\bf F})$ also has an
interior  point in $\overline{\Im({\bf F})}$. Since $\Im({\bf F})$ and
$\overline{\Im({\bf F})}$ are two linear subspaces of $X$,
and $\Im({\bf F})$ is dense in $\overline{\Im({\bf F})}$, we have
that $\Im({\bf F})=\overline{\Im({\bf F})}$.  By $({\bf i})$ of
Definition \ref{llzd6}, $\Im({\bf F})=\overline{\mbox{span}}\mathbf M$  is a closed finite
codimensional subspace of   $X$.  Hence, ${\bf (2)}$
holds.
\endpf

\medskip

\medskip

Furthermore,  we  have the following equivalent results.

\begin{theorem}\label{t2***}
The  following  assertions are equivalent:

\medskip

\noindent  {\bf (1) }   $0\in \mbox{Int}\ (\mathbf M)$.

\medskip

\noindent  {\bf (2) } $\Im({\bf F})=X$.
\medskip

\noindent {\bf (3)}  There  exists a   constant $\cC_1>0$ $($independent of $\phi\in X')$,  such
that
\begin{equation}\label{221***}
|\phi|_{X'}\leq \cC_1 |{\bf F}^*(\phi)|_{V'},\quad\forall\;  \phi\in X'.
\end{equation}
\end{theorem}

\noindent {\bf Proof. } It is obvious that {\bf (2)} is equivalent to {\bf (3)}.

\smallskip

Next,  assume that {\bf (1)} is true. Then there exists a $\rho_0>0$ such that
$B_X(0, \rho_0)\subseteq \mathbf M$.  For any $x\in X$,
we can find a sufficiently large $\gamma_0>0$ such that $\frac{1}{\gamma_0} x\in B_X(0, \rho_0)$.
Hence, there exists a $v^*\in V$ satisfying  $|v^*|_V\leq 1$ and
${\bf F}(v^*)=\frac{1}{\gamma_0} x$, which
 indicates that $x={\bf F}(\gamma_0 v^*)\in \Im({\bf F})$, proving that    {\bf (2)} holds.

On the other hand,  if {\bf (2)} holds, then ${\bf F}$ is an open mapping.
Consequently, the set $\mathcal{M}^*\deq
\{ {\bf F}(v)\in X\  |\   |v|_V< 1 \}
$ is open. Since $0\in \mbox{Int}(\mathcal{M}^*)$ and $\mathcal{M}^*\subseteq \mathcal{M}$,
we have that $0\in \mbox{Int}(\mathbf M)$.
\endpf

\medskip

\medskip

In Theorem \ref{tt1},   to obtain
 the Fritz John condition for a solution $\bar u$ to the problem {\bf (P)},
the condition ${\bf  (H_3)}$, i.e., the finite codimensionality of
$\pi_{X}
\big(\mathcal{V}(\bar u)\big)-E$, is required.  It suffices to show that
$\pi_{X}
\big(\mathcal{V}(\bar u)\big)$ is finite codimensional in $X$.
Meanwhile,  in Corollary \ref{coro2.1}, the condition (\ref{cqKKT}), i.e.,
$0\in\mbox{\rm  Int}\big(
\pi_{X}
\big(\mathcal{V}(\bar u)\big)-[E-f(\bar u)]\big)$, guarantees  the KKT
condition for  $\bar u$.
Since $f(\bar u)\in E$, it suffices to require
$0\in\mbox{\rm  Int}\big(
\pi_{X}
\big(\mathcal{V}(\bar u)\big)\big)$.

\smallskip

For the problem {\bf (P)}, in the special case where $\mathcal{K}=V$, $\mathcal{V}(\bar u)=\mbox{Var}_V(f_0, f)(\bar u)$,
and $f_0$ and $f$ are continuously differentiable,  we have
$$
\pi_{X}
\big(\mathcal{V}(\bar u)\big)=\Big\{
f'(\bar u)v\in X\ \Big|\ |v|_V\leq 1 \Big\}.
$$

Thus, by Theorem \ref{t2} and Theorem \ref{t2***},
when the conditions  ${\bf  (H_1)}$ and ${\bf  (H_2)}$ in Theorem \ref{tt1} hold,
both the Fritz John condition and the KKT condition  for  $\bar u$
can be derived by checking suitable {\it a priori} estimates.
Setting ${\bf F}=f'(\bar u)$,  the following results hold true.

\begin{itemize}

\item The finite codimensionality of $\pi_{X}
\big(\mathcal{V}(\bar u)\big)$ is equivalent to the existence of   a Banach space $W$, a  compact
operator $\mathcal G: X'\to W$, and a constant $\cC>0$ $($independent of $\phi\in X')$, such that
\begin{equation}\label{add1*}
\displaystyle |\phi|_{X'}\leq
\cC\big(|f'(\bar u)^*(\phi)|_{V'} +
|\mathcal G(\phi)|_{W}\big),  \quad\forall\;
\phi\in X'.
\end{equation}

\item The condition   that $0\in\mbox{\rm  Int}\big(
\pi_{X}
\big(\mathcal{V}(\bar u)\big)\big)$ is equivalent to the existence of   a constant $\cC_1>0$ $($independent of $\phi\in X')$, such that
\begin{equation}\label{add2*}
\displaystyle |\phi|_{X'}\leq
\cC_1 |f'(\bar u)^*(\phi)|_{V'},  \quad\forall\;
\phi\in X'.
\end{equation}
\end{itemize}

The estimates (\ref{add1*}) and (\ref{add2*})  provide   new
qualification conditions
 for  the Fritz John  and the KKT conditions  in the problem {\bf (P)}.
Below, we apply  these estimates to specific optimization problems.

\subsection{An optimal control problem for an evolution equation with end-point state  constraint}

Let  $ \mathbb V$ and $\mathbb X$ be two
Hilbert spaces,
$\mathcal{U}_2(0,T)= L^2(0, T; \mathbb V)$ and $E=\{y_1\}$ for a  fixed  $y_1\in\mathbb X$.
Given $F: [0, T]\times\mathbb X\times \mathbb V\rightarrow
\mathbb X$ and $g_0: [0, T]\times \mathbb X\times
\mathbb V\rightarrow \dbR$,  consider the optimal control problem:
\begin{equation}\label{ocp1+}
\text{Minimize } J(u(\cdot))  \quad
\text{subject to } u(\cdot)\in \mathcal{U}_2(0,T) \text{ and } y(T;u)=y_1,
\end{equation}
where $J$ is defined by (\ref{cost fun 2.29})  and $y(\cdot;u)$ is the mild solution to (\ref{LL61})
corresponding to  $u(\cdot)\in \mathcal{U}_2(0,T)$.
By choosing
$$V=\mathcal{K}=\mathcal U_2(0, T),\q f_0(u)=J(u(\cdot))\q \mbox{and} \q f(u)=y(T;  u),$$
the optimal control problem (\ref{ocp1+}) becomes    a special case of
the  problem {\bf(P)}.

In the following, we assume that

\smallskip

\noindent ${\bf (A_{2})}$ {\it $F: [0, T]\times\mathbb X\times \mathbb V\rightarrow
\mathbb X$ and $g_0: [0, T]\times \mathbb X\times
\mathbb V\rightarrow \dbR$ are
strongly measurable with respect to $t$ in $(0,
T)$, and continuously  differentiable
with respect to  $(y, u)$ in $\mathbb X\times
\mathbb V$.
Moreover,  there exists a
positive constant   $L_3$,  such  that for any $(t, y, u)\in [0, T]\times
\mathbb X\times \mathbb V$,
\begin{eqnarray}\label{eq 3.4}
\left\{\
\begin{array}{ll}
|F(t, 0, u)|_{\mathbb
X}\leq L_3|u|_{\mathbb V},\\[2mm]
|F_y(t, y, u)|_{\mathcal{L}(\mathbb X)}+|F_u(t, y, u)|_{\mathcal{L}(\mathbb V; \mathbb X)}\leq L_3,\\[2mm]
|g_0(t, y, u)|\leq L_3(1+|y|^2_{\mathbb X}+|u|^2_{\mathbb V}),\\[2mm]
|g_{0, y}(t,
y, u)|_{\mathbb X' }+|g_{0, u}(t,
y, u)|_{\mathbb V' }\leq L_3(1+|y|_{\mathbb X}+|u|_{\mathbb V}).
\end{array}
\right.
\end{eqnarray}
In addition, there exist
positive constants $\delta_2^*$,  $\delta_3^*$,  $L_4$,  $L_5$ and $L_6$, such  that
\begin{eqnarray}\label{eq 3.5}
\begin{array}{ll}
|F_y(t, y_1, u)-F_y(t, y_2, u)|_{\mathcal{L}(\mathbb X)}+
|F_u(t, y_1, u)-F_u(t, y_2, u)|_{\mathcal{L}(\mathbb V; \mathbb X)}
\leq L_4 |y_1-y_2|_{\mathbb X},  \\[2mm]
\q\q\q\q\qq\qq\qq \forall\;   t\in [0, T],  u\in \mathbb V \mbox{  and }
y_1, y_2\in\mathbb X \mbox{ with }|y_1-y_2|_\mathbb  X\leq \delta_2^*,
\end{array}
\end{eqnarray}
\begin{eqnarray}\label{eq 3.6}
\begin{array}{ll}
|F_y(t, y, u_1)-F_y(t, y, u_2)|_{\mathcal{L}(\mathbb X)}+
|F_u(t, y, u_1)-F_u(t, y, u_2)|_{\mathcal{L}(\mathbb V; \mathbb X)}
\leq L_5 |u_1-u_2|_{\mathbb V},\\[2mm]
\qq\qq\qq\qq\qq\qq\qq\qq\qq\qq \forall \;   t\in [0, T],  y\in \mathbb X \mbox{  and }
u_1, u_2\in \mathbb V,
\end{array}
\end{eqnarray}
and
\begin{eqnarray}\label{eq 3.6**}
\begin{array}{ll}
|g_{0, y}(t, y_1, u_1)-g_{0, y}(t, y_2, u_2)|_{\mathbb X'}+
|g_{0, u}(t, y_1, u_1)-g_{0, u}(t, y_2, u_2)|_{\mathbb V'}\\[2mm]
\leq L_6(|y_1-y_2|_{\mathbb X}+|u_1-u_2|_{\mathbb V}),\\[2mm]
\qq\qq\qq\qq\qq\qq \forall \;   t\in [0, T],  y_1, y_2\in \mathbb X \mbox{ with }
|y_1-y_2|_{\mathbb  X}\leq \delta_3^* \mbox{  and }
u_1, u_2\in \mathbb V.
\end{array}
\end{eqnarray}
}

Under the condition   ${\bf (A_{2})}$, \eqref{LL61} is well-posed and $J(u(\cdot))<+\infty$ for any $u(\cdot)\in \mathcal{U}_2(0,T)$.
For any $u(\cdot), v(\cdot)\in \mathcal{U}_2(0,T)$, consider the  following  system:
\begin{eqnarray}\label{LL62}\left\{
\begin{array}{ll}
\ds \xi_t(t)=A\xi(t)+F_y(t,  y(t; u),
u(t))\xi(t)+ F_u(t, y(t; u), u(t)) v(t) ,
\q  t\in [0, T], \\  \ns\ds \xi(0)=0.
\end{array}\right.
\end{eqnarray}
By ${\bf (A_{2})}$ and Example \ref{remark2.2},    for  any $u(\cdot)\in \mathcal U_2(0, T)$,
$$
\begin{array}{ll}
\ds\mbox{Var}_{\mathcal U_2(0, T)}(f_0,f)(u)=\Big\{  (
\xi^0 , \xi(T;v))^\top\in \dbR \times \mathbb X\
\Big|\ \xi(\cd;v) \mbox{ is the solution to }
\eqref{LL62}\mbox{ and } \\[+0.5em]
\qq\qq\qq\qq   \ds\xi^0\!=\!\int^T_0\!\! \Big(\langle
g_{0,y}(t, y(t;u),  u(t)), \xi(t;v)\rangle_{\mathbb{X}',
\mathbb{X}}+
\langle  g_{0,u}(t,  y(t;u), u(t)),   v(t) \rangle_{\mathbb V', \mathbb V}\Big)dt,~~ \\
\ds\q\q\q\q\qq\q\q\q\q\q\q\q\q\q\q\q\q\q\q \ \    \mbox{ for   some }
v(\cdot)\in\mathcal{U}_2(0,T)\mbox{ with } |v|_{L^2(0,  T;  \mathbb
V)}\leq  1 \Big\}.
\end{array}
$$
Then, we have
$$\pi_{\mathbb R}(\mbox{Var}_{\mathcal U_2(0, T)}(f_0,f)(u))=\mbox{Var}_{\mathcal U_2(0, T)}f_{0}(u)\    \mbox{ and }\ \pi_X(\mbox{Var}_{\mathcal U_2(0, T)}(f_0,f)(u))=\mbox{Var}_{\mathcal U_2(0, T)}f(u),$$
with
$$
\begin{array}{ll}\ds
	\mbox{Var}_{\mathcal U_2(0, T)}f_{0}(u)\!=\!\Big\{   \xi^0\!\in\! \mathbb R \
	\Big|\   \xi^0\!=\! \int^T_0 \! \big(\langle
	g_{0,y}(t, y(t;u),  u(t)), \xi(t;v)\rangle_{\mathbb{X}',
		\mathbb{X}}\\
\ns\ds	\qq\qq\qq\qq \qq\qq \qq \ \  +
	\langle  g_{0,u}(t,  y(t;u), u(t)),   v(t) \rangle_{\mathbb V', \mathbb V}\big)dt,\\
\ns\ds	\qq\qq\qq\qq \qq\q  \ \   \xi(\cdot;v )  \mbox{ solves (\ref{LL62}) for }  v(\cdot)\in\mathcal{U}_2(0,T)
	\  \mbox{with }  |v|_{L^2(0,  T;  \mathbb V)}\leq  1 \Big\},
\end{array}
$$
and
$$
\begin{array}{ll}\ds
\mbox{Var}_{\mathcal U_2(0, T)}f(u)=\big\{  \xi(T;v) \!\in\! \mathbb{X}\
\big|\   \xi(\cdot;v )  \mbox{ solves (\ref{LL62}) for }  v(\cdot)\in\mathcal{U}_2(0,T)
 \  \mbox{with }  |v|_{L^2(0,  T;  \mathbb V)}\leq  1 \big\}.
\end{array}
$$

In particular, for the optimal pair $(\bar u, \bar
y)$, letting $ \hat\xi(\cdot;v )$ be the
solution to the equation
\begin{eqnarray}\label{LLL62***}
\left\{
\begin{array}{ll}
\ds \hat\xi_t(t)=A \hat\xi(t)+F_y(t, \bar{y}(t),
\bar{u}(t))\hat\xi(t)+ F_u(t, \bar{y}(t), \bar{u}(t)) v(t), &
t\in [0, T], \\  \ns\ds \hat\xi(0)=0,
\end{array}\right.
\end{eqnarray}
we have
\begin{eqnarray}\label{var f_0 f EOCP}
&&\mbox{Var}_{\mathcal U_2(0, T)}f_{0}(\bar u)\!=\!\Big\{   \hat\xi^0\!\in\! \mathbb R \
	\Big|\    \hat\xi^0\!=\! \int^T_0 \! \big(\langle
	g_{0,y}(t, \bar y(t),  \bar u(t)),  \hat\xi(t;v)\rangle_{\mathbb{X}',
		\mathbb{X}}\nonumber \\
&&\qq\qq\qq\qq\qq\qq\qq   \  \ +
	\langle  g_{0,u}(t, \bar y(t), \bar u(t)),   v(t) \rangle_{\mathbb V', \mathbb V}\big)dt, \\
&&	\qq\qq\qq\qq\qq\q   \  \  \hat\xi(\cdot;v)  \mbox{ solves (\ref{LLL62***}) for }  v(\cdot)\in\mathcal{U}_2(0,T)
	\  \mbox{with }  |v|_{L^2(0,  T;  \mathbb V)}\leq  1 \Big\}\nonumber
\end{eqnarray}
and
\begin{equation}\label{var f EOCP}
\mbox{Var}_{\mathcal U_2(0, T)}f(\bar u)\!=\!\big\{  \hat\xi(T;v) \!\in\! \mathbb{X}\
\big|\   \hat\xi(\cdot;v)  \mbox{ solves (\ref{LLL62***}) for }  v(\cdot)\!\in\!\mathcal{U}_2(0,T)
\mbox{ with }  |v|_{L^2(0,  T;  \mathbb V)}\!\leq\!  1 \big\}.
\end{equation}

Let $\mathcal{V}(u)=\mbox{Var}_{\mathcal U_2(0, T)}(f_0,f)(u)$, and $\xi(\cdot; v)$ and $\hat\xi(\cdot; v)$ be respectively the solutions to (\ref{LL62}) and (\ref{LLL62***}) with respect to $v(\cdot)\in
\mathcal{U}_2(0,T) $.
By (\ref{eq 3.4}), we have the following estimates:
\begin{equation}\label{3.1-eq1}
\begin{cases}\ds
\ |y(\cdot; u)|_{C([0, T]; \mathbb X)}\leq \mathcal{C}(1+|u|_{L^2(0,  T;  \mathbb V)}),\\
\ns\ds
\ |y(\cdot;  u)-\bar y(\cdot)|_{C([0, T]; \mathbb X)}\leq \mathcal{C}
|u-\bar u|_{L^2(0,  T;  \mathbb V)}, \\
\ns\ds
\ |\hat\xi(\cdot; v)|_{C([0, T]; \mathbb X)}\leq \mathcal{C}|v|_{L^2(0,  T;  \mathbb V)}.
\end{cases}
\end{equation}
Further, by (\ref{eq 3.5})--(\ref{eq 3.6**}), for any $u , v\in L^2(0,  T;  \mathbb V)$
with $u$ sufficiently close to $\bar{u}$,
\begin{equation}\label{3.1-eq1**}
|\xi^0-\hat\xi^0|+|\xi(\cdot; v) -\hat\xi(\cdot; v)|_{C([0, T]; \mathbb X)}\leq \mathcal{C}
|u-\bar u|_{L^2(0,  T;  \mathbb V)}|v|_{L^2(0,  T;  \mathbb V)}.
\end{equation}
Here $\mathcal{C}$ is a positive constant depending on
$\delta_{2}^{*}$, $\delta_3^*$, $L_3,L_4$,  $L_5$, $L_6$  and $\bar u$,
 but independent of $u(\cdot)$ and $v(\cdot)$ in  $L^2(0,  T;  \mathbb V)$. Note that  when $|v|_{L^2(0,  T;  \mathbb V)}\le 1$, it follows from (\ref{3.1-eq1}) and (\ref{3.1-eq1**}) that ${\bf(H_1)}$  and ${\bf(H_2)}$ in Theorem
\ref{tt1} hold true. Then we have the following
necessary condition on $\bar u(\cdot)$ by Theorem
\ref{tt1}.

\begin{corollary}\label{corollary2.1}
Assume that ${\bf (A_{2})}$ holds,  let $\bar u$ be the optimal control of $(\ref{ocp1+})$  and let
 $\bar y$ be the corresponding optimal state. If
  $\mbox{\rm Var}_{\mathcal U_2(0, T)}f(\bar u)$ defined by $(\ref{var f EOCP})$  is finite
codimensional  in $\mathbb{X}$,   there
exists a non-zero pair $(z_0,
z)^\top\in [0,+\infty)\times \mathbb{X}'$ such that
\begin{eqnarray}\label{sec3-eq14}
H_u(t, \bar{y}(t), \bar{u}(t), z_0,
\psi(t))=0,\q \mbox{ a.e.
}t\in (0, T),
\end{eqnarray}
where $H$ is the Hamiltonian function defined by $(\ref{hamiltion})$
 and $\psi$ is the solution to $(\ref{LL65})$.

Moreover,

\medskip

\noindent  $(i)$ If $z_0\neq 0$,  the above results hold with $(z_0,z)$ replaced by $(1,\frac{z}{z_0})$.

\noindent  $(ii)$  If $z\neq0$,   there
exists a sequence $\{ u^k(\cdot)\}_{k=1}^\infty
\subseteq \mathcal{U}_2(0, T)$,   converging
to $\bar u(\cdot)$ in $L^2(0,  T;  \mathbb V)$   as
$k\rightarrow \infty$, such that
$$
\begin{cases}\ds
 y(T; u^k)\neq y_1  \text{ and } J(u^k)<J(\bar{u})\  \hbox{for sufficiently large $k$}, \\
\ns\ds
\ \lim\limits_{k\rightarrow\infty} y(T; u^k) = y_1, \\
\ns\ds
\ \lim\limits_{k\rightarrow\infty}
J(u^k)=J(\bar{u}).
\end{cases}
$$
Furthermore,
$$
\big(\hat z,  y(T; u^k)-y_1\big)_{\mathbb X}>0\ \ \hbox{for sufficiently large $k$},
$$
where $\hat z \in \mathbb X$ is the element
corresponding to $z \in \mathbb X'$ via the
Riesz-Fr\'echet isomorphism.
\end{corollary}

\noindent  {\bf  Proof. } By Theorem \ref{tt1}, there
exists a non-zero pair $(z_0,
z)^\top\in [0,+\infty)\times \mathbb{X}'$ with $z\in \mathcal{N}_{E}(\bar y(T))$, such that %
\begin{equation}\label{LLL63+}
z_0 \hat\xi^0+\langle z, \hat\xi(T;v)\rangle_{\mathbb X',\mathbb X}\geq 0,\q \forall\; (\hat\xi^0, \hat\xi(T;v))^\top
\in \mbox{Var}_{\mathcal U_2(0, T)}(f_0,f)(\bar u)
\end{equation}
with $\mbox{Var}_{\mathcal U_2(0, T)}(f_0,f)(\bar u) $ defined by (\ref{var f_0 f EOCP}).
Then, by the duality relationship between the  system (\ref{LLL62***})  and the adjoint system (\ref{LL65}), we  obtain from (\ref{LLL63+}) that
$$
\displaystyle\int^T_0 \big\langle
-F_u(t, \bar{y}(t), \bar{u}(t))^*\psi(t)
+z_0 g_{0, u}(t, \bar{y}(t), \bar{u}(t)),
v(t)\big\rangle_{\mathbb  V', \mathbb V}dt\geq  0,
\q  \forall \;  v(\cdot)\in\mathcal U_2(0, T),
$$
which implies (\ref{sec3-eq14}). The rest of Corollary \ref{corollary2.1} follows from  Theorem \ref{tt1} directly.
\endpf

\begin{remark}
The conditions $\eqref{eq 3.5}$--$\eqref{eq 3.6**}$  ensure  the validity of ${\bf (H_1)}$ in Theorem $\ref{tt1}$. The condition  ${\bf (A_{11})}$ from
 $\cite[Page~185]{LLZ}$ was used to study  optimal control problems for evolution equations with
end-point constraints; however,
  ${\bf (A_{11})}$ is weaker than the conditions in ${\bf (A_{2})}$ of this paper
    and   is  insufficient to guarantee ${\bf (H_1)}$.
\end{remark}

For any $\phi_T\in \mathbb X'$, consider the
following  equation corresponding to
(\ref{LLL62***}):\vspace{-1mm}
\begin{eqnarray}\label{LLL70}\left\{
\begin{array}{ll}
\ds \phi_t(t)=-A^*\phi(t)-F_y(t, \bar{y}(t), \bar{u}(t))^*\phi(t), \q  t\in [0,
T],&\\  \ns
\phi(T)=\phi_T.&
\end{array}\right.\vspace{-2mm}
\end{eqnarray}
The duality relationship between (\ref{LLL62***}) and
(\ref{LLL70})  yields\vspace{-1mm}
$$
f'(\bar{u})^*\phi_T=F_u(t, \bar{y}(t), \bar{u}(t))^*\phi,\q\forall\;  \phi_T\in \mathbb X'.
$$

\begin{corollary}\label{co!}
The following  assertions are equivalent:

\medskip

\noindent  {\bf (1)} The set  $\mbox{\rm Var}_{\mathcal U_2(0, T)}f(\bar{u})$ is  finite codimensional in $\mathbb X$.

\medskip

\noindent  {\bf (2)} There exists a finite
codimensional subspace $\wt{\mathbb X}\subseteq \mathbb X'$, such
that for any $\phi_T\in \wt{\mathbb X}$, the solution $\phi(\cdot)$ to $(\ref{LLL70})$ satisfies
$$
|\phi_T|_{\mathbb X'}\leq \cC |F_u(\cdot,
\bar{y}(\cdot),
\bar{u}(\cdot))^*\phi|_{L^2(0, T; \mathbb V')}.
$$

\medskip
\noindent {\bf (3)}  There  exists  a Banach space $W$ and a compact
operator $\mathcal G$  from  $\mathbb X'$ to  $W$,
such  that  for any $\phi_T\in \mathbb X'$, the solution $\phi(\cdot)$ to $(\ref{LLL70})$
 satisfies
$$
|\phi_T|_{\mathbb X'}\leq \cC\big(|F_u(\cdot,
\bar{y}(\cdot),
\bar{u}(\cdot))^*\phi|_{L^2(0, T; \mathbb V')} +
|\mathcal G\phi_T|_{W}\big).
$$
\end{corollary}

This result has been applied to  analyze
 the finite codimensionality condition for wave and heat equations, as described in \cite{LLZ}. The associated {\it a priori} estimates have been readily  verified.
Moreover, if the linear  wave equations in  \cite{LLZ} are exactly controllable,
the compact operator $\mathcal{G}=0$ in the  estimate to  {\bf (3)}
of Corollary \ref{co!}, yielding
 the KKT condition $(z_0\neq 0)$
in (\ref{sec3-eq14}).

\subsection{An optimal control problem for an elliptic equation   with pointwise state constraint}

This subsection investigates   an elliptic optimal control  problem with a
pointwise  state constraint to illustrate  how different control and state spaces
affect finite codimensionality.

Let $G\subseteq \dbR^n$ ($n\in\dbN$) be a
 bounded domain with a smooth boundary $\Gamma$,
and let  $\hat X$ and $\hat U$ be two Hilbert spaces.  We examine
 two distinct  cases:

\ms

 (I)
$\hat X=\hat U=L^2(G)$;

\ms

(II) $\hat X=H_0^1(G)$ and
$\hat U=H^{-1}(G)$.

\ms

Suppose that $(a^{i
j}(\cdot))_{1\leq i, j\leq n}\in
W^{1, \infty}(G;\dbR^{n\times n})$   is symmetric and
uniformly positive definite in $G$.
Consider the  elliptic
equation: \vspace{-1mm}
\begin{eqnarray}\label{el1}
\left\{
\begin{array}{lll}
\ds-\sum\limits_{i, j=1}^n  (a^{i j}(x)y_{x_i})_{x_j}=F(x)y+u &\mbox{ in }G,\\
\ds y=0 &\mbox{ on }\Gamma,
\end{array}\right.
\end{eqnarray}
where $u(\cdot)\in \hat U$ is the
control variable, $y(\cdot)\in \hat
X$ is the state variable, and $F(\cdot)\in L^\infty(G)$ satisfies  $F(x)\leq 0$ in $G$.
For a  given $y_d(\cdot)\in \hat X$, define \vspace{-2mm}
$$
J\big(u(\cdot)\big)=\displaystyle\frac{1}{2}|y(\cdot)-y_d(\cdot)|^2_{\hat X}+\frac{1}{2}|u(\cdot)|^2_{\hat U}
\quad\mbox{and}\quad
E=\big\{  y(\cdot)\in \hat X\  \big|\  y(x)\geq
0 \mbox{ a.e.  in } G  \big\}.
$$
 In both cases  (I) and  (II),
(\ref{el1}) is well-posed. The optimal control
problem considered in this subsection is
\begin{equation}\label{ocp2}
\text{Minimize } J(u(\cdot)) \quad
\text{subject to }  y(\cdot;
u)\in E\mbox{  and }u(\cdot) \in
\hat U.\vspace{-2mm}
\end{equation}

\smallskip

Assume that $\bar{u}(\cdot)$ is an optimal
control of the problem (\ref{ocp2}) and denote
by $\bar{y}(\cdot)$  the corresponding optimal state.
The optimal control problem in (\ref{ocp2})  becomes   a
special  optimization problem  {\bf(P)}  by  choosing
$$
V=\mathcal{K}=\hat U, \q X=\hat X, \q f_0(u)=J(u(\cdot))\q \mbox{and}\q
f(u)=y(\cdot;  u),
$$
where $y(\cdot;u)$ is the solution to
(\ref{el1}) corresponding  to $u(\cdot)\in \hat U$.

Consider the
following elliptic equation:
\begin{eqnarray}\label{el3+}
\left\{
\begin{array}{lll}
\ds-\sum_{i, j=1}^n  (a^{i j} \xi_{x_i})_{x_j}=F(x)\xi
+v  &\mbox{ in
}  G,\\\ns
\ds \xi =0 &\mbox{ on
} \G.
\end{array}\right.
\end{eqnarray}
We have that, for any $u(\cdot)\in \hat U$ with the corresponding state $y(\cdot; u)$ to (\ref{el1}),
\begin{eqnarray*}
&&\ds\mbox{Var}_{\hat U}(f_0,f)( u ) =\Big\{ (
\xi^0,\xi(\cdot; v))^\top\in \dbR \times
\hat X \  \Big|\
\xi(\cdot; v) \mbox{ solves
(\ref{el3+})  and }\\
&&\qq\qq\qq\qq\ \   \xi^0=
\big(y(\cdot; u)-y_d, \xi(\cdot; v)\big)_{\hat X}+(u, v)_{\hat U},   \mbox{ for   }v(\cdot)\in \hat{U}
\mbox{  with }|v|_{\hat{U}}\leq  1  \Big\},
\end{eqnarray*}
where $(\cdot, \cdot)_{\hat X}$ and $(\cdot, \cdot)_{\hat U}$, respectively,  denote
 the inner products on the Hilbert spaces $\hat X$ and  $\hat U$.
For the optimal control $\bar u(\cdot)$ with the corresponding optimal state $\bar y(\cdot)$,
it is easy to show that
\begin{eqnarray*}
&&\ds\mbox{Var}_{\hat U}(f_0,f)(\bar{u}) =\Big\{ (
\hat\xi^0, \xi(\cdot;v))^\top\in \dbR \times
\hat X \  \Big|\
\xi(\cdot; v) \mbox{ solves
(\ref{el3+})  and }\\
&&\qq\qq\qq\qq \ \   \hat\xi^0=
\big(\bar y-y_d, \xi(\cdot; v)\big)_{\hat X}+(\bar u, v)_{\hat U}  \mbox{ for some  }v\in \hat{U}
\mbox{  with }|v|_{\hat{U}}\leq  1  \Big\},
\end{eqnarray*}
and for any $u\in \hat U$, $f'(u)v=\xi(\cdot; v)$,
where $\xi(\cdot; v)$ is the  solution to (\ref{el3+}) corresponding
to $v(\cdot)$.

In addition, for any $h\in \hat X'=\hat U$, $\hat\f\deq f'(\bar{u})^* h
 \in\hat U'=(\hat X')'$ is the corresponding element for the solution $\f\in \hat X$ to    the following elliptic
equation via the
canonical isomorphism:
\begin{eqnarray}\label{el3}
\left\{
\begin{array}{lll}
\ds-\sum_{i, j=1}^n  (a^{i j}(x)\varphi_{x_i})_{x_j}=F(x)\varphi
+h &\mbox{ in
}  G,\\\ns
\ds \varphi =0 &\mbox{ on
} \G.
\end{array}\right.
\end{eqnarray}

Choose $\mathcal V(\cdot)=\mbox{Var}_{\hat U}(f_0, f)(\cdot)$.
It is easy to show that
the conditions ${\bf (H_1)}$ and ${\bf (H_2)}$ in Theorem \ref{tt1} hold.
Furthermore, consider the
following  elliptic equation:
\begin{eqnarray*}\left\{
\begin{array}{ll}\ds
-\sum_{i, j=1}^n  (a^{i j} \psi_{x_i})_{x_j}=F(x)\psi+z+
\bar w &\mbox{ in
} G,\\\ns
\psi =0 &\mbox{ on
} \Gamma,
\end{array}
\right.
\end{eqnarray*}
where $\bar w\in \hat X'$ is  the corresponding element of $z_0 (\bar y-y_d)\in\hat X$ via the
Riesz-Fr\'echet isomorphism with $z_0\in\mathbb R$.

As a corollary of Theorem \ref{tt1},  if $\mbox{Var}_{\hat U}f(\bar{u})$
is finite codimensional in $\hat X$,
then there exists  a non-zero pair
$(z_0, z)^\top\in [0, \infty)\times
\cN_E(\bar{y})$,  such
that
$\psi\in \hat X$  is the corresponding element of $-z_0\bar{u}\in\hat U=X'$   via the
Riesz-Fr\'echet isomorphism. Note that $z\in \cN_E(\bar{y})$,  if and only if,
$
\langle z, y\rangle_{\hat{X}', \hat{X}}
\leq 0, \ \forall\; y\in E$ and $\langle z, \bar{y}\rangle_{\hat{X}', \hat{X}}
= 0$. Hence, in the case  (I), $z(x)\leq 0$ and $z(x)\bar{y}(x)=0$ a.e. $x\in G$;
in the case (II), $\langle z, \varphi\rangle_{H^{-1}(G), H^1_0(G)}\leq 0$,
$\forall\; \varphi\in H^1_0(G)$ with $\varphi(x)\geq 0$ a.e. $x\in G$, and
$\langle z, \bar{y}\rangle_{H^{-1}(G), H^1_0(G)}= 0$.

Furthermore, by Theorem \ref{t2},  the finite codimensionality  of
the set $\mbox{Var}_{\hat U}f(\bar{u})$ in $\hat X$ reduces  to
the following estimate:
\begin{equation}\label{el2}
|h|_{\hat X'}\leq \cC\big(|\hat\varphi|_{\hat U'}+
|\mathcal{G}h|_{W}\big)=\cC\big(|\varphi|_{\hat X}+
|\mathcal{G}h|_{W}\big),\q\forall \;h\in \hat X',
\end{equation}
where  $\mathcal G$ is a  compact  operator  from
$\hat X'$ to a Banach  space  $W$ and $\varphi(\cdot)$
 is the  solution to  (\ref{el3}). Here and in the sequel of this subsection, $\cC$ denotes a positive constant independent of $h\in \hat X'$.

\medskip

In  the following,  we choose different spaces
$\hat X$ and $\hat U$ to study finite
codimensionality of $\mbox{Var}_{\hat U}f(\bar{u})$, based  on the  {\em a priori}  estimate
(\ref{el2}).

\medskip

\noindent {\bf 1)}  Let $\hat X=\hat U=L^2(G)$. Then
$\mbox{Var}_{\hat U}f(\bar{u})$ is not finite codimensional in $\hat X$.

\medskip

 Indeed, if $\mbox{Var}_{\hat U}f(\bar{u})$ is  finite codimensional in $\hat X$, then by (\ref{el2}), \begin{equation}\label{9.4-eq2}
|h|_{L^2(G)}\leq \cC\big(|\varphi|_{L^2(G)}+
|\mathcal{G}h|_{W}\big),\q \forall \;h\in
L^2(G).
\end{equation}
Write
$$
S=\Big\{ \varphi\in H^2(G)\cap H_0^1(G)\ \Big|\
\varphi \mbox{ is the solution to \eqref{el3}
for some }h\in L^2(G) \Big\}.
$$
By the classical $L^2$ estimate for elliptic
equations and (\ref{9.4-eq2}), $S$ is an
infinite-dimensional subspace of $H^2(G)\cap
H_0^1(G)$ and\vspace{-2mm}
\begin{equation}\label{9.4-eq3}
|\varphi|_{H^2(G)}\leq \cC|h|_{L^2(G)}\leq  \cC\big(|\varphi|_{L^2(G)}+
|\mathcal{G}h|_{W}\big),\q \forall \;h\in
L^2(G).
\end{equation}

Let $\{\varphi_k\}_{k=1}^\infty\subseteq S$ with
$|\varphi_k|_{H^2(G)}=1$.  Set
$$h_{k}=-\sum\limits_{i, j=1}^n  (a^{i j}(x)\varphi_{k,
x_i})_{x_j}-F(x)\varphi_{k}.$$ Then, it follows that
$$
\begin{array}{ll}
|h_{k}|_{L^2(G)}
\leq \cC|
\varphi_{k}|_{H^2(G)}+|F|_{L^\infty(G)}|\varphi_{k}|_{L^2(G)}\leq \cC.\end{array}
$$
This implies that
$\{h_{k}\}_{k=1}^\infty$  has  a
weakly convergent subsequence
$\{h_{k_j}\}_{j=1}^\infty$ in $L^2(G)$; and $\{\varphi_k\}_{k=1}^\infty$ has a strongly convergent subsequence
$\{\varphi_{k_j}\}_{j=1}^\infty$ in $L^2(G)$. Then,
\begin{equation}\label{elliptic}
\{\varphi_{k_j}\}_{j=1}^\infty  \mbox{  is a
Cauchy sequence in } L^2(G).
\end{equation}
Since $\cG$ is a compact
operator,\vspace{-2mm}
\begin{equation}\label{9.4-eq5}
\{\cG h_{k_j}\}_{j=1}^\infty  \mbox{  is a
Cauchy sequence in }W.
\end{equation}
From \eqref{9.4-eq3}--\eqref{9.4-eq5},
$\{\varphi_{k_j}\}_{j=1}^\infty$ is a Cauchy
sequence in $H^2(G)\cap H_0^1(G)$. This
contradicts  the fact that $S$ is an infinite-dimensional
subspace of $H^2(G)\cap H_0^1(G)$. Therefore,
$\mbox{Var}_{\hat U}f(\bar{u})$ is not finite
codimensional in $L^2(G)$.

\medskip

\noindent {\bf 2)} Choose $\hat X=H_0^1(G)$ and $\hat U=H^{-1}(G)$.
 Then the set $\mbox{Var}_{\hat U}f(\bar{u})$ is  finite codimensional in $\hat X$.

\medskip
By the classical $L^2$ theory for
elliptic equations, we
have
\begin{equation}\label{remark!}
|h|_{H^{-1}(G)}\leq \cC |\varphi|_{H^1_0(G)}, \q\forall\; h\in H^{-1}(G),
\end{equation}
where $\varphi\in H^1_0(G)$ is the weak solution to (\ref{el3}) associated to $h\in H^{-1}(G)$.
Then,   (\ref{el2}) holds true with $\mathcal{G}=0$. Consequently, by Theorem \ref{t2},
$\mbox{Var}_{\hat U}f(\bar{u})$ is  finite
codimensional  in $H^1_0(G)$.
\begin{remark}
This example provides  us with a method of choosing  suitable control
 and state spaces to guarantee  finite codimensionality conditions. In
fact, by this method, it is also easy   to show that if  $\hat
X=H^2(G)\cap H_0^1(G)$ and $\hat U=L^2(G)$, then   $\mbox{\rm
Var}_{\hat U}f(\bar{u})$ is  finite codimensional in $\hat X$. But if
   $\hat X=H_0^1(G)$ and $\hat U=L^2(G)$, then
   $\mbox{\rm Var}_{\hat U}f(\bar{u})$ is  not  finite codimensional in $\hat X$.
\end{remark}
\begin{remark}
In the case {\bf 2)} where $\hat X=H_0^1(G)$ and $\hat U=H^{-1}(G)$,
$(\ref{remark!})$ implies that the compact operator $\mathcal{G}=0$ in the corresponding estimate to  {\bf (3)}
of Corollary  $\ref{co!}$.
At this moment, we can obtain the KKT condition for the optimal control, i.e.,  $z_0\neq 0$.
\end{remark}

\subsection{An optimal control problem for  stochastic differential equations  with  end-point  state constraint}

\subsubsection{Formulation of the problem and a necessary
condition for optimal controls}

Let $\left(\Omega, \mathcal{F}, \{\mathcal{F}_t\}_{t\geq 0}, \mathbb
P\right)$ be a complete filtered probability space, on which a
one-dimensional standard Brownian motion $\{B(t)\}_{t\geq 0}$ is
defined, such that $\{\mathcal{F}_t\}_{t\geq 0}$ is the natural
filtration generated by $B(\cdot)$, augmented by all the $\mathbb
P$-null sets in $\mathcal{F}$. Denote by $\mathbb F$ the progressive
$\si$-field (in $[0, T ]\times \O$) with respect to
$\{\mathcal{F}_t\}_{t\geq 0}$. For any $\alpha,\beta\geq 1$ and
$t\in [0,T]$,   denote by $L_{\mathcal{F}_{t}}^{\beta}(\Omega;
\dbR^n)$ the set of all $\dbR^n$-valued,
$\mathcal{F}_{t}$-measurable random variables $\zeta$ with $\mathbb
E|\zeta|_{\dbR^n}^{\beta}<\infty$, by $L_{\mathbb{F}}^\beta(\Omega; C([0,T];\dbR^n))$ the space of $\dbR^n$-valued, $\mathbb{F}$-progressively measurable continuous stochastic processes $\eta$ such that $\mathbb{E}\big(\sup_{0\leq t\leq T}|\eta(t)|_{\dbR^n}^\beta\big)<\infty$,
and by $L_{\mathbb F}^{\beta}(\Omega;
L^{\alpha}(0,T; \dbR^n))$ the set of all $\dbR^n$-valued, $\mathbb
F$-progressively measurable stochastic processes $\eta$ with $
\mathbb{E}~\big(\int_{0}^{T}|\eta(t,\omega)|_{\dbR^n}^{\alpha}dt\big)
^{\beta/\alpha} <\infty$. When $\alpha=\beta$, it is simply
denoted by $L_{\mathbb{F}}^{\beta}(0,T;\dbR^n)$. As usual, when the
context is clear, we omit the argument $\omega\in \O$ in the corresponding
variables and processes.

\ss

Consider the  controlled stochastic
differential equation:
\begin{eqnarray}\label{1LLZ}
\left\{
\begin{array}{ll}
dy(t)=\Big[A_1(t)y(t)+C_1(t)u(t)+\eta_1(t)\Big]dt&\\ \ns
\quad\quad\quad\quad+\Big[A_2(t)y(t)+C_2(t)u(t)+\eta_2(t)\Big]dB(t), \q  t\in [0, T],&\\ \ns
y(0)=y_0,&
\end{array}
\right.
\end{eqnarray}
where $u(\cdot)$ is the control variable and
$y(\cdot)$ is the state variable. Let
$$
J(u(\cdot))=\mathbb{E}\displaystyle\int^T_0 g_0(t, y(t), u(t))dt,
\quad\forall\; u(\cdot)\in L^2_\mathbb{F}(0, T; \dbR^m).
$$

Assume that the
following condition holds.

\smallskip

\noindent ${\bf (A_{3})}$
{\it For $i=1, 2$, $\eta_i(\cdot)\in L^\infty_{\mathbb{F}}(0, T; \mathbb{R}^n)$,
$A_i(\cdot)\in L^\infty_{\mathbb{F}}(0, T; \mathbb{R}^{n\times n})$ and
$C_i(\cdot)\in L^\infty_{\mathbb{F}}(0, T; \mathbb{R}^{n\times m})$;
and for $g_0:
[0, T]\times \dbR^n\times
\dbR^m\times\Omega\rightarrow \dbR$,
$g_0(\cd, y, u,\cd): [0, T]
\times\Omega\rightarrow \dbR$ are $\mathbb F$-progressively
measurable for any $y\in\mathbb{R}^n$ and $u\in\mathbb{R}^m$,
 and
$g_0(t,\cd,\cd,\omega)$  are  continuously
differentiable for a.e. $(t,\o)\in [0,T]\times\O$.  Moreover, there is  a  positive constant $L_7$,
 such that for    a.e.
$(t, \omega)\in[0,T]\times \Omega, \; \forall\;  y\in\mathbb R^n \mbox{ and }
  u\in\mathbb R^m$,
\begin{eqnarray*}
\left\{
\begin{array}{ll}
|g_0(t, y, u)|\leq L_7(1+|y|^2_{\mathbb R^n}+|u|^2_{\mathbb R^m}),\\[2mm]
|g_{0, y}(t, y, u)|_{\mathbb R^n}+|g_{0, u}(t, y, u)|_{\mathbb R^m}\leq L_7(1+|y|_{\mathbb R^n}+|u|_{\mathbb R^m}),
\end{array}
\right.
\end{eqnarray*}
and
\begin{eqnarray*}\begin{array}{ll}
|g_{0, y}(t, y_1, u_1)-g_{0, y}(t, y_2, u_2)|_{\mathbb R^n}+
|g_{0, u}(t, y_1, u_1)-g_{0, u}(t, y_2, u_2)|_{\mathbb R^m}\\[2mm]
\leq L_{7}(|y_1-y_2|_{\mathbb R^n}+|u_1-u_2|_{\mathbb R^m}),\\[2mm]
\q\q\qq\qq \mbox{a.e.} \  (t,\omega)\in[0,T]\times \Omega,\quad
 \forall\; y_1, y_2\in \mathbb R^n  \mbox{  and }
 u_1, u_2\in \mathbb R^m.
\end{array}
\end{eqnarray*}
}

Under the condition ${\bf (A_{3})}$, the equation \eqref{1LLZ} is well-posed (see, e.g., \cite[Section
3.1]{LZ} and \cite{YZ}) and  $J(u(\cdot))<+\infty$  for any $u(\cdot)  \in L^2_\mathbb{F}(0, T;
\dbR^m)$.
For a given $y_T\in
L^2_{\mathcal{F}_T}(\Omega; \dbR^n)$, consider the following stochastic optimal control problem:
\begin{equation}\label{ocp3}
\text{Minimize } J(u(\cdot)) \quad \text{subject to } y(T; u)=y_T
\mbox{ and }u(\cdot)\in L^2_\mathbb{F}(0, T; \dbR^m).
\end{equation}
Suppose that $\bar{u}(\cdot)$
solves the optimal
control problem in (\ref{ocp3}) and
 $\bar{y}(\cdot)$ is the corresponding state of  (\ref{1LLZ}).
Let us give a necessary
condition of $\bar{u}$ by
 Theorem \ref{tt1}.

\smallskip

In the optimization  problem  {\bf(P)}, choose
$$
V=\mathcal{K}=L^2_\mathbb{F}(0, T; \dbR^m), \ \
X=L^2_{\mathcal{F}_T}(\Omega; \dbR^n),\ \
E=\{y_T\},\ \
f_0(u)=J(u(\cdot))\ \mbox{ and }\
f(u)=y(T; u),
$$
where $y(\cdot)=y(\cdot; u)$ is the solution to (\ref{1LLZ})
corresponding  to $u(\cdot)\in L^2_\mathbb{F}(0, T; \dbR^m)$.

For any $u, v\in L^2_\mathbb{F}(0, T;
\dbR^m)$,  consider the stochastic differential
equation:
\begin{equation}\label{LLZ2+}
\left\{
\begin{array}{ll}
d\hat\xi(t)=[A_1(t)\hat\xi(t)+C_1(t) v(t)]dt
+[A_2(t)\hat\xi(t)+C_2(t) v(t)]dB(t),\q t\in [0, T], \\
\ns\ds \hat\xi(0)=0.
\end{array}
\right.
\end{equation}
Under the condition ${\bf (A_{3})}$,
the equation (\ref{LLZ2+}) is also well-posed.  Moreover, we have
\begin{eqnarray*}
&&\3n\3n\3n\ds\mbox{Var}_{L^2_\mathbb{F}(0, T; \dbR^m)}(f_0,f)(u ) =\Big\{ (
\xi^0,\hat\xi(T;v))^\top\in \dbR \times L^2_{\mathcal{F}_T}(\Omega;
\dbR^n) \ \Big|\
\hat\xi(\cdot;v)\mbox{ solves }(\ref{LLZ2+})\mbox{ and }\\
&&\q\q\q\q\ \qq\qq\qq   \xi^0\!=\mathbb E\!\!\int^T_0\!\!
\big[g_{0, y}(t, y(t;u),
u(t))^\top\hat\xi(t;v)\!+\! g_{0, u}(t, y(t;u),
u(t))^\top v(t)\big]dt, \\
&&\q\q\q\q\q\q\q\q\ \qq\qq\qq \mbox{for } v(\cdot)\in
L^2_\mathbb{F}(0, T; \dbR^m) \mbox{ with }|v|_{L^2_\mathbb{F}(0, T;
\dbR^m)}\!\leq\!  1 \Big\},
\end{eqnarray*}
and $f^{\prime}(u)v=\hat\xi(T;v)$.
Further, for simplicity of notation,  we write
$$
  A_3(t)\deq g_{0, y}(t, \bar{y}(t),
\bar{u}(t))\quad \mbox{ and }\quad  C_3(t)\deq g_{0, u}(t,
\bar{y}(t), \bar{u}(t)).
$$
Then,
\begin{eqnarray}\label{var f_0f SOCP}
\begin{array}{rl}
&\displaystyle\3n\3n\3n\ds\mbox{Var}_{L^2_\mathbb{F}(0, T; \dbR^m)}(f_0,f)(\bar{u}) =\Big\{ (\hat
\xi^0,\hat\xi(T; v))^\top\in \dbR \times L^2_{\mathcal{F}_T}(\Omega;
\dbR^n) \ \Big|\
\hat\xi(\cdot;v)\mbox{ solves }(\ref{LLZ2+})\mbox{ and }\\
&\displaystyle\q\q\q\q\q\qq\q\qq\      \hat\xi^0\!=\mathbb E\!\int^T_0\!
\big[A_3(t)^\top\hat\xi(t;v)\!+\! C_3(t)^\top v(t)\big]dt,\\
&\displaystyle\q\q\q\q\q\qq\qq\q\q\q\q\q \mbox{ for } v(\cdot)\in
L^2_\mathbb{F}(0, T; \dbR^m)\mbox{  with }|v|_{L^2_\mathbb{F}(0, T;
\dbR^m)}\!\leq\!  1 \Big\},
\end{array}
\end{eqnarray}
and
\begin{eqnarray}\label{varf for socp}
\begin{array}{rl}
&\displaystyle\mbox{Var}_{L^2_\mathbb{F}(0, T; \dbR^m)}f(\bar{u})\!=\!\Big\{  \hat\xi(T;v)\!\in\!
 L^2_{\mathcal{F}_T}(\Omega; \dbR^n)\, \big|\,
\hat\xi(\cdot;v)\mbox{ solves (\ref{LLZ2+})},\\
&\displaystyle\q\q\q\q\q\q\q\q\q\q\q\q\q\mbox{ for }v(\cdot)\!\in\!
L^2_\mathbb{F}(0, T; \dbR^m)\mbox{ with }
|v|_{L^2_\mathbb{F}(0, T; \dbR^m)}\!\leq\! 1
\Big\}.
\end{array}
\end{eqnarray}

Choose $\mathcal
V(\cdot)=\mbox{Var}_{L^2_\mathbb{F}(0, T; \dbR^m)}(f_0, f)(\cdot)$.
By the condition   ${\bf (A_3)}$, it is easy to show that
the conditions ${\bf (H_1)}$ and ${\bf (H_2)}$ in Theorem \ref{tt1} hold.
As  a
corollary of Theorem \ref{tt1},
 one gets the following necessary
condition for  the optimal  control $\bar u(\cdot)$.

\begin{corollary}\label{cc}
Assume that ${\bf (A_{3})}$ holds. Let $\bar u$ be the optimal control of $(\ref{ocp3})$ and $\bar y$ be the corresponding optimal state. If $\mbox{\rm Var}_{L^2_\mathbb{F}(0, T; \dbR^m)}f(\bar{u})$ defined as $(\ref{varf for socp})$ is finite
codimensional  in $L^2_{\mathcal{F}_T}(\Omega;
\dbR^n)$, then  there exists  a non-zero pair $(z_0,
z)^\top\in[0,+\infty)\times
L^2_{\mathcal{F}_T}(\Omega; \dbR^n)$,
such that
\begin{eqnarray}\label{3.3 stochastic mp}
\mathcal{H}_{u}(t, \bar{y}(t), \bar{u}(t), z_0,
\psi(t),\Psi(t))=0,\q \mbox{\rm
a.e. }(t,\omega)\in[0,T]\times \Omega,
\end{eqnarray}
where
\begin{eqnarray}\label{hamiltion of socp}
\begin{array}{rl}
&\displaystyle\mathcal{H}(t, y, u, z_0,  \psi,\Psi)\deq (\psi, A_1(t)y+C_1(t)u+\eta_1(t))_{\dbR^n} \\[2mm]
&\displaystyle\qq\qq\qq\qq\q+(\Psi, A_2(t)y+C_2(t)u+\eta_2(t))_{\dbR^n}-z_0 g_0(t, y,
u), \\[2mm]
&\displaystyle\qq\qq\qq\qq\q\forall\;  (t, y, u,z_0, \psi,\Psi, \omega)^\top\in (0,
T)\times \dbR^n\times \dbR^m\times \dbR \times \dbR^n\times \dbR^n\times \O,
\end{array}
\end{eqnarray}
and  $(\psi,\Psi)$ solves  the   backward
stochastic differential  equation:
\begin{eqnarray}\label{LLZ4}\left\{
	\begin{array}{ll}
		d\psi(t)=[z_0 A_3(t)-A_1(t)^\top
		\psi(t)-A_2(t)^\top \Psi(t)]dt+\Psi(t) dB(t), & t\in [0, T], \\
		\ns \psi(T)=-z.
	\end{array}
	\right.
\end{eqnarray}
Moreover,

\medskip

\noindent  $(i)$ If $z_0\neq 0$,  the above results hold with $(z_0,z)$ replaced by $(1,\frac{z}{z_0})$.

\noindent  $(ii)$  If $z\neq0$,   there
exists a sequence $\{ u^k(\cdot)\}_{k=1}^\infty
\subseteq L^2_\mathbb{F}(0, T; \dbR^m)$,   converging
to $\bar u(\cdot)$ in $L^2_\mathbb{F}(0, T; \dbR^m)$   as
$k\rightarrow \infty$, such that
$$
\begin{cases}\ds
 \ y(T; u^k)\neq y_T  \text{ and } J(u^k)<J(\bar{u})\ \text{ for sufficiently large $k$}, \\
\ns\ds
\ \lim\limits_{k\rightarrow\infty} y(T; u^k) = y_T, \\
\ns\ds
\ \lim\limits_{k\rightarrow\infty}
J(u^k)=J(\bar{u}).
\end{cases}
$$
Furthermore,
$$
\big(z,  y(T; u^k)-y_T\big)_{L^2_{\mathcal{F}_T}(\Omega; \dbR^n)}>0, \ \text{ for sufficiently large $k$}.
$$
\end{corollary}

\noindent  {\bf  Proof. } By Theorem \ref{tt1},
there exists  a non-zero pair $(z_0,
z)^\top\in\mathbb{R}\times
L^2_{\mathcal{F}_T}(\Omega; \dbR^n)$,  such that
\begin{equation}\label{LLZ3}
z_0 \hat\xi^0+\mathbb{E} (z^\top \hat\xi(T;v))\geq
0,\q  \forall\; (\hat\xi^0, \hat\xi(T;v))^\top \in  \mbox{\rm
Var}_{L^2_\mathbb{F}(0, T; \dbR^m)}(f_0, f)(\bar{u}),
\end{equation}
where $\mbox{\rm
Var}_{L^2_\mathbb{F}(0, T; \dbR^m)}(f_0, f)(\bar{u})$ is defined by (\ref{var f_0f SOCP}).
By the classical well-posedness result for backward stochastic
differential equations (e.g., \cite[Section 4.1]{LZ}), (\ref{LLZ4})
admits a unique solution $ (\psi(\cdot), \Psi(\cdot))^\top\in L^2_\mathbb F(\Omega;
C([0, T]; \dbR^n))\times
  L^2_\mathbb F(0, T; \dbR^n)$.
Further,   by the duality relationship between
(\ref{LLZ2+}) and (\ref{LLZ4}),  (\ref{LLZ3}) indicates  $C_1^{\top} \psi+C_2^\top \Psi-z_0C_3=0 \mbox{\rm
\ a.e. }(t,\omega)\in[0,T]\times \Omega$. Hence,  (\ref{3.3 stochastic mp}) holds true.
\endpf

\subsubsection{The finite codimensionality of $\mbox{Var}_{L^2_\mathbb{F}(0, T; \dbR^m)}f(\bar{u})$}

To begin with, we  introduce  the  following backward stochastic differential
equation:
\begin{equation}\label{LLZ6}
\left\{
\begin{array}{ll}
d\phi(t)=-\big[A_1(t)^\top\phi(t)+A_2(t)^\top
\Phi(t)\big]dt+\Phi(t) dB(t), &t\in  [0,T],
\\ \ns\ds \phi(T)=\phi_T.
\end{array}
\right.
\end{equation}
By the duality relationship between  (\ref{LLZ2+}) and (\ref{LLZ6}),
 and the definition of  $f'(\bar{u})$,  it is easy  to show that
$$
f'(\bar{u})^*(\phi_T)=
C_1(\cdot)^\top\phi(\cdot)+C_2(\cdot)^\top\Phi(\cdot),\q \forall \;
\phi_T\in  L^2_{\mathcal{F}_T}(\Omega; \dbR^n),
$$
where $(\phi(\cdot), \Phi(\cdot))$ is the solution to  (\ref{LLZ6}) with the terminal
datum $\phi_T$. Then, by Theorem \ref{t2}, one has the following
results.
\begin{corollary}\label{cllzz}
The following  assertions are equivalent:

\ss

\noindent  {\bf (1)} The set  $\mbox{\rm Var}_{L^2_\mathbb{F}(0, T; \dbR^m)}f(\bar{u})$ is
finite codimensional in
$L^2_{\mathcal{F}_T}(\Omega; \dbR^n)$.

\ss

\noindent  {\bf (2)} There is a finite
codimensional closed subspace $\wt{\mathbb X}\subseteq
L^2_{\mathcal{F}_T}(\Omega; \dbR^n)$, such that
for  any $\phi_T\in \wt{\mathbb X}$, the solution $(\phi(\cdot), \Phi(\cdot))$ to $(\ref{LLZ6})$ satisfies

$$
\mathbb{E} |\phi_T|^2_{\dbR^n} \leq\displaystyle
\cC\mathbb{E}\int^T_0
|C_1(t)^\top\phi(t)+C_2(t)^\top\Phi(t)|^2_{\dbR^m}dt. $$

\ss

\noindent {\bf (3)}  There  is  a  Banach  space
$W$ and a compact operator $\mathcal G$  from
$L^2_{\mathcal{F}_T}(\Omega; \dbR^n)$ to  $W$,
such  that
for  any $\phi_T\in L^2_{\mathcal{F}_T}(\Omega;  \dbR^n)$, the solution $(\phi(\cdot), \Phi(\cdot))$ to $(\ref{LLZ6})$ satisfies
\begin{equation}\label{8.20-eq0}
\mathbb{E} |\phi_T|^2_{\dbR^n} \leq\displaystyle
\cC\Big(\mathbb{E}\int^T_0
|C_1(t)^\top\phi(t)+C_2(t)^\top\Phi(t)|^2_{\dbR^m}dt+|\mathcal
G \phi_T|_{W}^2\Big).
\end{equation}
\end{corollary}

Next,  we consider a special case, i.e.,
$C_2(\cdot)\equiv C_2\in\dbR^{n\times m}$, and
define a linear  operator $\cG:
L^2_{\mathcal{F}_T}(\Omega; \dbR^n)\to \dbR^n$
by
\begin{equation}\label{8.20-eq2}
\cG \phi_T = \phi(0),\q \forall\; \phi_T \in
L^2_{\mathcal{F}_T}(\Omega; \dbR^n),
\end{equation}
where  $\phi(\cdot)$ is the solution component  to  (\ref{LLZ6})
with the terminal datum $\phi_T$. By the well-posedness of (\ref{LLZ6})
(e.g., \cite[Section 4.1]{LZ}),  $\cG$ is a bounded
linear operator. Then we have the following result.
\begin{theorem}\label{LLZZ**}
The assertion {\bf (3)} in Corollary $\ref{cllzz}$  holds for the
operator $\cG$ defined by $(\ref{8.20-eq2})$,
 if and only if
$\rank(C_2)=n$.
\end{theorem}

\noindent{\bf Proof.}
First,  we assume  that $\rank(C_2)=n$.  By the well-posedness of
\eqref{LLZ6},
 there is a  positive constant $\cC$,  such that for
any $\phi_{T}\in L^2_{\mathcal{F}_T}(\Omega;
\dbR^n)$,
$$
|\phi(0)|_{\dbR^n}\leq
\cC|\phi_{T}|_{L^2_{\mathcal{F}_T}(\Omega;
\dbR^n)}.
$$
This  implies that $\cG(\widetilde\cM)$ is a bounded  set
in  $\dbR^n$  whenever $\widetilde\cM\subseteq
L^2_{\mathcal{F}_T}(\Omega; \dbR^n)$ is a
bounded set, and   $\cG$ is compact.
Set
\begin{equation}\label{8.20-eq6.1}
\hat\ell(t)\deq C_1(t)^\top\phi(t)+C_2^\top\Phi(t),\q \mbox{\rm a.e.
}(t,\omega)\in[0,T]\times \Omega.
\end{equation}
Since $\rank(C_2)=n$, there is a matrix $\wt C_2\in\dbR^{n\times m}$,
such that $\wt C_2 C_2^\top \Phi(t) = \Phi(t)$. Then we have
\begin{equation}\label{8.20-eq12}
\Phi(t)  =\wt C_2  \hat\ell(t)  -\wt C_2 C_1(t)^\top\phi(t),\q
\mbox{\rm a.e. }(t,\omega)\in[0,T]\times \Omega.
\end{equation}
If  a pair $(\phi(\cdot),\Phi(\cdot))$ solves \eqref{LLZ6}, by
\eqref{8.20-eq12}, we obtain the following forward equation:
\begin{eqnarray}\label{8.20-eq5}
\left\{\!\!
\begin{array}{ll}
d\phi(t)= -\big[A_1^\top(t)\phi(t)+ A_2^\top (t)\wt C_2 \hat\ell(t) -
A_2^\top (t)\wt C_2 C_1^\top\phi(t)\big]dt\\[+0.3em]
\q\q\q\;\;  +\big[\wt C_2 \hat\ell(t) - \wt C_2
C_1(t)^\top\phi(t)\big] dB(t), & t\in [0,T],
\\
\ns\ds \phi (0)=\phi(0).
\end{array}
\right.
\end{eqnarray}
By the  well-posedness results for stochastic differential equations
(e.g., \cite[Section 3.1]{LZ}), there is a positive constant $\cC$,
such that\vspace{-2mm}
\begin{equation*}\label{8.20-eq13}
\mathbb{E}|\phi(T)|^2_{\dbR^n} \leq \cC
\Big(\mathbb{E}\int^T_0 |
\hat\ell(t)|^2_{\dbR^m}dt+|\phi(0)|_{\dbR^n}^2\Big).
\end{equation*}
This, together with \eqref{8.20-eq6.1}, implies
\eqref{8.20-eq0}.

\ms

On the other hand, we assume that the
assertion {\bf (3)} in Corollary $\ref{cllzz}$  holds for the
operator $\cG$ defined by $(\ref{8.20-eq2})$.
  Consider the
following stochastic differential equation:
\begin{eqnarray}\label{8.20-eq6}
\left\{
\begin{array}{ll}
d\psi(t)= -[A_1(t)^\top\psi(t)+ A_2(t)^\top r(t)]dt+r(t) dB(t), &
t\in [0,T],
\\
\ns\ds \psi (0)=\psi_0,
\end{array}
\right.
\end{eqnarray}
where $\psi_0\in \dbR^n$ and $r(\cdot)\in L^2_\mathbb F(0,T;\dbR^n)$.
Clearly,  if $\psi(\cdot)$ is a solution to \eqref{8.20-eq6} with the
initial datum $\psi_0$,  then $(\phi(\cdot),\Phi(\cdot))\deq(\psi(\cdot),r(\cdot))$ is a solution to
\eqref{LLZ6} with the terminal datum $\phi(T)=\psi(T)$. By
\eqref{8.20-eq0} and the definition of $\mathcal G$,  we have the
following estimate  for any solution to (\ref{8.20-eq6}):
\begin{equation}\label{8.20-eq7}
\begin{array}{ll}\ds
\mathbb{E} |\psi(T)|^2_{\dbR^n} \leq
\cC\Big(\mathbb{E}\int^T_0
|C_1(t)^\top\psi(t)+C_2^\top
r(t)|^2_{\dbR^m}dt+|  \psi_0 |_{\dbR^n}^2\Big),\\
\ns\ds\hspace{5.1cm}\forall\;  \psi_0\in
\dbR^n\mbox{ and } r(\cdot)\in L^2_\mathbb
F(0,T;\dbR^n).
\end{array}
\end{equation}

If $\rank(C_2)<n$,  we can find a  $\hat
r\in\dbR^n\setminus\{0\}$, such that $C_2^\top \hat r=0$.  For any
$k\in \dbN$ with $1/k<T$,  let $r_k(\cdot)=  \chi_{[T-1/k,T]}(\cdot)\sqrt{k}\hat r$, and
$\psi_k(\cdot)$ be the corresponding solution to \eqref{8.20-eq6} for $r(\cdot)=r_k(\cdot)$ and
$\psi_0=0$,   where $\chi_{[T-1/k,T]}(\cdot)$ denotes the characteristic
function on $[T-1/k,T]$.
Then, $\displaystyle\mathbb{E}\int^T_{0}
|r_k(t)|^2_{\dbR^n}dt=|\hat r|^2_{\dbR^n}$, and
by the well-posedness of the equation (\ref{8.20-eq6}) (e.g.,
\cite[Section 3.1]{LZ}), we have that
$\psi_k(t)=0$ in $[0, T-1/k]$, and
it follows from \eqref{8.20-eq7} that
\begin{equation}\label{8.20-eq8}
\begin{array}{ll}\ds
\mathbb{E} |\psi_k(T)|^2_{\dbR^n} \leq \cC
\mathbb{E}\int^T_{T-1/k}
|C_1(t)^\top\psi_k(t)|^2_{\dbR^m}dt.
\end{array}
\end{equation}
By the well-posedness of \eqref{8.20-eq6} again, we
have  that
$$
\displaystyle\mathbb{E}\sup\limits_{0\leq t\leq T} |\psi_k(t)|^2_{\dbR^n} \leq \cC
\mathbb{E}\int^T_{0}
|r_k(t)|^2_{\dbR^n}dt=\cC|\hat r|^2_{\dbR^n},
$$
which, together with  \eqref{8.20-eq8}, implies that
\begin{equation}\label{8.20-eq10}
\begin{array}{ll}\ds
\mathbb{E} |\psi_k(T)|^2_{\dbR^n} \leq
\frac{\cC}{k}|\hat r|^2_{\dbR^n},
\end{array}
\end{equation}
with  $\cC$  independent of $k$.
Meanwhile,  $(\psi_k(\cdot),r_k(\cdot))$ can be regarded  as a
solution to \eqref{LLZ6}.  By the
well-posedness of \eqref{LLZ6}, there is a
positive constant $\cC$, such that
$$
\begin{array}{ll}\ds
\displaystyle\mathbb |\hat  r|^2_{\dbR^n}=\mathbb E\int_{T-1/k}^T|r_k(t)|^2_{\dbR^n}dt \leq
\cC\mathbb{E} |\psi_k(T)|^2_{\dbR^n}.
\end{array}
$$
This contradicts   \eqref{8.20-eq10} for a sufficiently large
$k\in\dbN$.  Hence,  $\rank(C_2)=n$ and this finishes the proof.
\endpf

\medskip

\begin{remark}
It is well known that when $C_2$ is an invertible  matrix,  the
following  controlled system is not exactly controllable $($e.g.,
$\cite[Proposition\; 6.3]{LZ})$:
\begin{eqnarray}\label{dd0}
\left\{
\begin{array}{ll}
d\xi= A_1(t)\xi dt+\big[A_2(t)\xi+C_2u(t)\big]dB(t),\q t\in [0, T],&\\
\ns \xi(0)=\xi_0,&
\end{array}
\right.
\end{eqnarray}
where $u$ is the control variable and $\xi$ is
the state variable. It implies that
\begin{eqnarray*}
\mathcal{R}_T\!\!\!&\deq& \!\!\!\Big\{ \xi(T; u)\in
L^2_{\mathcal{F}_T}(\Omega; \dbR^n)\  \Big|\ \xi(\cdot; u)\mbox{ is the solution
to }(\ref{dd0})\mbox{ for some } u(\cdot)\!\in\! L^2_\mathbb{F}(0, T;
\dbR^m) \Big\}\\
\!\!\!&\varsubsetneq&\!\!\! L^2_{\mathcal{F}_T}(\Omega; \dbR^n).
\end{eqnarray*}
Hence, the Robinson constraint
qualification does  not hold and the following estimate fails  for solutions
to $(\ref{LLZ6})$:
$$
\mathbb{E} |\phi_T|^2_{\dbR^n} \leq\displaystyle
\cC\mathbb{E}\int^T_0
|C_2^\top\Phi(t)|^2_{\dbR^m}dt,\q \forall\;
\phi_T\in L^2_{\mathcal{F}_T}(\Omega; \dbR^n).
$$
Theorem $\ref{LLZZ**}$ implies  that when  $C_2$ is an invertible  matrix,  the
above inequality will  hold, only  if an extra term with
a compact operator $\mathcal G$ is added into it.  This
weaker estimate is enough to obtain a Fritz John
condition for  the optimal control
 problem $(\ref{ocp3})$.
\end{remark}

\begin{example}
Let us consider the following  controlled  stochastic differential equation:
\begin{equation}\label{sto ode sy1}
\left\{
\begin{array}{lll}\ds
dy(t) = \big[A_1(t)y(t) + C_1(t)u(t)\big]dt +
 \big[A_2(t)y(t) + C_2 u(t)\big]dB (t), &t\in  [0, T],\\
\ns\ds y(0)=y_0.
\end{array}
\right.
\end{equation}
The
system \eqref{sto ode sy1} is indeed a
special case of \eqref{1LLZ}. It is a
model of investment in a financial
market $($e.g., $\cite{Merton})$.  The constant matrix  $C_2$ means  that
volatilities of the stocks are
constant. This is reasonable  if $T$ is
not very large $($e.g., $\cite{Merton1})$.
By Theorem $\ref{LLZZ**}$, in order to guarantee
the existence  of  a non-zero pair
$(z_0, z)^\top\in\mathbb{R}\times
L^2_{\mathcal{F}_T}(\Omega; \dbR^n)$ as
the Lagrange multiplier,  a sufficient
condition is $\rank(C_2)=n$. This
condition means that there are enough
numbers of stocks in the market, which
holds for a modern financial market.
\end{example}

\section{Further comments} \label{sec-5}

This paper focuses on establishing an enhanced Fritz John (first-order necessary) condition for optimization problems in infinite-dimensional spaces using a finite codimensionality method.
This approach  offers a unified framework for dealing with
 various types of optimal control problems with different state constraints.

In the field of infinite-dimensional optimization theory, many constraint qualifications have been introduced to ensure the existence of nontrivial or normal Lagrange multipliers (e.g., see \cite{BS} and \cite{Borgens2020}). However, validating these constraint qualifications is generally challenging.

To illustrate the difficulty,  let us consider a special
case where
\begin{equation}\label{5.1}
\mathcal{K}= V\quad \text{ and } \quad E=\{0\},
\end{equation}
with $V$ being a Banach space, and both $f_0$	
and $f$ being continuously differentiable  in the problem {\bf (P)}. In this case, $\mathcal{Z}(\bar u) =\Im (f^\prime(\bar u))$ (see (\ref{c(bar u)}) for the definition of $\mathcal{Z}(\bar u)$).  In \cite{BS}, it was required that $\Im(f'(\bar{u}))$ has a nonempty relative interior, which is equivalent to the condition that $\Im (f^\prime(\bar u))$ is closed in $X$.

The following proposition demonstrates  that, similar to the finite codimensionality condition discussed in Section 3, the closed range condition can also be established through certain {\em a priori} estimates. However, validating these estimates is not easier than proving the associated estimates for finite codimensionality.

Let $V$ be a reflexive Banach space,  $X$ be a Hilbert space
and  ${\bf F}\in \mathcal{L}(V; X)$.  We have the following result.

\begin{proposition}\label{prop21}
The following assertions are equivalent:

\ss

\noindent {\bf (1)} $\Im({\bf F})$ is closed  in $X$.

\ss

\noindent {\bf (2)} There exists a positive constant $\mathcal{C}$,  such that
\begin{equation}\label{2p2}
|h|_{(\overline{\Im({\bf F})})'}\leq
\mathcal{C} |{\bf F}^*(h)|_{V'},\quad\forall\;
h\in \big(\overline{\Im({\bf F})}\big)'.
\end{equation}

\noindent {\bf (3)}  There exists a positive constant $\mathcal{C}$,  such that
\begin{equation}\label{2p3}
|h|_{X'}\leq \mathcal{C}
\Big(|{\bf F}^*(h)|_{V'}+\Big|\Pi_{\big({\overline{\Im({\bf F})}}^{\prime}\big)^\bot}h\Big|_{X'}\Big),
\quad\forall\;
h\in X'.
\end{equation}
\end{proposition}

\noindent {\bf Proof.}  First,  we prove the equivalence   between
{\bf (1)} and {\bf (2)}. Set $\tilde X_0=\overline{\Im({\bf F})}$. Then,  $\tilde  X_0$ is
a  Banach space  under the norm $|\cdot|_X$. Consider the mapping ${\bf F}:
V\rightarrow \tilde X_0$ and the identity mapping $\mathbb I_{\tilde X_0}:
\tilde X_0\rightarrow  \tilde X_0$.
It is easy to check that ${\bf F}\in \mathcal{L}(V;
\tilde  X_0)$ and $\mathbb I_{\tilde X_0}\in \mathcal{L}(\tilde X_0)$. If  $\Im({\bf F})$ is
closed  in $X$,  it  follows that $\Im(\mathbb I_{\tilde X_0})\subseteq
\Im({\bf F})$. By the range  comparison theorem,   there exists a positive
constant $\mathcal  C$,   such that
$$|h|_{\tilde X_0'}\leq \mathcal{C} |{\bf F}^*(h)|_{V'},\quad\forall\; h\in \tilde X_0',
$$
which means that {\bf (2)} holds.

\ss

Conversely,  if {\bf (2)} is true, then $\Im(\mathbb I_{\tilde X_0})\subseteq
 \Im({\bf F})$, which implies
{\bf (1)}.

\medskip

Next, we prove that {\bf (2)} implies  {\bf (3)}.  Since
$\tilde X_0=\overline{\Im({\bf F})}$ is a closed  subspace  of the Hilbert space $X$,
we have the direct sum decomposition $X'=\tilde X_0' \oplus
\big(\tilde X_0'\big)^\bot$. Hence, for any $h\in  X'$,
$$
h= \Pi_{\tilde X_0'}h+
\Pi_{\big(\tilde X_0'\big)^\bot}h.
$$
By {\bf (2)}, it holds that
\begin{eqnarray*}
|h|_{X'}&\leq& \big| \Pi_{\tilde X_0'}h\big|_{X'}+\big|
\Pi_{\big(\tilde X_0'\big)^\bot}h\big|_{X'} \leq \mathcal C
\big|{\bf F}^*\big(\Pi_{\tilde X_0'}h \big)\big|_{V'}
+\big|\Pi_{\big(\tilde X_0'\big)^\bot}h \big|_{X'}\\
&\leq&  \mathcal C\big( |{\bf F}^*(h)|_{V'}+ \big|{\bf F}^*\Big(
\Pi_{\big(\tilde X_0'\big)^\bot}h\big)\big|_{V'}\big)+\big|
\Pi_{\big(\tilde X_0'\big)^\bot}h\big|_{X'}\\
&\leq& \mathcal C \big(|{\bf F}^*(h)|_{V'}+ \big|
\Pi_{\big(\tilde X_0'\big)^\bot}h\big|_{X'}\big).
\end{eqnarray*}
Conversely,  it is obvious that {\bf (3)} implies {\bf (2)}.   \endpf

\begin{remark}
From the proof  of Proposition $\ref{prop21}$,  the equivalence between  {\bf (1)}  and  {\bf
(2)}  continues to hold for any Banach space $X$.
\end{remark}
\begin{remark}\label{remark2.6*}
By $\cite[\mbox{Theorems }1.1.1 \mbox{ and }   1.1.2]{H}$
and $\cite[\mbox{Lemma }   4]{P}$,   if $V$ and $X$
are two Hilbert  spaces, and ${\bf F}:
\mathcal{D}({\bf F})\subseteq V\rightarrow X$ is a
linear closed densely defined operator,  the
following assertions are equivalent:

\medskip

\noindent $(a)$ $\Im({\bf F})$ is closed in $X$.

\medskip

\noindent $(b)$ There exists a positive constant $\mathcal{C}$,  such that
$$
|h|_{V}\leq \mathcal{C} |{\bf F}(h)|_X,\quad\forall\; h\in  \mathcal{D}({\bf F})\cap
 \overline{\Im({\bf F}^*)},
$$
where ${\bf F}^*:   \mathcal{D}({\bf F}^*)\rightarrow V$ denotes  the  adjoint  operator of ${\bf F}$.

\medskip

\noindent $(c)$ $\Im({\bf F}^*)$ is closed in $V$.

\medskip

\noindent $(d)$ There exists a positive constant $\mathcal{C}$,  such that
$$
|h|_{X}\leq \mathcal{C} |{\bf F}^*(h)|_V,\quad\forall\; h\in
\mathcal{D}({\bf F}^*)\cap \overline{\Im({\bf F})}.
$$

\medskip

\noindent $(e)$ There exists a positive constant $\mathcal{C}$,  such that
$$
|h|_{X}\leq \mathcal{C}
\big(|{\bf F}^*(h)|_V+|\Pi_{\ker({\bf F}^*)}h|_X\big),\quad\forall\;
h\in
\mathcal{D}({\bf F}^*),
$$
where $\ker({\bf F}^*)$ denotes the null space of ${\bf F}^*$.

It is straightforward to
 show that when $V$ and $X$ are Hilbert  spaces, and ${\bf F}\in\mathcal{L}(V; X)$,
every assertion in Proposition $\ref{prop21}$ is equivalent to one of the
conditions $(a)$, $(b)$, $(c)$, $(d)$  and $(e)$ above. \end{remark}

In  the infinite-dimensional optimization
problem  {\bf (P)}, if  $V$ is a reflexive Banach space,  $X$ is a Hilbert space,
$\mathcal{K}= V$, $E=\{0\}$,
and $f$ and $f^0$ are continuously
differentiable,  then by Proposition \ref{prop21},    the assertion that
\begin{equation}\label{2p5}
\Im(f'(\bar{u})) \mbox{  is closed in  }  X
\end{equation}
is equivalent to  the  estimate
\begin{equation}\label{2p6+}
|h|_{\big(\overline{\Im(f'(\bar{u}))}\big)'}\leq
\mathcal{C} |f'(\bar{u})^*(h)|_{V'},\quad\forall\;
h\in \Big(\overline{\Im(f'(\bar{u}))}\Big)',
\end{equation}
or the estimate
\begin{equation}\label{2p6}
|h|_{X'}\leq \mathcal{C} \Big(|f'(\bar{u})^* h|_{V'}
+\Big|\Pi_{\big({\overline{\Im(f'(\bar{u}))}}^{\prime}\big)^{\bot}} h\Big|_{X'}\Big),\quad\forall\;
h\in X'.
\end{equation}

Although  the condition related to  nonempty relative interior for $\mathcal{Z}(\bar{u})$
 is generally weaker than the finite codimensionality condition,
  it is difficult to verify   directly. Even in the special case \eqref{5.1}, where the nonemptiness of the relative interior of $\mathcal{Z}(\bar{u})$ is equivalent to the closedness of $\Im(f'(\bar{u}))$ in $X$, it is often challenging to precisely characterize $ \big(\overline{\Im(f'(\bar{u}))}\big)'$
  or $\big({\overline{\Im(f'(\bar{u}))}}^{\prime}\big)^\bot$. This
  makes the conclusions in Proposition \ref{prop21} (i.e. (\ref{2p6+}) and (\ref{2p6})) difficult to apply in practice. As a result, to ensure the closedness of $\Im(f'(\bar{u}))$, one has to assume   stronger sufficient conditions, such as the surjectivity of $f'(\bar{u})$ or the finite codimensionality of $\Im(f'(\bar{u}))$. An illustrative example is provided below  for further explanation.

\ss

Let $G\subseteq\dbR^n$ be a convex
bounded domain with a smooth
boundary $\Gamma$  and $T>2\sup\limits_{x,y\in
G}|x-y|_{\dbR^n}$. Put $Q=G\times(0,T)$
and $\Sigma=\Gamma\times(0,T)$. Assume
that $G_0$  is a nonempty open subset
of $G$, and denote by $\chi_{G_0}$ the
characteristic function  of $G_0$.
Consider the following controlled wave equation:
\begin{eqnarray}\label{31}\left\{
\begin{array}{ll}
y_{tt}-\Delta y+ay =\chi_{G_0}  u  &\mbox{ in }
Q,\\
\ns\ds y=0 &\mbox{ on }  \Sigma,\\ \ns\ds
y(0, x)=y_0(x), \ y_t(0, x)=y_1(x)  &\mbox{ in }  G,
\end{array}
\right.
\end{eqnarray}
with the cost functional
\begin{equation}\label{31cost}
J(u(\cdot))=\ds\frac{1}{2}\ds\int_Q
\big[q(x, t)|y(x,t)|^2+r(x,t)|u(x,
t)|^2\big]dxdt,
\end{equation}
where  $u\in L^2(Q)$ is  the control variable,
 $(y, y_t)$ is the state variable,
 $(y_0,
y_1)^\top\in H^1_0(G)\times L^2(G)$ is the
initial value, and $a, q, r\in
L^\infty(Q)$.

\ss

 For a given $(y^0_d, y^1_d)^\top\in H^1_0(G)\times
L^2(G)$, let
$$
\mathcal{U}_{ad}=\big\{  u(\cdot)\in  L^2(Q)\
\big|\ \mbox{the   solution }y(\cdot;u)\mbox{ to
}\eqref{31} \mbox{  satisfies that } (y(T;u),
y_t(T;u))=(y^0_d, y^1_d)  \big\}.
$$
Assume that $\bar{u}(\cdot)\in \mathcal{U}_{ad}$  is an optimal control of the following
optimal control problem:
$$
\text{Minimize } J(u(\cdot)) \quad
\text{subject to }u(\cdot)\in
\mathcal{U}_{ad}.
$$

This optimal control problem is a
special case of the optimization
problem {\bf (P)} obtained by setting
$$\begin{array}{ll}
&V=\mathcal{K}=L^2(Q),  \ \ X=H^1_0(G)\times
L^2(G),  \ \ E=\{(y^0_d, y^1_d)\},\\[3mm]
&f_0(u)=J(u(\cdot))\ \mbox{ and }\  f(u)=(y(\cdot, T; u),  y_{t}(\cdot, T; u))^\top,
\end{array}
$$
where  $y(\cdot;u)$ is the solution to
(\ref{31}) associated to $u(\cdot)$.
  Further, it
is  easy to show that
$$
\begin{array}{rl}
&\Im(f'(\bar{u}))=\big\{  (y(\cdot, T;u), y_t(\cdot, T;u))^\top\in
H^1_0(G)\!\times\! L^2(G)\ \big|\ y(\cdot;u)\mbox{ is the
solution to } (\ref{31})\mbox{ with }
\\[3mm]
&\quad\qq\qq\qq\qq\qq\qq\qq\qq\qq\qq\qq (y_0,  y_1)=(0,  0)\mbox{ for
some } u\in L^2(Q)  \big\}.
\end{array}
$$

By Proposition \ref{prop21},  the closedness of $\Im(f'(\overline{u}))$
is equivalent to the following estimate:
\begin{eqnarray}\label{2**}
\begin{array}{ll}
&\displaystyle|(\phi_1, \phi_2)|_{L^2(G)\times H^{-1}(G)}\leq
\mathcal{C} \Big(|\phi|_{L^2(G_0\times(0, T))}
+\Big|\Pi_{\big( {\overline{\Im(f'(\bar{u}))} }^{\prime}\big)^\bot}(\phi_1,
\phi_2)\Big|_{L^2(G)\times
H^{-1}(G)}\Big),\\[4mm]
&\qq\qq\qq\qq\qq\qq\qq\qq\qq\qq\qq\forall\; (\phi_1, \phi_2)^\top\in L^2(G)\times
H^{-1}(G);
\end{array}
\end{eqnarray}
and by Theorem \ref{t2}, the finite
codimensionality of $\Im(f'(\bar{u}))$ is
equivalent to the following estimate:
\begin{equation}\label{8-18-eq2}
\displaystyle|(\phi_1, \phi_2)|_{L^2(G)\times
H^{-1}(G)}\leq \mathcal{C}
 |\phi|_{L^2(G_0\times(0, T))},
\q \forall\; (\phi_1, \phi_2)^\top\in \mathcal{H}_0,
\end{equation}
where  $\mathcal{H}_0$ is a
 finite codimensional closed subspace of
$L^2(G)\times H^{-1}(G)$ and $\phi$  is the solution to
\begin{eqnarray}\label{32}\left\{
\begin{array}{lll}
\ds\phi_{tt}-\Delta \phi +a\phi=0  &\mbox{ in
}Q,\\
\ns\ds\phi=0 &\mbox{ on }\Sigma,\\
\ns\ds\phi(T, x)=\phi_1(x),\  \phi_t(T, x)=\phi_2(x) &\mbox{
in }G.
\end{array}
\right.
\end{eqnarray}

\ms

We have the following results.

\ss

(1) In the  case where $a=0$, if   $(G, T,
G_0)$ fulfills the geometric control
condition (e.g., \cite{BLR}), then
$f'(\bar{u})$ is surjective, which
implies the closedness of
Im$(f'(\bar{u}))$.

\smallskip

Indeed, under the geometric
control condition,  the system (\ref{31})  (with $a=0$)
is exactly controllable.  Hence, any solution $\phi$ to (\ref{32})  (with $a=0$)
satisfies the  estimate:
\begin{eqnarray}\label{8-18-eq1}
\begin{array}{ll}
&\displaystyle|(\phi_1, \phi_2)|_{L^2(G)\times
H^{-1}(G)}\leq \mathcal{C}
 |\phi|_{L^2(G_0\times(0, T))}, \q \forall\;
(\phi_1, \phi_2)^\top\in L^2(G)\times H^{-1}(G),
\end{array}
\end{eqnarray}
which implies the estimate  \eqref{2**}.

Meanwhile,  if
$(G,T,G_0)$ does not fulfill the
geometric control condition,
\eqref{2**} fails  to hold. Indeed,
since  the  system \eqref{31} (with
$a=0$) is approximately controllable at
time $T$ (e.g.,\cite{Laurent}),
$\Im(f'(\bar{u}))$ is dense in $X$.  Because the geometric control
condition does not hold,  the
 system \eqref{31} (with $a=0$) is not
exactly controllable,  and therefore, $\Im(f'(\bar{u}))$ is not
closed. This shows that when $a=0$ and the time $T$ is
sufficiently large, the closedness of $\Im (f'(\bar{u}))$ is
equivalent to the geometric control condition. In this particular instance,
 the conditions  guaranteeing  the closedness and
finite codimensionality of $\Im(f'(\bar{u}))$ are identical.

\medskip

(2) In the general case where $a\in  L^\infty(Q)$, if $(G, T, G_0)$
fulfills the geometric control condition,  $\Im(f'(\bar{u}))$
is  finite codimensional in $X$, which  implies the closedness of
Im$(f'(\bar{u}))$.

\smallskip

As proved in \cite[Proposition 6.2]{LLZ},  the finite
codimensionality condition for $\Im(f'(\bar{u}))$ holds
whenever  $(G,T,G_0)$ fulfills the geometric control  condition.
However, under the same  condition, proving the surjectivity of
$f'(\bar{u})$ or the estimate (\ref{8-18-eq1}) remains open. In
addition, verifying   the estimate $(\ref{2**})$ or
the closedness of $\Im(f'(\bar{u}))$ directly currently  appears  difficult.

\appendix

\section{Proofs of (\ref{appendix}) and (\ref{appendix_add})}\label{A}

In this appendix, we prove (\ref{appendix}) and (\ref{appendix_add})  in
Example \ref{remark2.2}.

\medskip

\noindent{\bf Proof of
(\ref{appendix}).}
We first prove that  for any $e\in \mathcal{K}$,
\begin{equation}\label{1.21-eq2}
\big\{ f'(e; v)\in X\ \big|\  v\in
\mathcal{R}_{\mathcal{K}}(e)\cap B_V(0, 1)
\big\}\subseteq \mbox{\rm Var}_{\mathcal{K}}f(e).
\end{equation}
Indeed,  for any $v\in \mathcal R_{\mathcal{K}}(e)\cap B_V(0,
1)$,  there exists a $\alpha>0$ and a $y\in \mathcal{K}$, such that
$v=\alpha(y-e)$. Let $\ell = \frac{1}{2\a}$.   Then
$$ e+\ell
v=e+\ell\alpha(y-e)=\ell\alpha y+(1-\ell\alpha)e\in \mathcal{K}.
$$
Take a sequence $\{h_k\}_{k=1}^\infty\subseteq (0, \ell)$ with
$\lim\limits_{k\rightarrow \infty}h_k=0$ and set $e_k=e+h_k v$ for
$k\in\dbN$. Then $\{e_k\}_{k=1}^\infty\subseteq \mathcal{K}$ and
$|e_k-e|_V=h_k|v|_V\leq h_k$. Consequently,
$$
f'(e;v)=\lim\limits_{k\rightarrow\infty}\displaystyle\frac{f(e+h_k
v)-f(e)}{h_k}
=\lim\limits_{k\rightarrow\infty}\displaystyle\frac{f(e_k)-f(e)}{h_k}\in
\mbox{\rm Var}_{\mathcal{K}}f(e).
$$
Thus,  \eqref{1.21-eq2} holds.

For any $v\in  \cT_{\mathcal{K}}(e)\cap B_V(0, 1) $, there exists a
sequence $\{v_n\}_{n=1}^\infty\subseteq \mathcal
R_{\mathcal{K}}(e)$  such that
$\lim\limits_{n\rightarrow\infty}v_n=v$. Without loss of generality,
we may assume that $|v_n|_V\leq 1$ for every $n\in\dbN$. In fact, if
$|v|_V<1$, then $|v_n|_V<1$ for sufficiently large $n\in\dbN$. On
the other hand,  if $|v|_V=1$, we set $\tilde
v_n=\displaystyle\frac{v_n}{|v_n|_V}$. It is easy to show that
$\tilde v_n\in \mathcal R_{\mathcal{K}}(e)\cap B_V(0, 1)$ and
$\lim\limits_{n\rightarrow\infty}\tilde v_n=v$.

By the locally Lipschitz continuity of $f$, we get that
\begin{eqnarray*}
&&|f'(e; v_n)-f'(e; v)|=\Big|\lim\limits_{h\rightarrow 0^+} \frac{f(e+h v_n)-f(e)}{h}-
\lim\limits_{h\rightarrow 0^+} \frac{f(e+h v)-f(e)}{h}
\Big|\\
&&\leq C|v_n-v|_V,
\end{eqnarray*}
which, together with $\lim\limits_{n\rightarrow\infty}v_n=v$, implies
$$
\lim\limits_{n\rightarrow\infty}f'(e; v_n)=f'(e;  v).
$$
Since $f'(e;  v_n)\in\mbox{\rm Var}_{\mathcal{K}}f(e)$, there exists  a sequence
$\{h_k^n\}_{k=1}^\infty\subseteq (0,+\infty)$ with
$\lim\limits_{k\rightarrow\infty} h_k^n=0$ and a sequence
$\{e_k^n\}_{k=1}^\infty\subseteq\mathcal{K}$, such that
\vspace{-2mm}
$$
|e_k^n-e|_V\leq h_k^n\q\mbox{and}\q f'(e;v_n)
=\lim\limits_{k\rightarrow \infty}\displaystyle\frac{f(e_k^n)-f(e)}{h_k^n}.
$$
Therefore,  we may find a sequence $\{h_k\}_{k=1}^\infty\subseteq
(0,+\infty)$ with $\lim\limits_{k\rightarrow\infty} h_k=0$  and a
sequence $\{e_k\}_{k=1}^\infty\subseteq\mathcal{K}$, such that
\vspace{-2mm}
$$
|e_k-e|_V\leq h_k\q\mbox{and}\q f'(e; v)
=\lim\limits_{k\rightarrow \infty}\displaystyle\frac{f(e_k)-f(e)}{h_k}.
$$
This implies that\vspace{-2mm}
$$
 f'(e; v)
\in\mbox{\rm Var}_{\mathcal{K}}f(e).
$$
This, together with the arbitrariness of $v\in
\cT_{\mathcal{K}}(e)\cap B_V(0, 1)$, implies (\ref{appendix}).
\endpf

\medskip

\medskip

\noindent{\bf Proof of
(\ref{appendix_add}).}
By \eqref{appendix},  it suffices to prove that   for any $e\in \mathcal{K}$,
\begin{equation}\label{1.21-eq1}
\mbox{\rm Var}_{\mathcal{K}}f(e)\subseteq\big\{ f'(e;v)\in X\
\big|\  v\in  \cT_{\mathcal{K}}(e)\cap B_V(0, 1)
\big\}.
\end{equation}

We first handle the case where $f$ is Fr\'echet differentiable. For
any $e\in\mathcal{K}$ and $\xi\in\mbox{\rm Var}_{\mathcal{K}} f(e)$,  there
exists a sequence $\{h_k\}_{k=1}^\infty\subseteq(0,+\infty)$ with
$\lim\limits_{k\rightarrow\infty} h_k=0$  and a sequence
$\{e_k\}_{k=1}^\infty\subseteq\mathcal{K}$, such
that\vspace{-2mm}
\begin{equation}\label{newA3}
|e_k-e|_V\leq h_k\q\mbox{and}\q \xi=\lim\limits_{k\rightarrow \infty}\displaystyle\frac{f(e_k)-f(e)}{h_k}.
\end{equation}
Since $\displaystyle\Big|\frac{e_k-e}{h_k}\Big|_V\leq 1$,  there are
subsequences  of $\{h_k\}_{k=1}^\infty$ and $\{e_k\}_{k=1}^\infty$
(denoted still by themselves for simplicity), and a $v\in B_V(0, 1)$,
such that\vspace{-2mm}
\begin{equation}\label{1.21-eq3}
\displaystyle\frac{e_k-e}{h_k} \mbox{  converges  weakly to }v \mbox{ in }V, \q \mbox{ as }k\rightarrow  \infty.
\end{equation}
Since $ \cT_{\mathcal{K}}(e) =\overline{\bigcup\limits_{\alpha\geq0}
\alpha\big(\mathcal{K} -e\big)}$,  we have that
$\displaystyle\frac{e_k-e}{h_k}\in  \cT_{\mathcal{K}}(e) $. Noting
that $ \cT_{\mathcal{K}}(e) $ is closed and convex, we obtain that  $v\in
 \cT_{\mathcal{K}}(e) $.  Moreover, by the definition of  Fr\'{e}chet differentiability,
$$
\frac{f(e_k)-f(e)}{h_k}
=f'(e)\left(\frac{e_k-e}{h_k}\right)+\frac{o(|e_k-e|_V)}{h_k}.
$$
Note that $f'(e)\Big(\displaystyle\frac{e_k-e}{h_k}\Big)$ converges weakly to $f'(e)v$ and the remainder term $\frac{o(|e_k-e|_V)}{h_k}$ converges strongly to $0$,
as $k\rightarrow \infty$. This,  together  with \eqref{newA3}, implies that
$$\xi=\lim\limits_{k\rightarrow\infty}\displaystyle\frac{f(e_k)-f(e)}{h_k}
=f'(e)v.$$
Hence,  \eqref{1.21-eq1} is proved.

Next, we deal with the case where $V$ is finite dimensional. Note
that the convergence in \eqref{1.21-eq3} becomes strong convergence now, and
by the locally  Lipschitz continuity of $f$,
$$
\xi=\lim\limits_{k\rightarrow\infty}\displaystyle\frac{f(e_k)-f(e)}{h_k}
=  \lim\limits_{k\rightarrow\infty}\displaystyle\frac{f(e+h_k v)-f(e)}{h_k}  +\lim\limits_{k\rightarrow\infty}\displaystyle\frac{f(e+h_k \frac{e_k-e}{h_k})-f(e+h_k v)}{h_k} =f'(e;  v).
$$
This proves \eqref{1.21-eq1}.
\ss
\endpf

\section{Proof of Proposition  \ref{add**}}\label{B}

In this appendix, we give a proof of Proposition \ref{add**}.
The whole proof is divided into two parts.

\medskip

{\bf Step 1}. In this step,  we prove  that
\begin{equation}\label{1.27-eq1}
\overline{\mbox{span}}\big\{ d-x_0\in X\  \big|\  d\in D\big\}
=\overline{\mbox{span}}\big\{ d_1-d_2\in X\  \big|\  d_1, d_2\in D\big\},\q \forall\; x_0\in \coh D.
\end{equation}

Indeed, for any
$\xi\in \mbox{span}\big\{ d_1-d_2\in X\  \big|\  d_1, d_2\in D\big\}$,
there exists a  $k\in\dbN$,
$\alpha_1, \cdots, \alpha_k\in \dbR$ and
$d_{11}, \cdots, d_{1k}$, $d_{21}, \cdots, d_{2k}\in D$,  such that
$$
\xi=\sum\limits_{i=1}^k  \alpha_i(d_{1i}-d_{2i})=\sum\limits_{i=1}^k  \alpha_i(d_{1i}-x_0)
-\sum\limits_{i=1}^k  \alpha_i(d_{2i}-x_0)
 \in \overline{\mbox{span}}\big\{ d-x_0\in X\  \big|\  d\in D\big\}.
$$
Consequently,
\begin{equation}\label{1.27-eq2}
\overline{\mbox{span}}\big\{ d_1-d_2\in X\  \big|\  d_1, d_2\in D\big\}
\subseteq \overline{\mbox{span}}\big\{ d-x_0\in X\  \big|\  d\in D\big\}.
\end{equation}

On the other hand, for any $\xi\in  \mbox{span}\big\{ d-x_0\in X\  \big|\  d\in D\big\}$,
there exists a  $k\in\dbN$, $\alpha_1, \cdots, \alpha_k\in\dbR$ and
$d_1, \cdots, d_k\in D$,  such that
$
\xi=\sum\limits_{i=1}^k \alpha_i (d_i-x_0).
$
Notice that $x_0\in \coh D$.  We can find a sequence $\{x_j\}_{j=1}^{\infty}\subseteq \mbox{co} D$ such that $x_j\to x_0$ as $j\to\infty$. From the definition of the convex hull, for any $j\in\dbN$,  there exists a
 $k_j\in\dbN$,
$\beta_{j 1}, \cdots, \beta_{j k_j}\in[0, 1]$ and $d_{j 1}, \cdots, d_{j k_j}\in D$,  such that
$
\sum\limits_{\ell=1}^{k_j} \beta_{j \ell}=1    \mbox{ and }
x_j=\sum\limits_{\ell=1}^{k_j}  \beta_{j \ell} d_{j \ell}.
$
Then, it follows that \vspace{-3mm}
\begin{eqnarray*}
\xi&=&\lim\limits_{j\rightarrow \infty}\sum\limits_{i=1}^k \alpha_i (d_i-x_j)
=\lim\limits_{j\rightarrow \infty}
 \sum\limits_{i=1}^k \alpha_i \(d_i-\sum\limits_{\ell=1}^{k_j}  \beta_{j \ell} d_{j \ell}\)\\
&=&\lim\limits_{j\rightarrow \infty}\sum\limits_{i=1}^k\sum\limits_{\ell=1}^{k_j}
 \alpha_i \beta_{j \ell}(d_i-d_{j \ell})\in\overline{\mbox{span}}\big\{ d_1-d_2\in X\  \big|\  d_1, d_2\in D\big\}.
 \end{eqnarray*}
This yields
\begin{equation}\label{1.27-eq3}
\overline{\mbox{span}}\big\{ d-x_0\  \big|\  d\in D\big\}\subseteq
\overline{\mbox{span}}\big\{ d_1-d_2\  \big|\  d_1, d_2\in D\big\}.
\end{equation}
Combining \eqref{1.27-eq2} and \eqref{1.27-eq3}, we obtain \eqref{1.27-eq1}.

\ms

{\bf Step 2}. In this step, we prove that if  for some $x_0^1\in \coh D$ and $\d>0$, there exists
 a $w_1\in\coh\big(D - x_0^1\big)$  such that
\begin{equation}\label{1.27-eq7}
B_{X}(w_1, \delta)\cap
\overline{\mbox{span}}\big\{ D-x_0^1\big\}\subseteq
\coh\big( D-x_0^1\big),
\end{equation}
then for any $x_0^2\in \coh D$,
\begin{equation}\label{1.27-eq4}
B_{X}(w_2, \delta)\cap
\overline{\mbox{span}}\big\{ D-x_0^2\big\}\subseteq
\coh\big( D-x_0^2\big),
\end{equation}
with $w_2=w_1+x_0^1-x_0^2
\in\coh\big( D-x_0^2\big)$.

\smallskip

Indeed, for any $\eta\in B_{X}(w_2, \delta)\cap
\overline{\mbox{span}}\big\{ D-x_0^2\big\}$, it holds that
$|\eta-w_2|_X\leq \delta$. Hence,  for $\xi=\eta-x_0^1+x_0^2$,
$$|\xi-w_1|_X=|\eta-x_0^1+x_0^2-w_1|_X=|\eta-w_2|_X\leq \delta.$$
This implies
\begin{equation}\label{1.27-eq5}
\xi\in B_{X}(w_1, \delta).
\end{equation}

Meanwhile, since $\eta\in \overline{\mbox{span}}\big\{ D-x_0^2\big\}$, there is a sequence $\{\eta_j\}_{j=1}^{\infty}\subseteq \mbox{span}\big\{ D-x_0^2\big\}$ so that $\eta_j\to \eta$ as $j\to \infty$. Then, for any $j\in\dbN$,  there exists a  $k_j\in\dbN$,
$\beta_{j 1}, \cdots, \beta_{j  k_j}\in\dbR$ and $d_{j 1}, \cdots, d_{j k_j}\in D$, such that
$$
\eta_j\deq\sum\limits_{\ell=1}^{k_j} \beta_{j \ell}(d_{j \ell}-x_0^2)\in
\mbox{span}\{D-x_0^2\}\  \mbox{ and }\  \eta=\lim\limits_{j\rightarrow\infty}\eta_j.
$$
Consequently,
\vspace{-2mm}
\begin{eqnarray*}
&&\xi=\eta-x_0^1+x_0^2=
\lim\limits_{j\rightarrow\infty}\Big[
\sum\limits_{\ell=1}^{k_j} \beta_{j \ell}(d_{j \ell}-x_0^2)
-x_0^1+x_0^2
\Big]\\
&&\quad=\lim\limits_{j\rightarrow\infty}\Big[
\sum\limits_{\ell=1}^{k_j} \beta_{j \ell}(d_{j \ell}-x_0^1)-
\sum\limits_{\ell=1}^{k_j} \beta_{j \ell}(x_0^2-x_0^1)
+(x_0^2-x_0^1)
\Big].
\end{eqnarray*}
Since $x_0^2\in\coh D$, we get that $x_0^2-x_0^1\in  \overline{\mbox{span}}\big\{ D-x_0^1\big\}.$
Then, the above equality implies that  $\xi\in \overline{\mbox{span}}\big\{ D-x_0^1\big\}$. This, together with \eqref{1.27-eq5}, implies that
\begin{equation}\label{1.27-eq8}
\xi\in B_{X}(w_1, \delta)\cap
	\overline{\mbox{span}}\big\{ D-x_0^1\big\}.
\end{equation}
Combining \eqref{1.27-eq7} and \eqref{1.27-eq8}, we find that
$
\xi\in \coh\big( D-x_0^1\big)=\coh D-x_0^1,
$
which implies that
$\eta\in \coh D-x_0^2=\coh\big( D-x_0^2\big).$
This, together with the arbitrariness of $\eta\in B_{X}(w_2, \delta)\cap
\overline{\mbox{span}}\big\{ D-x_0^2\big\}$, yields \eqref{1.27-eq4}.
\endpf

\section{Proof of Lemma \ref{codim converg}}\label{C}

In this appendix, we give a proof of Lemma \ref{codim converg}.
The whole proof is divided into four parts.

\medskip

{\bf Step 1. }  First, by Definition  \ref{llzd6} and the finite codimensionality of
$D$,   there exists a $z_0\in  \overline{\mbox{co}}D$  such that
$\overline{\mbox{co}}(D-z_0)$ has at least  an  interior point
$w_0$ in the subspace $\overline{\mbox{span}}\{D-z_0\}$.
Hence,  for some $\sigma>0$,  it holds that
$$
B_X(w_0, \sigma)\cap\overline{\mbox{span}}\{D-z_0\}\subseteq
\overline{\mbox{co}}(D-z_0).
$$
This implies that
\begin{equation}\label{appendix1}
B_X(0, \sigma)\cap\big(\overline{\mbox{span}}\{D-z_0\}-w_0\big)\subseteq
\overline{\mbox{co}}(D-z_0)-w_0.
\end{equation}

On the other hand,  for any $z\in \overline{\mbox{span}}\{D-z_0-w_0\}$,
there exists a sequence $\{z_j\}_{j=1}^\infty\subseteq \mbox{span}\{D-z_0-w_0\}$,
such that
$
z_j\rightarrow z \mbox{ in  }X, \mbox{ as } j\rightarrow\infty.
$
For any $j\in\dbN$,  there exists a positive integer $k\in\dbN$,
$w^\ell_j\in D-z_0$ and $\alpha^\ell_j\in\dbR$ $(\ell=1, 2, \cdots, k)$,
such that
$$
z_j=\sum\limits_{\ell=1}^{k}\alpha^\ell_j(w^\ell_j-w_0)
=\sum\limits_{\ell=1}^{k}\alpha^\ell_j w^\ell_j+
\big(1-\sum\limits_{\ell=1}^{k}\alpha^\ell_j\big)w_0-w_0
\in \overline{\mbox{span}}\{D-z_0\}-w_0.
$$
Then, $z\in \overline{\mbox{span}}\{D-z_0\}-w_0$, and therefore,
$$
\overline{\mbox{span}}\{D-z_0-w_0\}\subseteq  \overline{\mbox{span}}\{D-z_0\}-w_0.
$$
This, together with (\ref{appendix1}),  means that
\begin{equation}\label{appendix2}
B_X(0, \sigma)\cap\overline{\mbox{span}}\{D-z_0-w_0\}\subseteq
\overline{\mbox{co}}(D-z_0)-w_0.
\end{equation}

\medskip

{\bf  Step 2. }
By  the finite codimensionality of $D$,
we have that
$\overline{\mbox{span}}\{D-z_0\}$ is a finite codimensional
subspace of $X$. Then,  by Proposition  \ref{add**},  the subspace
$\overline{\mbox{span}}\{D-z_0-w_0\}$  is also finite  codimensional, since $z_0+w_0\in
 \overline{\mbox{co}}D$.

Hence,  there exists a $\nu\in\dbN$ and  linearly independent
$\tilde x_1, \cdots, \tilde x_\nu\in X
\setminus\overline{\mbox{span}}\{D-z_0-w_0\}$,  such that
$$
X=\overline{\mbox{span}}
\{D-z_0-w_0\}+\span\{\tilde x_1, \cdots, \tilde x_\nu\}.
$$
Set
$$
L_j=\overline{\mbox{span}}
\{D-z_0-w_0\}+\span\{\tilde x_1, \cdots, \tilde x_{j-1},  \tilde x_{j+1},
\cdots,  \tilde x_\nu\}, \quad 1\leq j\leq \nu.
$$
Then,  $\tilde x_j\notin L_j$ and by the Hahn-Banach theorem,  there
exist $f_1, \cdots, f_\nu\in X'$, such that
$$
 f_j(\tilde x_j)=1\  \mbox{ and }\
L_j\subseteq\mathcal{N}(f_j)
=\Big\{
z\in X\  \Big|\
f_j(z)=0\Big\}, \quad \forall\; j=1, \cdots, \nu.
$$
It follows that
$
\overline{\mbox{span}}
\{D-z_0-w_0\}\subseteq \bigcap\limits_{j=1}^{\nu}\mathcal{N}(f_j).
$

On the other hand, for any $x\in  \bigcap\limits_{j=1}^{\nu}\mathcal{N}(f_j),$
there exists an  $x_0\in \overline{\mbox{span}}
\{D-z_0-w_0\}$ and $\beta_1, \cdots, \beta_\nu\in\dbR$, such that
$
x=x_0+\sum\limits_{j=1}^\nu \beta_j \tilde x_j.
$
Then for any $j=1, \cdots, \nu$,
$$
0=f_j(x)=f_j(x_0)+\sum\limits_{i=1}^{\nu}\beta_i f_j(\tilde x_i)=\beta_j.
$$
Hence, $x=x_0\in \overline{\mbox{span}}
\{D-z_0-w_0\}$, which  implies
$
\overline{\mbox{span}}
\{D-z_0-w_0\}=\bigcap\limits_{j=1}^{\nu}\mathcal{N}(f_j).
$

\medskip

{\bf Step 3. } Set   $\kappa=\max\big\{|f_{j}|_{X'}\  \big|\ j=1,\cdots,\nu\big\}$ and
$$
K=\Big\{\sum\limits_{i=1}^{\nu} \alpha_i\tilde x_i\in X\  \Big|\
|\alpha_i|\leq \kappa ,  i=1, \cdots, \nu  \Big\}.
$$
Then $K$ is  a compact set in $X$.

Let $x\in B_X(0, 1)$.
For any $j=1, \cdots, \nu$, by the definition of $f_j$,
$f_j(x-\sum\limits_{i=1}^{\nu} f_i(x)\tilde x_i)=0$. Hence,
$x-\sum\limits_{i=1}^{\nu} f_i(x)\tilde x_i\in
\bigcap\limits_{j=1}^{\nu}\mathcal{N}(f_j)=\overline{\mbox{span}}
\{D-z_0-w_0\}$.  In addition, since $|f_j(x)|\leq |f_j|_{X'}|x|_X\leq \kappa$, we have $ \sum\limits_{i=1}^{\nu} f_i(x)\tilde x_i\in K$. Therefore,
$$
x=\Big[x-\sum\limits_{i=1}^{\nu} f_i(x)\tilde x_i\Big]+
\sum\limits_{i=1}^{\nu} f_i(x)\tilde x_i\in \overline{\mbox{span}}
\{D-z_0-w_0\}+K.
$$

Note that $K$ is a compact set. There is  a $\gamma>0$, such that
$$\Big|x-\sum\limits_{i=1}^{\nu} f_i(x)\tilde x_i\Big|_{X}\le
|x|_{X}+\Big|\sum\limits_{i=1}^{\nu} f_i(x)\tilde x_i\Big|_{X}\le
 1+\gamma.$$
It implies that
$
B_X(0, 1)\subseteq\overline{\mbox{span}}
\{D-z_0-w_0\}\cap B_X(0, 1+\gamma)+K.
$
Letting $0<\rho\le\min\{\frac{\sigma}{ 1+\gamma},1\}$, we have that
$$
B_X(0, \rho)\subseteq\overline{\mbox{span}}
\{D-z_0-w_0\}\cap B_X(0, \sigma)+K.
$$
By  (\ref{appendix2}),
$$
B_X(0, \rho)\subseteq\overline{\mbox{co}}(D-z_0)-w_0+
K=\overline{\mbox{co}}D-z_0-w_0+K.
$$

{\bf Step 4. } For any $y\in B_X(0, \rho)$,
$$
y=x+z, \quad\mbox{with }x\in \overline{\mbox{co}}D-z_0-w_0
\mbox{ and }z\in K.$$
By the assumptions on $\Lambda_k$,
$\left\langle \Lambda_{k}, x\right\rangle_{X', X}
\geq-\varepsilon_{k}$ for all $(x,k) \in D\times \dbN$, which implies  $\left\langle \Lambda_{k}, x\right\rangle_{X', X}
\geq-\varepsilon_{k}$ for all $(x,k) \in \overline{\mbox{co}}D\times \dbN$. Therefore,
$$
\langle \Lambda_k, y\rangle_{X', X}=
\langle \Lambda_k, x\rangle_{X', X}+\langle \Lambda_k, z\rangle_{X', X}
\geq -\varepsilon_k-
\langle \Lambda_k,  z_0+w_0\rangle_{X', X}+
\langle \Lambda_k, z\rangle_{X', X}.
$$
Since $-y\in B_X(0, \rho)$, it holds that
\begin{eqnarray*}
|\Lambda_k|_{X'}
&\leq& \frac{1}{\rho}\varepsilon_k+
\frac{1}{\rho}\Big|\langle \Lambda_k, z_0+w_0\rangle_{X', X}\Big|
+\frac{1}{\rho}\sup\limits_{z\in K}
\Big|\langle \Lambda_k, z\rangle_{X', X}\Big|\\
&\leq& \frac{1}{\rho}\varepsilon_k+
\frac{1}{\rho}\Big|\langle \Lambda_k, z_0+w_0\rangle_{X', X}\Big|
+ \frac{\kappa}{\rho}\sum\limits_{i=1}^\nu
\Big|\langle \Lambda_k, \tilde x_i\rangle_{X', X}\Big| .
\end{eqnarray*}
If $\Lambda=0$,   then
$
\lim\limits_{k\rightarrow\infty}|\Lambda_k|_{X'}=0$.
This contradicts   the fact that $|\Lambda_k|_{X'}\geq \delta>0.$
\endpf

\section{Proof of Case 2 in Example  \ref{remark 2.4}}\label{D}

In this appendix, we provide a proof for {\bf Case 2} in Example  \ref{remark 2.4}.

\medskip

{\bf Step 1. }  First, if $g(\bar{u})<0$, then by  the continuity of $g$, there exists  $r_0>0$ such that
$g(u)<0$ for any $u\in B_V(\bar{u}, r_0)$.  Hence, $\bar{u}, u\in $ Int $\mathcal{K}$, and
$$
\cT_{\mathcal{K}}(\bar u)=\cT_{\mathcal{K}}(u)=V,\quad \mbox{ for all }
u\in B_V(\bar{u}, r_0)\cap\mathcal{K}.
$$
Thus,  the assertion (\ref{2.2-eq1}) holds.

\medskip

{\bf Step 2. }  Next, assume that $g(\bar{u})=0$.
Due to the continuous differentiability of $g$ and the assumption $g'(\bar{u})\neq 0$,
we choose $r_0$ sufficiently small so that
$$g'(w)\neq 0 \quad \mbox{ for
any }  w\in B_V(\bar{u}, r_0)\cap\mathcal{K} \mbox{ satisfying  } g(w)=0.$$
Then, we will prove that  for the aforementioned $w$,
\begin{equation}\label{appendixd1}
\cT_{\mathcal{K}}(w)=\Big\{ v\in V\ \Big|\  \langle g'(w), v\rangle_{V', V}\leq 0 \Big\}.
\end{equation}
Indeed, for any $v\in \mathcal{R}_{\mathcal{K}}(w)$,
there exist $\alpha\geq 0$ and $\nu \in \mathcal{K}$ such that
$v=\alpha(\nu-w)$. Since $g$ is continuously differentiable and convex,  $g(w)=0$, and $g(\nu)\leq 0$, we obtain
$\langle g'(w), \nu-w\rangle_{V', V}\leq g(\nu)-g(w)\leq 0$.
Hence,
$$
\langle g'(w), v\rangle_{V', V}=\alpha \langle g'(w), \nu-w \rangle_{V', V}\leq 0,
$$
which implies
$$
\mathcal{R}_{\mathcal{K}}(w)\subseteq \Big\{
v\in V\ \Big|\  \langle g'(w), v\rangle_{V', V}\leq 0
\Big\}.
$$
Since the set on the right-hand side is closed, it follows that
$$\cT_{\mathcal{K}}(w)\subseteq\Big\{ v\in V\ \Big|\  \langle g'(w), v\rangle_{V', V}\leq 0 \Big\}.$$

Next, for any $v\in V$ satisfying $\langle g'(w), v\rangle_{V', V}< 0$,
the  differentiability of $g$ yields
$g(w+\tau v)=g(w)+\tau \langle g'(w), v\rangle_{V', V}+o(\tau)
=\tau \langle g'(w), v\rangle_{V', V}+o(\tau)$, where $o(\tau)$ denotes   a higher-order infinitesimal with respect to $\tau$.
Therefore, $g(w+\tau v)<0$ for a sufficiently small $\tau>0$.
This indicates $w+\tau v\in \mathcal{K}$ and consequently,
$$
v=\displaystyle\frac{(w+\tau v)-w}{\tau}\in \mathcal{R}_{\mathcal{K}}(w)\subseteq \mathcal{T}_{\mathcal{K}}(w).
$$

Now,  suppose that $v\in V$ satisfies  $\langle g'(w), v\rangle_{V', V}= 0$.
For any $\epsilon>0$, define $v_\epsilon=v-\epsilon \tilde{G}$, with $\tilde{G}\in V$ being
 the element corresponding
to $g'(w) \in V'$ via the Riesz-Fr\'echet
isomorphism.
Then,  $$
\langle g'(w), v_\epsilon\rangle_{V', V}=-\epsilon \langle g'(w), \tilde{G}\rangle_{V', V}=-\epsilon
|g'(w)|^2_{V'}<0.
$$
By the preceding argument, $v_\epsilon\in \mathcal{T}_{\mathcal{K}}(w)$.
Since $v_\epsilon\rightarrow v$ in $V$, we obtain $v\in \mathcal{T}_{\mathcal{K}}(w)$. This implies that
$$
\Big\{ v\in V\ \Big|\  \langle g'(w), v\rangle_{V', V}\leq 0 \Big\}\subseteq\cT_{\mathcal{K}}(w).
$$
Therefore, (\ref{appendixd1})
holds.

\medskip

{\bf Step 3. } Let $u\in B_V(\bar{u}, r_0)\cap\mathcal{K}$. If
$g(u)<0$, then $\mathcal{T}_{\mathcal{K}}(u)=V$,  and the desired inclusion (\ref{2.2-eq1}) is trivial.

Assume now that $g(u)=0$. Then, $g'(u)\neq 0$. We define
$$
\displaystyle n_1=\frac{g'(\bar{u})}{|g'(\bar{u})|_{V'}}
\quad\mbox{and}\quad
n_2=\frac{g'(u)}{|g'(u)|_{V'}}.
$$
It holds that $|n_1|_{V'}=|n_2|_{V'}=1$ and
(\ref{appendixd1}) gives
$$
\cT_{\mathcal{K}}(\bar u)=\Big\{ v\in V\ \Big|\  \langle n_1, v\rangle_{V', V}\leq 0 \Big\}\quad\mbox{ and }\quad
\cT_{\mathcal{K}}(u)=\Big\{ v\in V\ \Big|\  \langle n_2, v\rangle_{V', V}\leq 0 \Big\}.
$$
For any $v\in \cT_{\mathcal{K}}(\bar u)\cap B_V(0, 1)$,  we have
$\langle n_1, v\rangle_{V', V}\leq 0$ and $|v|_V\leq 1$.
Define
$$\hat{\alpha}=\max\{0, \langle n_2, v\rangle_{V', V}\}\quad\mbox{and}\quad
\hat{w}=v-\hat{\alpha} \hat{n}_2,$$
where $\hat{n}_2\in V$ is the element corresponding
to $n_2 \in V'$ via the Riesz-Fr\'echet
isomorphism.
Note that $ \langle n_2, \hat{w}\rangle_{V', V}=\langle n_2, v\rangle_{V', V}-\hat{\alpha}\leq 0.$
Thus,
$\hat{w}\in \mathcal{T}_{\mathcal{K}}(u)$.

Next, we will show that $\hat{w}\in B_V(0, 1)$. Indeed,
if $\langle n_2, v\rangle_{V', V}\leq 0$,
then $\hat{\alpha}=0$. It follows that $\hat{w}=v$ and hence, $|\hat{w}|_V\leq 1$.
If $\langle n_2, v\rangle_{V', V}> 0$,
then $\hat{\alpha}=\langle n_2, v\rangle_{V', V}$. This indicates
$$
|\hat{w}|_V^2=|v-\hat{\alpha} \hat{n}_2|^2_V
=|v|^2_V-2\hat{\alpha} \langle n_2, v\rangle_{V', V}+\hat{\alpha}^2
=|v|^2_V-\hat{\alpha}^2\leq 1.
$$
Therefore, $\hat{w}\in  \mathcal{T}_{\mathcal{K}}(u)\cap B_V(0, 1).$

\medskip

{\bf Step 4. }  For any $v\in \cT_{\mathcal{K}}(\bar u)\cap B_V(0, 1)$, since $\langle n_1, v\rangle_{V', V}\leq 0$,
we have
$$
\langle n_2, v\rangle_{V', V}\leq \langle n_2-n_1, v\rangle_{V', V}\leq |n_2-n_1|_{V'}.
$$
Note that $|v-\hat{w}|_V=\hat{\alpha}$.
Then,
$$
|v-\hat{w}|_V=\max\{0, \langle n_2, v\rangle_{V', V}\}\leq |n_2-n_1|_{V'}.
$$
This demonstrates  that
\begin{equation}\label{appendix22}
\mathcal{T}_{\mathcal{K}}(\bar{u})\cap B_V(0, 1)
\subseteq \mathcal{T}_{\mathcal{K}}(u)\cap B_V(0, 1)+|n_2-n_1|_{V'} B_{V}(0, 1).
\end{equation}

Finally, we estimate  $|n_2-n_1|_{V'}$.  Recall that
$\displaystyle n_1=\frac{g'(\bar{u})}{|g'(\bar{u})|_{V'}}
\quad\mbox{and}\quad
n_2=\frac{g'(u)}{|g'(u)|_{V'}}.$
Then,
\begin{eqnarray*}
&&n_1-n_2=\frac{g'(\bar{u})}{|g'(\bar{u})|_{V'}}-\frac{g'(u)}{|g'(u)|_{V'}}
=\frac{g'(\bar{u})-g'(u)}{|g'(\bar{u})|_{V'}}+g'(u)\Big(\frac{1}{|g'(\bar{u})|_{V'}}-\frac{1}{|g'(u)|_{V'}}\Big).
\end{eqnarray*}
Thus,
\begin{eqnarray*}
&&|n_1-n_2|_{V'}\leq
\frac{|g'(\bar{u})-g'(u)|_{V'}}{|g'(\bar{u})|_{V'}}+
\frac{\Big||g'(u)|_{V'}-|g'(\bar{u})|_{V'}\Big|}{|g'(\bar{u})|_{V'}}
\leq 2\frac{|g'(\bar{u})-g'(u)|_{V'}}{|g'(\bar{u})|_{V'}}.
\end{eqnarray*}
By the assumption on $g$,
$$
|n_1-n_2|_{V'}\leq 2\frac{\ell(|\bar{u}-u|_V)}{|g'(\bar{u})|_{V'}}.
$$
This, together with (\ref{appendix22}), indicates
$$\displaystyle
\mathcal{T}_{\mathcal{K}}(\bar{u})\cap B_V(0, 1)
\subseteq \mathcal{T}_{\mathcal{K}}(u)\cap B_V(0, 1)+2\frac{\ell(|\bar{u}-u|_V)}{|g'(\bar{u})|_{V'}} B_{V}(0, 1).
$$
This finishes the proof of
{\bf Case 2} in Example  \ref{remark 2.4}.  \endpf

\end{document}